\documentclass[11pt]{amsart}
\usepackage{blindtext}
\usepackage{url}
\usepackage{tabularx}
\usepackage[utf8]{inputenc}
\usepackage[margin=2.5cm]{geometry}
\usepackage[utf8]{inputenc}
\usepackage{amsmath}
\usepackage{amsfonts}
\usepackage{color}
\usepackage{slashbox}
\usepackage{amssymb}
\usepackage{amsthm}
\usepackage{array}
\usepackage{rotating}
\usepackage{physics}
\usepackage{tikz}
\usepackage{braids}
\usepackage{tikz-cd}
\tikzset{%
	symbol/.style={%
		draw=none,
		every to/.append style={%
			edge node={node [sloped, allow upside down, auto=false]{$#1$}}}
	}
}
\usepackage[all]{xy}
\usepackage{cite}
\usepackage{mathdots}
\usepackage{pstricks}
\usepackage{mathrsfs}
\usepackage{setspace}
\usepackage{hyperref}
\theoremstyle{theorem}
\newtheorem{theorem}{Theorem}

\numberwithin{theorem}{subsection}
\newtheorem{lemma}[theorem]{Lemma}
\newtheorem{proposition}[theorem]{Proposition}
\newtheorem{corollary}[theorem]{Corollary}
\theoremstyle{definition}
\newtheorem{definition}[theorem]{Definition}
\newtheorem{example}[theorem]{Example}
\theoremstyle{remark}
\newtheorem{remark}[theorem]{Remark}
\theoremstyle{notation}
\newtheorem{notation}[theorem]{Notation}
\newcommand{\Gr}{\mathbf{Gray}}
\newcommand{\A}{\mathcal{A}}
\newcommand{\B}{\mathcal{B}}

\newcommand{\PsG}{\mathbf{Psmnd}(\mathbf{Gray})}

\newcommand{\Cat}{\mathbf{Cat}}
\newcommand{\Past}{{\mathbf{Psmnd}\left(\mathbf{Gray}\right)}_\text{ast}}
\title{Eilenberg-Moore Bicategories for Opmonoidal Pseudomonads}
\author{Adrian Miranda}
\thanks{The material is based on work supported by EPSRC under grant EP/V002325/2. I would also like to thank Nicola Gambino for helpful discussions while I was conducting this research and for advice while I was preparing this paper.}
\address{Department of Mathematics, University of Manchester, United Kingdom}
\email{adrian.miranda@manchester.ac.uk}

\begin{document}
	
	\begin{abstract}
		\noindent We analyse compatibility between monads and monoidal structures in the two-dimensional setting. We describe sufficient conditions for monoidal structures to lift to the Eilenberg-Moore pseudoalgebras. We then extend these results to braids, syllapses and symmetries. To achieve these results we define the $\mathbf{Gray}$-tensor product of pseudomonads, and examine its interaction with the Eilenberg-Moore construction.
	\end{abstract}
	
	\maketitle

	\tableofcontents

\section{Introduction}\label{Introduction}

\subsection{Context and goals}
\noindent Symmetric monoidal bicategories have applications to combinatorics \cite{Fiore Gambina Hyland Winskel generalised species}, two-dimensional models of differential linear logic \cite{Fiore Gambino Hyland Monoidal Bicategories differential linear logic and analytic functors, From Thin Concurrent Games to Generalized Species of Structures, Galal Bicategorical Finite Nondeterminism}, and three-dimensional topological quantum field theories \cite{Schommer Piers Classification of 2D ETFTs, BDSV Extended 3D bordism theory of modular objects, BDSV Modular categories as representations of the 3-dimensional bordism 2-category}. They are however intricate structures, typically difficult to even construct in practice \cite{Hansen Shulman Constructing Symmetric Monoidal Bicategories Functorially}. Recent research has focused on simplifying this complexity via coherence results \cite{Crans Generalised Centers, Gurski Loop Spaces, Gurski Osorno Infinite Loop Spaces and coherence for symmetric monoidal bicategories} and studying structures which naturally live in a symmetric monoidal bicategory \cite{monoidal bicategories and hopf algebroids}. 
\\
\\
 \noindent Many complicated mathematical structures can be described as algebras for monads over categories of simpler structures. Given a monad whose underlying category has some extra structure, it is of interest to describe how that extra structure might be compatible with the monad so that it should lift to the category of algebras  \cite{Bruguieres Lack Virelizier, Hasegawa Lemay Star Autonomous Monads, Traced Monads and Hopf Monads, Moerdijk Monads on Tensor Categories, Pastro Street closed categories star autonomy and monoidal comonads}. In this paper we pursue a study of such compatibility in the two-dimensional setting, where this extra structure is that of a symmetric monoidal bicategory. As such our results provide a streamlined and simpler method for constructing new symmetric monoidal structures for pseudoalgebras, from known structures on bases of pseudomonads. 
\\
\\
\noindent We introduce symmetric opmonoidal pseudomonads in Definitions \ref{definition semi strict opmonoidal pseudomonads} and \ref{Definition braided opmonoidal pseudomonad} to precisely capture the appropriate compatibility between the two structures. These categorify the symmetric opmonoidal monads of \cite{McCrudden Opmonoidal Monads} whose natural families of maps $\chi_{X, Y}: S(X\otimes Y) \to SX \otimes SY$ and $\iota: SI \to I$ mediate interaction between the monad and symmetric monoidal structure, and satisfy axioms which ensure that these structures lift to categories of algebras. This notion of compatibility differs from others which have been considered in the two-dimensional setting such as strengths \cite{Hyland Power Pseudo commutative monads and pseudo closed categories, Pacquet Saville Strong Pseudomonads and Premonoidal Bicategories} which instead relate to extensions of structures to Kleisli bicategories.

\subsection{Summary of main results and contributions}

\noindent We prove a coherence result in Theorem \ref{theorem strictifying strong monoidal pseudomonads} to simplify our study from symmetric monoidal bicategories to symmetric $\mathbf{Gray}$-monoids without losing generality. These satisfy associativity and left and right unit laws on the nose, but have tensor products for which joint functoriality remains weak. Perhaps surprisingly, the lifted monoidal structure on $2$-categories of pseudoalgebras that we establish in Theorem \ref{monoidal bicategory structure on pseudoalgebras} turns out to be slightly weaker than this. It only satisfies associativity and left and right unit laws up to $2$-natural isomorphisms, with these data obeying the usual axioms for a monoidal category on the nose. This extra weakness comes from data in opmonoidal $2$-functors which mediate the usual axioms of opmonoidal functors.
\\
\\
\noindent One important stepping stone towards this result is Theorem \ref{Theorem Gray tensor of pseudomonads}, in which we describe a $\Gr$-tensor product of pseudomonads. Another important stepping stone is Theorem \ref{theorem semi-strict opmonoidal pseudomonads are pseudomonoids in the Gray monoid of pseudomonads} in which we describe a bijection between pseudomonoids in the resulting $\mathbf{Gray}$-monoid of pseudomonads, and semi-strict opmonoidal pseudomonads. Section \ref{Section Interaction between the Eilenberg-Moore construction and Gray tensor products} studies the interaction between the Eilenberg-Moore construction for pseudomonads and $\mathbf{Gray}$-tensor products. If $2$-natural transformations were closed under $\mathbf{Gray}$-tensor product, these results could be summarised as saying that the Eilenberg-Moore construction is a normal lax monoidal $2$-functor. This intuition is made precise in Theorem \ref{Theorem EM is normal lax Gray functor}, after a suitable closure of $2$-natural transformations under the $\Gr$ tensor product is described in Definition \ref{Definition Gray 2-nat}.
\\
\\
\noindent Finally, in Section \ref{Section symmetries for opmonoidal pseudomonads} we describe compatibility between pseudomonads and further structure such as braidings, syllapses and symmetry. Proposition \ref{Proposition Braidings Syllapses and Symmetries as data in Psmnd Gray} shows that such structures also give rise to higher dimensional cells of pseudomonads, while Theorem \ref{Theorem lifted braidings syllapses and symmetries} shows that they also give rise to braidings, syllapses and symmetries at the level of $2$-categories of pseudoalgebras.

\subsection{Key ideas and techniques}

\subsubsection{Strictification}

\noindent One of the major challenges in lifting monoidal, braid, syllepsis and symmetry structures to pseudoalgebras for pseudomonads is the plethora of complex data and axioms involved. In Subsection \ref{Subsection opmonoidal pseudomonads as pseudomonads in the tricategory of monoidal bicategories} we give a systematic description of this information using the framework of tricategories and their various higher dimensional transfors developed in \cite{Low dimensional structures formed by tricategories}. In Theorem \ref{theorem strictifying strong monoidal pseudomonads} we justify restricting our attention to simpler structures via which general results can still be established. This exploits the strictification results of tricategories and trihomomorphisms \cite{Gurski PhD}, as well as analogous results for higher dimensional maps between these \cite{Miranda strictifying operational coherences}.

\subsubsection{Alternative perspectives on compatibility between structures}

\noindent Even after restricting our attention to the semi-strict setting, a complex interplay between different structures remains. We reorganise these data into a form convenient for our goals in Subsection \ref{Section monoidal pseudomonads as pseudomonoids in pseudomonads}. In particular, the associativity and left and right unit constraints of the underlying opmonoidal $2$-functor give rise to a pseudomonoid structure on the pseudomonad. Viewing semi-strict opmonoidal pseudomonads as pseudomonoids in an appropriate $\mathbf{GRAY}$-monoid of pseudomonads conceptually motivates our strategy towards lifting monoidal structures to pseudoalgebras. 

\subsubsection{Ambistrict $2$-cells}

\noindent The $\mathbf{Gray}$-tensor product is not a limit, and does not automatically commute with Eilenberg-Moore constructions for pseudomonads just because the latter are limit constructions. A direct calculation of the interaction between $\mathbf{Gray}$-tensor products and Eilenberg-Moore constructions is necessary. We find that the associators and left and right unitors of semi-strict opmonoidal pseudomonads have the special property of satisfying middle-four interchange with arbitrary pseudomonad transformations on either side. These `ambistrict' $2$-cells are sent to $2$-natural transformations between $2$-categories of pseudoalgebras, and moreover the components of these $2$-natural transformations are also of pseudomorphisms of a particularly simple form. 
\subsubsection{$\mathbf{Gray}$-tensor products of $2$-natural transformations}

\noindent Even after simplifying and reorganising the data in opmonoidal pseudomonads, some complexity remains. Conjecturally, there should be a three-dimensional monoidal structure on pseudomonads and the Eilenberg-Moore construction $\mathbf{EM}: \PsG \to \Gr^\mathbf{2}$ should be lax monoidal with respect to this structure, with the $\mathbf{Gray}$-tensor product on the target. Not wanting to get bogged down in developing theory for three-dimensional monoidal categories, we seek a two-dimensional monoidal structure which receives the Eilenberg-Moore construction of pseudomonads. Although the $2$-cells we need are $2$-natural transformations, these fail to be closed under the $\mathbf{Gray}$ tensor product. Hence describing an appropriate two-dimensional monoidal structure requires quotienting certain pseudonatural transformations under a congruence generated by middle-four interchange. As such, even after we make precise the sense in which the Eilenberg-Moore construction $(\A, S) \mapsto \A^S$ for pseudomonads is lax monoidal, further calculation is needed to establish the lifted monoidal structure on $\A^S$.

\subsubsection{Faithfulness on $2$-cells of $U^S: \A^S \to \A$}

\noindent When we move to braidings, syllapses and symmetries we find that compatibility between pseudomonads and these structures gives rise to various cells of pseudomonads. These include arbitrary pseudomonad transformations, rather than just ambistrict ones, and they also include pseudomonad modifications. Seeing braided (resp. sylleptic, symmetric) opmonoidal pseudomonads as $E_{n}$-algebras in an appropriate symmetric monoidal $\mathbf{Gray}$-category of pseudomonads remains open. Nonetheless, each individual cell of pseudomonads is sent by $\mathbf{EM}$ to appropriate data on $2$-categories of pseudoalgebras, using the results of \cite{Formal Theory of Pseudomonads}. We are then left to check that these data on the $2$-category of pseudoalgebras satisfy the appropriate axioms for braiding, syllepsis or symmetry. But we find that each of these axioms follows immediately from the analogous axiom for the original braiding, syllepsis or symmetry structure on the base $2$-category of the pseudomonad. This uses faithfulness on $2$-cells of the $2$-functor $\A^S \to \A$ which forgets pseudoalgebra and pseudomorphism structure.

\subsection{Conventions}

\subsubsection{Diagrams and notation}

\noindent Some of the observations and proofs in this paper follow immediately after the relevant definitions are unwound, requiring no further intermediate steps. Others require detailed calculations involving large pasting diagrams. When calculations are required, the steps involved are typically just interchange axioms for $\mathbf{Gray}$-monoids, as recalled in Remark \ref{Remark Gray tensor product}. A notable exception to this pattern in in Subsection \ref{Subsection the Gray tensor product of pseudomonads}, where proofs also involve the axioms for pseudomonads or for cells between them. For space reasons, most of the calculations are explained in words; detailed pasting diagram chases are either omitted or deferred to the appendix. In such proofs, colour is used to distinguish new data appearing in each step. We use \textcolor{blue}{blue} for new objects and morphisms, and \textcolor{red}{red} for new $2$-cells. Although we will often be able to restrict our attention to structures that are strictly associative and unital, in Remark \ref{Remark explicit opmonoidal pseudomonad} we need to adopt the convention that composites in bicategories are to be read as parenthesised to the left $(hg)f$ unless otherwise specified. To fit diagrams onto the page without reducing the font size too much, we sometimes need to adopt ad-hoc notational conventions such as writing ordered pairs as columns rather than rows. Such conventions will be explained whenever they are used.

\subsubsection{Size}

\noindent While size is not a central issue considered in this paper. we do need to carefully distinguish between mathematical structures of different sizes. Specifically, we need to consider a two-dimensional monoidal structure whose objects are also two-dimensional monoidal structures. To avoid set theoretical paradoxes, we fix two universes of sets $\mathfrak{U}_{1}$ and $\mathfrak{U}_{2}$, with $\mathbb{N} \in \mathfrak{U}_{1} \in \mathfrak{U}_{2}$. We call the sets in $\mathfrak{U}_{1}$ \emph{small}, and the sets in $\mathfrak{U}_{2}$ \emph{large}. A small category is one whose set of objects and set of morphisms are both in $\mathfrak{U}_{1}$, while a large category is one whose set of objects and set of morphisms are in $\mathfrak{U}_{2}$. A small $2$-category is one whose set of objects is in $\mathfrak{U}_{1}$, and whose hom-categories are small, while a large $2$-category is one whose set of objects is in $\mathfrak{U}_{2}$ and whose hom-categories are large. Finally, $\mathbf{Gray}$ denotes the three dimensional structure consisting of small $2$-categories, $2$-functors, pseudonatural transformations, and modifications, while $\mathbf{GRAY}$ denotes the three dimensional structure whose objects are large $2$-categories and whose higher cells are analogous to those in $\mathbf{Gray}$.

\section{Review of symmetric opmonoidal monads}\label{Subsection review of symmetric opmonoidal monads}

\noindent Denote by $\mathbf{MonCat}_\text{oplax}$ the $2$-category of monoidal categories, opmonoidal functors and monoidal transformations, and the full-sub-$2$-category on strict monoidal categories as $\mathbf{strMonCat}_\text{oplax}$. The wide, locally full sub-$2$-categories spanned by strong or strict monoidal functors are denoted with the corresponding subscripts in place of $\text{oplax}$.

\begin{remark}\label{remark opmonoidal monads as monads in MonCat}
	Opmonoidal monads are the objects of a $2$-category $\mathbf{Mnd}\left(\mathcal{K}\right)$, described for arbitrary $2$-categories in \cite{Formal Theory of Monads}, with $\mathcal{K} = \mathbf{MonCat}_\text{oplax}$. More explicitly, an opmonoidal monad consists of a monoidal category $\left(\mathcal{V}, \bigoplus, I, \alpha, \lambda, \rho\right)$, an endofunctor $T: \mathcal{V} \rightarrow \mathcal{V}$, and natural transformations with components $\iota: TI \rightarrow I$, $ \chi_{{}_{X, Y}}: T\left(X \bigoplus Y\right) \rightarrow TX \bigoplus TY$, $\eta_{X}: X\rightarrow TX$ and $\mu_{X}: T^{2}X \rightarrow TX$, subject to the conditions that the data $(T, \chi, \iota): (\mathcal{V}, \bigoplus, I, \lambda, \alpha, \rho) \to (\mathcal{V}, \bigoplus, I, \lambda, \alpha, \rho)$ comprise a opmonoidal functor, the natural transformations $\eta: 1_{V} \Rightarrow T$ and $\mu: T^2 \Rightarrow T$ are monoidal, and the data $(T, \eta, \mu)$ comprise a monad on the category $\mathcal{V}$.
\end{remark}

\noindent The perspective on opmonoidal monads described in Remark \ref{remark opmonoidal monads as monads in MonCat} motivates the definition of opmonoidal pseudomonads that we will give in Definition \ref{definition colax monoidal pseudomonad}. Opmonoidal monads have a coherence result, which says that they are equivalent in $\mathbf{Mnd}(\mathbf{MonCat_\text{oplax}})$ to one whose underlying monoidal category is strict. Moreover, if the underlying opmonoidal functor is strong, then the opmonoidal monad is moreover equivalent to another opmonoidal monad whose underlying endofunctor is also strict. The proof is via standard strictification techniques analogous to the proof we will give for Theorem \ref{theorem strictifying strong monoidal pseudomonads}.

\begin{remark}\label{Opmonoidal monads as pseudomonoids in Mnd Cat}
	There is another perspective on opmonoidal monads. The $2$-category $\mathbf{Mnd}\left(\mathbf{Cat}\right)$ has products given by taking products of each piece of data in $\mathbf{Cat}$. As observed in Example 3.7 of \cite{McCrudden Opmonoidal Monads}, pseudomonoids in the resulting monoidal $2$-category are precisely opmonoidal monads while strict monoids in $(\mathbf{Mnd}(\mathbf{Cat}), \times, \mathbf{1})$ are precisely opmonoidal monads on strict monoidal categories. In the two-dimensional setting, associativity and left and right unit laws only hold weakly for opmonoidal $2$-functors. As such, semi-strict opmonoidal pseudomonads only correspond to pseudomonoids in an appropriate $\mathbf{Gray}$-monoid of pseudomonads, to be described in Subsections \ref{Subsection the Gray tensor product of pseudomonads} and \ref{Section monoidal pseudomonads as pseudomonoids in pseudomonads}.
\end{remark} 

\noindent The perspective on opmonoidal monads outlined in Remark \ref{Opmonoidal monads as pseudomonoids in Mnd Cat} streamlines the lifting of monoidal structures to categories of algebras. Consider the $2$-functor $\mathbf{EM}: \mathbf{Mnd}\left(\mathbf{Cat}\right) \hookrightarrow \mathbf{Cat}^\mathbf{2}$, which acts on objects by sending a monad to the forgetful functor from its Eilenberg-Moore category of algebras. This $2$-functor $\mathbf{EM}$ is a strictly unital strong monoidal $2$-functor. Since strong monoidal $2$-functors preserve pseudomonoids, an opmonoidal pseudomonad has a monoidal structure on its category of algebras. Moreover, the monadic functor $U^S: A^S \rightarrow A$ is strict monoidal, so if the base is a strict monoidal category then the lifted monoidal structure on $A^S$ is also strict.
	\\
	\\
\noindent One way to see that Eilenberg-Moore objects for monads commute with cartesian products is to observe that both constructions are examples of limits. In the two-dimensional setting, the $\mathbf{Gray}$-tensor product of pseudomonads is no longer a limit. More delicate arguments relating the constructions are given in Subsection \ref{Subsection normal lax monoidality of EM}. The comparisons $\A^S \otimes \B^T \to {\left(\A\otimes \B\right)}^{(S\otimes T)}$ turn out to not even be equivalences, and as such the Eilenberg-Moore construction for pseudomonad turns out to only be normal lax monoidal.
\\
\\
\noindent Similar arguments extend to braided and symmetric opmonoidal monads. These correspond to braidings and symmetries in the sense of \cite{monoidal bicategories and hopf algebroids} on the corresponding pseudomonoids in the monoidal $2$-category $(\mathbf{Mnd}(\Cat), \times, \mathbf{1})$. In particular, the condition for $\beta$ to be well-defined as a $2$-cell of monads is equally the condition for $(S, \chi, \iota)$ to be braided as an opmonoidal functor. In the two-dimensional setting, there should also be a correspondence between braidings, syllapses and symmetries for opmonoidal pseudomonads, and for similar structure for pseudomonoids in the $\Gr$-category of pseudomonads $\PsG$. We establish enough of this correspondence, and enough compatibility between the symmetry on $\PsG$ and the Eilenberg-Moore construction, to lift these structures to $2$-categories of pseudoalgebras.

\section{Perspectives on compatibility between pseudomonads and monoidal structures}\label{Section alternative perspective}

\noindent The aim of this Section is to introduce opmonoidal pseudomonads via a categorification of the perspective on opmonoidal monads in Remark \ref{remark opmonoidal monads as monads in MonCat}, and to show that such structures are also suitably equivalent to pseudomonoids in a suitable monoidal sub-$2$-category of $\mathbf{Psmnd}(\mathbf{Gray})$, also to be described. In Section \ref{Section lifting monoidal structures}, we will use this to show that monoidal structures lift to $2$-categories of pseudoalgebras for semi-strict opmonoidal pseudomonads.
\\
\\
\noindent Theorem \ref{theorem strictifying strong monoidal pseudomonads} allows us to restrict our attention to simpler semi-strict opmonoidal pseudomonads without loss of generality up to biequivalence. This involves monoids and morphisms of monoids over the $\Gr$ tensor product, which we recall in Remark \ref{Remark Gray tensor product}. The $\mathbf{Gray}$-tensor product of pseudomonads is established in Theorem \ref{Theorem Gray tensor of pseudomonads}. Then in Theorem \ref{theorem semi-strict opmonoidal pseudomonads are pseudomonoids in the Gray monoid of pseudomonads}, we describe how semi-strict opmonoidal pseudomonads are certain pseudomonoids in this monoidal structure, categorifying the perspective on opmonoidal monads in Remark \ref{Opmonoidal monads as pseudomonoids in Mnd Cat}.

\begin{remark}\label{Remark Gray tensor product}
	Recall (Definition 1.2.1 of \cite{Miranda PhD}) that if $\A$ and $\B$ are $2$-categories, their $\mathbf{Gray}$-tensor product $\A \otimes \B$ can be presented as having
	
	\begin{itemize}
		\item generating $n$-cells given by pairs $(\Phi, \Psi)$ where $\Phi$ is a $k$-cell of $\A$, $\Psi$ is a $l$-cell of $\B$ and $k+l=n$, and
		\item relations corresponding to $(\Phi, -)$ being a cell in the $2$-category $\mathbf{Gray}(\B, \A\otimes \B)$ of the same dimension as $\Phi \in \A$, and $(?, \Psi)$ being a cell in the $2$-category $\mathbf{Gray}(\A, \A\otimes \B)$ of the same dimension as $\Psi \in \B$.
	\end{itemize} 
	
	\noindent We refer to generators of the form $(\Phi, \Psi)$ where both $\Phi$ and $\Psi$ are $1$-cells as \emph{interchangers}, and often denote them as $\Phi_\Psi$. fix terminology for relations in this presentation so that we can easily refer to these relations in proofs. If $b: \Psi \to \Psi'$ is a morphism in $\B$ we will refer to the relations pertaining to 
	
	\begin{itemize}
		\item the composition axiom for the $2$-functor $(?, \Psi): \A \to \A \otimes \B$, or the composition axiom for the pseudonatural transformation $(?, b): (?, \Psi) \Rightarrow (?, \Psi')$, on composable pairs of morphisms in $\A$ as \emph{functoriality relations},
		\item the unit axiom for $(?, \Psi): \A \to \A \otimes \B$, or the unit axiom for $(?, b): (?, \Psi) \Rightarrow (?, \Psi')$ as the \emph{unit relations}, and
		\item the naturality axiom for $(?, b)$ on $2$-cells $\alpha: a \Rightarrow a'$ as \emph{naturality relations}.
	\end{itemize}
	
	\noindent Similar terminology will be used for relations pertaining to pseudonaturality of $(a, -): (\Phi, -)\Rightarrow (\Phi', -)$. For simplicity, we use coherence for monoidal categories and work as though the $\mathbf{Gray}$-tensor product is strictly associative and unital. As such, a similar presentation can be given for a ternary $\Gr$-tensor product $\A \otimes \B \otimes \mathcal{C}$, involving $(k+l+m)$-cell generators for a $k$-cell in $\A$, $l$-cell in $\B$ and $m$-cell in $\mathcal{C}$. We freely use the relationship between the $\mathbf{Gray}$-tensor product and cubical pseudofunctors out of cartesian products. A reader looking for background on the $\Gr$-tensor product should consult Chapter 3 of \cite{Gurski Coherence in Three Dimensional Category Theory}. 
\end{remark}

\subsection{Opmonoidal pseudomonads as pseudomonads in the tricategory of monoidal bicategories}\label{Subsection opmonoidal pseudomonads as pseudomonads in the tricategory of monoidal bicategories}

\noindent The data and axioms involved in our definition of opmonoidal pseudomonads and higher maps between them can be conceptually understood in terms of higher maps between tricategories \cite{Low dimensional structures formed by tricategories}, and higher maps between pseudomonads \cite{Formal Theory of Pseudomonads}. We list all of the data, but refrain from displaying the axioms, instead giving references to definitions from which these conditions can be recovered.
\\
\\
\noindent There is a large tricategory $\mathbf{Tricat}_\text{oplax}$ categorifying $\mathbf{Bicat}_\text{icon}$ and related to the tricategory described in \cite{Low dimensional structures formed by tricategories}. This tricategory has the following cells.

\begin{itemize}
	\item its objects are given by small tricategories.
	\item its morphisms are given by oplax trihomomorphisms $F: \mathfrak{A} \rightarrow \mathfrak{B}$, which are lax trihomomorphisms $F: \mathfrak{A}^\text{co} \rightarrow \mathfrak{B}^\text{co}$.
	\item its $2$-cells are given by pseudo-icons. These are similar to trinatural transformations, but have
	\begin{itemize}
		\item identity $1$-cell components on objects,
		\item $2$-cell components $p_{f}: Ff \Rightarrow Gf$ on morphisms, which need not be part of adjoint equivalences,
		\item invertible $3$-cell components mediating unit and composition laws, or the varying of a $1$-cell along a $2$-cell. 
	\end{itemize}
	\item its $3$-cells are pseudo-icon modifications. These are like trimodifications but have identity $2$-cell components on objects, and $2$-cell components on morphisms.
\end{itemize}

\noindent The relationship between $\mathbf{Tricat}_\text{oplax}$ described above and the tricategory of tricategories described in \cite{Low dimensional structures formed by tricategories} is that the former has $1$-cells given by oplax, rather than lax, homomorphisms of tricategories. We are interested in the restriction of $\mathbf{Tricat}_\text{oplax}$ to tricategories with one object, or \emph{monoidal bicategories}. Here we call the $1$-cells \emph{opmonoidal pseudofunctors}, or \emph{opmonoidal $2$-functors} if their underlying pseudofunctors are strict. Similarly, we call the $2$-cells in this tricategory \emph{monoidal pseudonatural transformations}, and the $3$-cells \emph{monoidal modifications}. 

\begin{notation}
	\hspace{1mm}
	\begin{itemize}
		\item The tricategory of monoidal bicategories, opmonoidal pseudofunctors, monoidal pseudonatural transformations and monoidal modifications will be denoted $\mathbf{MonBicat}_\text{oplax}$.
		\item The sub-$\mathbf{Gray}$-category of $\mathbf{MonBicat}_\text{oplax}$ consisting of $\mathbf{Gray}$-monoids and opmonoidal $2$-functors will be denoted $\mathbf{Gray}\text{-}\mathbf{Monoids}_\text{oplax}$.
		\item Monoidal bicategories $\left(\mathcal{V}, \bigoplus, I, \alpha, \lambda, \rho, \mathfrak{l}, \mathfrak{m}, \mathfrak{r}, \mathfrak{p}\right)$ will be abbreviated as $\left(\mathcal{V}, \bigoplus, I\right)$. 
		\item The invertible modification component $\mathfrak{p}$ will be called the \emph{pentagonator}, $\mathfrak{r}$ will be called the \emph{right unit mediator}, $\mathfrak{m}$ will be called the \emph{middle unit mediator}, $\mathfrak{l}$ will be called the \emph{left unit mediator}.
	\end{itemize}
\end{notation}

\begin{definition}\label{definition colax monoidal pseudomonad}
	\hspace{1mm}
	\begin{enumerate}
		\item An \emph{opmonoidal pseudomonad} is a pseudomonad in the tricategory $\mathbf{MonBicat}_\text{oplax}$,
		\item The tricategory of opmonoidal pseudomonads is defined to be the tricategory $\mathbf{Psmnd}\left(\mathfrak{K}\right)$ of \cite{Formal Theory of Pseudomonads}, where $\mathfrak{K} = \mathbf{MonBicat}_\text{oplax}$.
	\end{enumerate}
\end{definition}

\noindent See Section 5.4.2 of \cite{Miranda PhD} for pseudomonads in a tricategory.

\begin{remark}\label{Remark explicit opmonoidal pseudomonad}
	Opmonoidal pseudomonads on a monoidal bicategory  $\left(\mathcal{V}, \bigoplus, I\right)$ can be described more explicitly, as follows.
	
	\begin{itemize}
		\item An opmonoidal pseudomonad has an underlying endomorphism in $\mathbf{MonBicat}_\text{oplax}$. This consists of a pseudofunctor $T: \mathcal{V} \rightarrow \mathcal{V}$, two pseudonatural transformations $\chi: T.\bigoplus \rightarrow \bigoplus.\left(T \times T\right)$ and $\iota: TI \rightarrow I$, and three invertible modifications $\gamma$, $\delta$ and $\omega$ with $2$-cell components as in Equation \ref{Equation opmonoidal pseudofunctor 3-cell components}. These data satisfy the conditions needed for $\left(T, \iota, \chi, \gamma, \delta, \omega\right): \left(\mathcal{V}, \bigoplus, I\right) \rightarrow \left(\mathcal{V}, \bigoplus, I\right)$ to be an opmonoidal pseudofunctor. These conditions in turn correspond to the oplax trihomomorphism axioms as listed in Definition 4.10 of \cite{Gurski Coherence in Three Dimensional Category Theory}.
	\end{itemize}
	
	\begin{equation}\label{Equation opmonoidal pseudofunctor 3-cell components}
	\begin{tikzcd}[font=\fontsize{9}{6}]
		T\left(\left(X \bigoplus Y\right) \bigoplus Z\right) \arrow[rr, "\chi_{X\bigoplus Y{,}Z}"] 
		\arrow[dd, "T\alpha_{X{,}Y{,}Z}"']
		&& T\left(X\bigoplus Y\right)\bigoplus TZ \arrow[rr, "\chi_{X{,}Y}\bigoplus TZ"]
		\arrow[dd, Rightarrow, shorten = 5, "\omega_{X{,}Y{,}Z}"]
		&& \left(TX \bigoplus TY\right) \bigoplus TZ
		\arrow[dd, "\alpha_{TX{,}TY{,}TZ}"]
		\\
		\\
		T\left(X \bigoplus \left(Y \bigoplus Z\right)\right) \arrow[rr, "\chi_{X{,}Y\bigoplus Z}"']
		&& TX\bigoplus T\left(Y \bigoplus Z\right) \arrow[rr, "TX \bigoplus\chi_{Y{,}Z}"'] && TX \bigoplus \left(TY \bigoplus TZ\right)
	\end{tikzcd}
\end{equation}
	$$\begin{tikzcd}[font=\fontsize{9}{6}]
		T\left(X \bigoplus I\right) \arrow[rr, "\chi_{X{,}I}"] 
		\arrow[rrdd, "T\rho_{X}"']&& TX \bigoplus TI \arrow[rr, "TX \bigoplus \iota"]\arrow[dd, Leftarrow, shorten = 5, "\delta_{X}"]
		&& TX \bigoplus I
		\arrow[lldd, "\rho_{TX}"]
		\\
		\\
		&& TX
	\end{tikzcd}\begin{tikzcd}[font=\fontsize{9}{6}]
		T\left(I \bigoplus X\right) \arrow[rr, "\chi_{I{,}X}"] 
		\arrow[rrdd, leftarrow, "T\lambda_{X}"']&& TI \bigoplus TX \arrow[rr, "\iota \bigoplus TX"]
		\arrow[dd, Rightarrow, shorten = 5, "\gamma_{X}"]
		&& I \bigoplus TX
		\arrow[lldd,leftarrow, "\lambda_{TX}"]
		\\
		\\
		&& TX
	\end{tikzcd}$$
	\begin{itemize}
		\item An opmonoidal pseudomonad has a unit, which is a $2$-cell in $\mathbf{MonBicat}_\text{oplax}$. This consists of a pseudonatural transformation $\eta: 1_{\mathcal{V}} \rightarrow T$ and two invertible modifications $\eta_{0}$ and $\eta_{-, ?}$ with components as in Equation \ref{Equation monoidal pseudonaturality constraints eta}. These data satisfy the axioms required for $\left(\eta, \eta_{0}, \eta_{2, -, ?}\right)$ to be a monoidal pseudonatural transformation. These conditions in turn correspond to the axioms for $(\eta, \eta_{0}, \eta_{-, ?})$ to be a pseudo-icon, as listed in Definitions 3 and 5 of \cite{Low dimensional structures formed by tricategories}.
	\end{itemize}
	
\begin{equation}\label{Equation monoidal pseudonaturality constraints eta}
	\begin{tikzcd}[row sep = 20, font=\fontsize{9}{6}]
		I \arrow[rr, "\eta_{I}"] \arrow[rrdd, bend right = 30, "1_{I}"']
		&{}\arrow[dd, Leftarrow, shorten = 10, "\eta_{0}"]
		& TI \arrow[dd, "\iota"]
		\\
		\\
		&{}& I&&{}
	\end{tikzcd} \begin{tikzcd}[font=\fontsize{9}{6}]
		X \bigoplus Y \arrow[dd, "\eta_{X \bigoplus Y}"'] \arrow[rr, "\eta_{X}\bigoplus Y"]
		&{}\arrow[dd, Rightarrow, shorten = 5, "\eta_{2{,}X{,} Y}"]
		& TX \bigoplus Y
		\arrow[dd, "TX \bigoplus \eta_{Y}"]
		\\
		\\
		T\left(X\bigoplus Y\right)
		\arrow[rr, "\chi_{{}_{X{,}Y}}"'] &{}& TX \bigoplus TY
	\end{tikzcd}
\end{equation}
	\begin{itemize}
		\item An opmonoidal pseudomonad has a multiplication, which is a $2$-cell in $\mathbf{MonBicat}_\text{oplax}$. This consists of a pseudonatural transformation $\mu: T^2 \rightarrow T$ and two invertible modifications $\mu_{0}$ and $\mu_{-, ?}$ with components as in Equation \ref{Equation monoidal pseudonaturality constraints mu}. These data are subject to the axioms required for $\left(\mu, \mu_{0}, \mu_{2, -, ?}\right)$ to be a monoidal pseudonatural transformation. These conditions in turn correspond to the conditions required for $(\mu, \mu_{0}, \mu_{-, ?})$ to be a pseudo-icon.
		
	\end{itemize}
	\begin{equation}\label{Equation monoidal pseudonaturality constraints mu}
		\begin{tikzcd}[row sep = 29, font=\fontsize{9}{6}]
			T^{2}I\arrow[dd, "\mu_{I}"'] \arrow[rr, "T\iota"]
			&{}\arrow[dd, Rightarrow, shorten = 5, "\mu_{0}"]& TI\arrow[dd, "\iota"]
			\\
			\\
			TI \arrow[rr, "\iota"'] &{}&I&&{}
		\end{tikzcd}\begin{tikzcd}[font=\fontsize{9}{6}]
			T^{2}\left(X\bigoplus Y\right)
			\arrow[dd, "\mu_{X \bigoplus Y}"']
			\arrow[rr, "T\chi_{{}_{X{,}Y}}"]
			&& T\left(TX \bigoplus TY\right) \arrow[rr, "\chi_{{}_{TX{,}TY}}"] \arrow[dd, Rightarrow, shorten = 10, "\mu_{2{,}X{,}Y}"]
			&& T^{2}X \bigoplus T^{2}Y
			\arrow[d, "\mu_{X}\bigoplus TY"]
			\\
			&&&&TX \bigoplus T^{2}Y
			\arrow[d, "TX\bigoplus \mu_{Y}"]
			\\
			T\left(X \bigoplus Y\right) \arrow[rrrr, "\chi_{{}_{X{,}Y}}"'] &&{}&&
			TX \bigoplus TY
		\end{tikzcd}
	\end{equation}
	\begin{itemize}
		\item Their left unitor $\mathbf{l}: \mu.T\eta \Rightarrow 1_{T}$, right unitor $\mathbf{r}: 1_{T} \Rightarrow \mu.\eta_{T}$, and associator $\mathbf{a}: \mu.T\mu \Rightarrow \mu.\mu_{T}$, are each invertible modifications subject to conditions required for them to be monoidal modifications. These conditions correspond to the axioms for pseudo-icon modifications listed in Definition 6 of \cite{Low dimensional structures formed by tricategories}.
	\end{itemize}
	
	\noindent The data described above are moreover required to satisfy commutativity conditions amounting to the axioms needed for $(\mathcal{V}, T, \eta, \mu, \mathbf{a}, \mathbf{l}, \mathbf{r})$ to be a pseudomonad, as listed in Definition 3.1 of \cite{Marmolejo Pseudodistributive 1}. 
\end{remark}

\begin{definition}\label{definition semi strict opmonoidal pseudomonads}
	\hspace{1mm}
	\begin{enumerate}
		\item An opmonoidal pseudomonad is \emph{semi-strict} if its underlying monoidal bicategory is a $\mathbf{Gray}$-monoid, and its underlying pseudofunctor is a $2$-functor.
		\item The $\mathbf{Gray}$-category of semi-strict opmonoidal pseudomonads is defined to be the $\mathbf{Gray}$-category $\mathbf{Psmnd}\left(\mathfrak{K}\right)$, where $\mathfrak{K} = \mathbf{Gray}-\mathbf{Monoids}_\text{oplax}$.
		\item
		A morphism in $\mathbf{MonBicat}_\text{oplax}$ will be called \emph{strong} if its underlying pseudofunctor is strict and the pseudonatural transformations $\iota$ and $\chi$ are part of adjoint equivalences.
		\item A morphism in $\mathbf{MonBicat}_\text{oplax}$ will be called \emph{strict} if it is strong and moreover the pseudonatural transformations $\iota$ and $\chi$ are identities, and the modifications $\gamma$, $\delta$ and $\omega$ are also identities.
		\item A monoidal pseudonatural transformation $\left(\phi, \phi_{0}, \phi_{-, ?}\right): F \rightarrow G$ will be called \emph{semi-strict} if the corresponding pseudo-icon has unitor modifications $\phi_{0}$ and compositor modifications $\phi_{-, ?}$ which are both just identities.
	\end{enumerate}
\end{definition}

\noindent We now state the main theorem of this Subsection, which is higher dimensional analogue of Proposition \ref{strictification 1D setting}. The proof uses coherence for tricategories, and is deferred to Appendix \ref{Appendix Coherence}.

\begin{theorem}\label{theorem strictifying strong monoidal pseudomonads}
	\hspace{1mm}
	\begin{enumerate}
		\item 
		Every opmonoidal pseudomonad $(\mathcal{V}, \bigoplus, I, T, \eta, \mu, \mathbf{a}, \mathbf{l}, \mathbf{r})$ is biequivalent in the tricategory $\mathbf{Psmnd}(\mathbf{MonBicat}_\text{oplax})$ to a semi-strict opmonoidal pseudomonad.
		\item If the underlying endomorphism of $(\mathcal{V}, \bigoplus, I, T, \eta, \mu, \mathbf{a}, \mathbf{l}, \mathbf{r})$ is strong, then a biequivalent semi-strict opmonoidal pseudomonad can be chosen for part (1) in which moreover the underlying endomorphism is strict, and the unit and multiplication are semi-strict.
	\end{enumerate} 
\end{theorem}

\noindent As a consequence of Theorem \ref{theorem strictifying strong monoidal pseudomonads} part (1), once we show that monoidal structures underlying biequivalent semi-strict opmonoidal pseudomonads $(\overline{\mathcal{V}}, \overline{T})$ lift to $2$-categories of pseudoalgebras $\overline{\mathcal{V}}^{\overline{T}}$, it will follow that these lifted structures transport to the bicategory of pseudoalgebras $\mathcal{V}^T$ for the original opmonoidal pseudomonad $(\mathcal{V}, T)$. 

\subsection{The $\mathbf{Gray}$-tensor product of pseudomonads}\label{Subsection the Gray tensor product of pseudomonads}
See \cite{Formal Theory of Pseudomonads} for the structure of the $\mathbf{Gray}$-category of pseudomonads $\mathbf{Psmnd}(\mathbf{Gray})$. We restrict this to a $2$-category on which we describe a $\mathbf{Gray}$-monoid structure. Semi-strict opmonoidal pseudomonads will correspond to pseudomonoids in this $\mathbf{Gray}$-monoid.
\\
\\
\noindent Recall that the $k$-cells of a $\mathbf{Gray}$-category for $k \leq 2$ typically only form a sesquicategory \cite{Stell modelling term rewriting via sesquicategories}, rather than a $2$-category, since middle-dour interchange may only hold up to isomorphism. We describe an appropriate restriction of the $\mathbf{Gray}$-category $\mathbf{Psmnd}(\mathbf{Gray})$ to a $2$-category, $\mathbf{Psmnd}(\mathbf{Gray})_\text{ast}$. The objects and morphisms in $\mathbf{Psmnd}(\mathbf{Gray})_\text{ast}$ are the same as in $\mathbf{Psmnd}(\mathbf{Gray})$, but the $2$-cells are just those which are ambistrict, in the sense of Definition \ref{definition ambistrict 2-cell}, to follow.

\begin{definition}\label{definition ambistrict 2-cell}
	Let $\mathcal{A}$ be a sesquicategory and let $\psi: h \Rightarrow k: Y \rightarrow Z$ be a $2$-cell in $\mathcal{A}$.
	
	\begin{enumerate}
		\item $\psi$ is called \emph{strict} if for all $\phi: f \Rightarrow g: X \rightarrow Y$, the following diagram commutes in the hom-category $\mathcal{A}\left(X, Z\right)$.
		
		$$\begin{tikzcd}
			hf \arrow[d, "\psi.f"']
			\arrow[r, "h.\phi"]
			& hg\arrow[d, "\psi.g"]
			\\
			kf \arrow[r, "k.\phi"'] & kg
		\end{tikzcd}$$
		
		\item $\psi$ is called \emph{costrict} if it is strict as a $2$-cell $\psi: h \Rightarrow k: Z \rightarrow Y$ in $\mathcal{A}^\text{op}$.
		\item $\psi$ is called \emph{ambistrict} if it is both strict and costrict.
	\end{enumerate}
\end{definition}

\noindent In a $2$-category, every $2$-cell is ambistrict. It is easy to see that the class of ambistrict $2$-cells in a sesquicategory contains the identities and is closed under vertical composition and whiskering. As such, restricting a sesquicategory to its ambistrict $2$-cells results in a $2$-category.

\begin{example}\label{Example ambistrict 2-cells in Gray}
	As discussed in \cite{Icons}, a pseudonatural transformation is
	
	\begin{itemize}
		\item strict if and only if its $2$-cell components are identities,
		\item costrict if and only if its $1$-cell components are identities.
	\end{itemize} 
	
	\noindent In particular, the only ambistrict $2$-cells in the underlying sesquicategory of $\mathbf{Gray}$ are identity $2$-natural transformations.
\end{example}

\begin{proposition}\label{Proposition characterising ambistrict pseudomonad transformations}
	A pseudomonad transformation $\left(\phi, \tilde{\phi}\right): \left(G, g, u, m\right) \Rightarrow \left(H, h, v, n\right)$ is ambistrict if and only if its $2$-cell component $\phi: G \Rightarrow H$ is an identity $2$-natural transformation.
\end{proposition}

\begin{proof}
	Follows by inspecting the interchangers of $2$-cells in $\mathbf{Psmnd}\left(\mathbf{Gray}\right)$, which as described in the proof of Lemma 2.9 in \cite{Formal Theory of Pseudomonads} are given by interchangers in $\mathbf{Gray}$ of their underlying pseudonatural transformations.
\end{proof}

\begin{notation}
	Let $\mathbf{Psmnd}(\mathbf{Gray})_\text{ast}$ denote the $2$-category of pseudomonads, pseudomonad morphisms, and ambistrict pseudomonad transformations.
\end{notation}

\noindent The remainder of this Subsection builds towards the symmetric $\mathbf{Gray}$-monoid structure on $\PsG_\text{ast}$, which is established in Theorem \ref{Theorem Gray tensor of pseudomonads}. In our proofs we give full diagrammatic details for the constructions, but explain in words the steps in the proofs needed to check that things are well-defined. Due to space constraints, we only give details for equations between cells that are not of the highest dimension.
\\
\\
\noindent As part of the $\mathbf{Gray}$-monoid structure on $\PsG_\text{ast}$, we require a $2$-functor $\mathbf{Psmnd}(\Gr)_\text{ast} \otimes \mathbf{Psmnd}(\Gr)_\text{ast} \to \mathbf{Psmnd}(\Gr)_\text{ast}$. Lemma \ref{Gray tensor of two pseudomonads}, to follow, describes the behaviour of this $2$-functor on objects.

\begin{lemma}\label{Gray tensor of two pseudomonads}
	Let $(\A, S)$, $(\B, T)$ be pseudomonads in $\mathbf{Gray}$. Then there is a pseudomonad on $\A\otimes \B$ with underlying endomorphism $S\otimes T$, unit $\eta^{S}\otimes \eta^{T}$, multiplication $\mu^{S}\otimes \mu^{T}$, left and right unitors given by the modifications depicted below
	
	$$\begin{tikzcd}[font=\fontsize{9}{6}, column sep = 20]
		S\otimes T \arrow[rrdd, bend right = 20, "1"']
		\arrow[rr, "(\eta^{S}.1)\otimes 1"]
		&{}\arrow[dd, Rightarrow, shorten = 15, "\lambda^S\otimes 1"]
		& S^2\otimes T
		\arrow[rr, "1\otimes (1.\eta^T)"]
		\arrow[dd, "\mu^S \otimes 1"description]
		&{}\arrow[dd, Rightarrow,shorten =15, "{\mu^{S}}_{(1.\eta^{T})}"]
		& S^2\otimes T^2
		\arrow[dd, "\mu^{S}\otimes 1"]
		\\
		\\
		&{}&S\otimes T
		\arrow[rr, "1\otimes (1.\eta^{T})"]
		\arrow[rrdd, bend right = 20, "1"']
		&{}\arrow[dd, Rightarrow, shorten = 15, "1\otimes \lambda^{T}"]
		&S\otimes T^2
		\arrow[dd, "1\otimes \mu^{T}"]
		\\
		\\
		&&&{}&S\otimes T&{}
	\end{tikzcd}\begin{tikzcd}[font=\fontsize{9}{6}, column sep = 20]
		S\otimes T \arrow[rrdd, bend right = 20, "1"']
		\arrow[rr, "(1.\eta^{S})\otimes 1"]
		&{}\arrow[dd, Rightarrow, shorten = 15, "\rho^S\otimes 1"]
		& S^2\otimes T
		\arrow[rr, "1\otimes (\eta^T.1)"]
		\arrow[dd, "\mu^S \otimes 1"description]
		&{}\arrow[dd, Rightarrow,shorten =15, "{\mu^{S}}_{(\eta^{T}.1)}"]
		& S^2\otimes T^2
		\arrow[dd, "\mu^{S}\otimes 1"]
		\\
		\\
		&{}&S\otimes T
		\arrow[rr, "1\otimes (\eta^{T}.1)"]
		\arrow[rrdd, bend right = 20, "1"']
		&{}\arrow[dd, Rightarrow, shorten = 15, "1\otimes \rho^{T}"]
		&S\otimes T^2
		\arrow[dd, "1\otimes \mu^{T}"]
		\\
		\\
		&&&{}&S\otimes T
	\end{tikzcd}$$
	and associator given by the modification depicted below.
	
	$$\begin{tikzcd}[font=\fontsize{9}{6}]
		S^3\otimes T^3\arrow[rr, "(1.\mu^{S})\otimes 1"]\arrow[dd, "(\mu^{S}.1)\otimes 1"']
		&{}\arrow[dd, Rightarrow,shift right = 5,  shorten = 15, "\mathbf{a}^S\otimes 1"]
		&S^2\otimes T^3\arrow[rr, "1\otimes(1.\mu^T)"]\arrow[dd, "\mu^S \otimes 1"description]
		&{}\arrow[dd, Rightarrow, shorten = 15, "{(\mu^S)}_{(1.\mu^T)}"]
		&S^2\otimes T^2\arrow[dd, "\mu^S \otimes 1"]
		\\
		\\
		S^2\otimes T^3\arrow[dd, "1\otimes (\mu^{T}.1)"']\arrow[rr, "\mu^{S}\otimes 1"]
		&{}\arrow[dd, Rightarrow,shift right = 5,  shorten = 15, "{(\mu^{S})}_{(\mu^{T}.1)}"]
		&S\otimes T^3\arrow[dd, "1\otimes (\mu^{T}.1)" description]\arrow[rr, "1\otimes (1. \mu^{T})"]&{}\arrow[dd, Rightarrow, shorten = 15, "1\otimes \mathbf{a}^{T}"]
		&S\otimes T^2\arrow[dd, "1\otimes \mu^T"]
		\\
		\\
		S^2\otimes T^2\arrow[rr, "\mu^S \otimes 1"']
		&{}&S\otimes T^{2}\arrow[rr, "1\otimes \mu^T"']&{}&S\otimes T
	\end{tikzcd}$$
\end{lemma}

\begin{proof}
	We give details for the unit law in Appendix \ref{Appendix proof of Gray tensor product of two pseudomonads}. The associativity law follows by similar techniques, using interchange and the associativity laws for both pseudomonads. 
\end{proof}

\noindent As part of defining the remaining data of the $2$-functor $\mathbf{Psmnd}(\Gr)_\text{ast} \otimes \mathbf{Psmnd}(\Gr)_\text{ast} \to \mathbf{Psmnd}(\Gr)_\text{ast}$, we need to define $2$-functors out of each factor when a pseudomonad is fixed in the other factor. Lemma \ref{Lemma - tensor B T 2-functor}, to follow, defines such a $2$-functor when the pseudomonad in the second factor is fixed. Following this, Lemma \ref{Lemma A S tensor - 2-functor} describes the analogous $2$-functor when the first factor is fixed.

\begin{lemma}\label{Lemma - tensor B T 2-functor}
	Consider a pseudomonad transformation 
	$$\begin{tikzcd}[font=\fontsize{9}{6}]
		(\A{,}S)\arrow[rr, bend left = 30, "(G{,}g{,}g_{0}{,}g_{2})"name=A]\arrow[rr, bend right = 30, "(G{,}g'{,}g_{0}'{,}g_{2}')"'name=B] &&(\A'{,}S')\arrow[from=A, to=B, Rightarrow, shorten = 10, "\left(\phi{,}\tilde{\phi}\right)"]
	\end{tikzcd}$$  

\noindent and $(\B, T)$ be another pseudomonad. Let $(\A\otimes\B, S\otimes T)$ and $(\A'\otimes\B, S'\otimes T)$ be the pseudomonads constructed as in Lemma \ref{Gray tensor of two pseudomonads}. 
	
	\begin{enumerate}
		\item There is a morphism of pseudomonads $(G\otimes \B, g\otimes \B): (\A\otimes\B, S\otimes T)\to (\A'\otimes\B, S'\otimes T)$ whose data are given as follows.
		\begin{itemize}
			\item The underlying morphism is given by the $2$-functor $G\otimes 1_{\B}: \A \otimes \B \to \A'\otimes \B$.
			\item The underlying $2$-cell $g\otimes T$ is given by
			 
			$$ \begin{tikzcd}[font=\fontsize{9}{6}]
				(S'\otimes T).(G\otimes 1_{\B})\arrow[r, equal]& S'G\otimes T \arrow[r, "g\otimes T"] &GS\otimes T \arrow[r, equal] &(G\otimes 1_{\B}).(S\otimes T)
			\end{tikzcd}$$
			
			\item The unit given as a pasting of $g_{0}$ with the interchanger of $\eta_{G}^{S'}$ and $\eta^{T}$, and has components on $(X, Y) \in \mathcal{A} \otimes \mathcal{B}$ as depicted below.
			
			$$\begin{tikzcd}[font=\fontsize{9}{6}, column sep = 40]
				(GX{,}Y)\arrow[dd, shift left =3, Rightarrow, shorten = 7, "(g_{0{,}X}{,}Y)"]
				\arrow[rr, "(G\eta_{X}{,}Y)"]
				\arrow[dd, bend right = 50, "(\eta_{GX}{,}Y)"']
				&{}
				&(GS{,}Y)
				\arrow[rrdd, "(GSX{,}\eta_{Y})"]
				\arrow[dd, Rightarrow, shorten = 5, "(g_{X}{,}\eta_{Y})"]
				\\
				\\
				(S'GX{,}Y)
				\arrow[rr, "(SG'X{,}\eta_{Y})"']
				\arrow[rruu, "(g_{X}{,}Y)"']
				&{}&(S'GX{,}TY)
				\arrow[rr, "(g_{X}{,}TY)"']
				&{}& (GSX{,}TY)
			\end{tikzcd}$$
			
			\item The multiplication given as a pasting of $g_{2}$ with the interchanger of $g$ and $\mu^{T}$, and has components on $(X, Y) \in \mathcal{A} \otimes \mathcal{B}$ as depicted below. 
			
			$$\begin{tikzcd}[font=\fontsize{9}{6}]
				({S'}^{2}GX{,}T^2Y)\arrow[rr,"(S'g_X{,}T^2Y)"]\arrow[dd, "(\mu_{GX}{,}T^2Y)"']
				&&(S'GSX{,}T^2Y)\arrow[rr, "(g_{SX}{,}T^2Y)"]\arrow[dd, Rightarrow, shorten = 15, "(g_{2{,}X}{,}T^2Y)"]
				&&(GS^2X{,}T^2Y)\arrow[dd, "(G\mu_{X}{,}T^2Y)"]
				\\
				\\
				({S'}GX{,}T^2Y)\arrow[rrrr, "(g_X{,}T^2Y)"']
				\arrow[dd, "({S'}GX{,}\mu_Y)"']
				&&{}\arrow[dd, Rightarrow, shorten = 15, "(g_{X}{,}\mu_Y)"]
				&&(GSX{,}T^2Y)\arrow[dd, "(GSX{,}\mu_Y)"]
				\\
				\\
				({S'}GX{,}TY)\arrow[rrrr, "(g_X{,}TY)"']&&{}&&(GSX{,}TY)
			\end{tikzcd}$$
		\end{itemize} 
		\item Let $(G\otimes \B, g'\otimes \B): (\A\otimes \B, S \otimes T) \to (\A'\otimes \B')$ be defined similarly to part (1). Then there is a pseudomonad transformation $(\phi\otimes \B, \tilde{\phi}\otimes T):(G\otimes \B, g\otimes \B) \to (G\otimes \B, g'\otimes \B)$, which is ambistrict if $\left(\phi, \tilde{\phi}\right)$ is ambistrict.
		\item This assignment extends to a $2$-functor $- \otimes (\B, T): \Past \to \Past$, which is the identity if $\B$ is terminal.
	\end{enumerate}
\end{lemma}

\begin{proof}
	Part (1) uses the axioms for $(G, g, g_{0}, g_{2})$, and the interchange laws for the $2$-cell components of the pseudomonads and of $g_{0}$ and $g_{2}$. Part (2) uses the axioms for the pseudomonad transformation $(1_{G}, \tilde{\phi})$. Part (3) follows from the functoriality aspects in the presentation of the $\mathbf{Gray}$-tensor product.
\end{proof}

\begin{lemma}\label{Lemma A S tensor - 2-functor}
	Consider a pseudomonad transformation 
	$$\begin{tikzcd}[font=\fontsize{9}{6}]
		(\B{,}T)\arrow[rr, bend left = 30, "(H{,}h{,}h_{0}{,}h_{2})"name=A]\arrow[rr, bend right = 30, "(H'{,}h'{,}h_{0}'{,}h_{2}')"'name=B] &&(\B'{,}T')\arrow[from=A, to=B, Rightarrow, shorten = 10, "(\psi{,}\tilde{\psi})"]
	\end{tikzcd}$$  
	
	\noindent and $(\A, S)$ be another pseudomonad. Let $(\A\otimes\B, S\otimes T)$ and $(\A\otimes\B', S\otimes T')$ be the pseudomonads constructed as in Lemma \ref{Gray tensor of two pseudomonads}. 
	
	\begin{enumerate}
		\item There is a morphism of pseudomonads $(\A\otimes H, \A\otimes h): (\A\otimes\B, S\otimes T)\to (\A\otimes\B', S\otimes T')$ whose data are given as follows.
		\begin{itemize}
			\item The underlying morphism is given by the $2$-functor $1_\A \otimes H: \A \otimes \B \to \A\otimes \B'$.
			\item The underlying $2$-cell $S\otimes h$ is given by
			
			$$ \begin{tikzcd}
				(S\otimes T').(1_\A\otimes H)\arrow[r, equal]& S\otimes T'H  \arrow[r, "S\otimes h"] &S\otimes HT \arrow[r, equal] &(1_\A\otimes H).(S\otimes T)
			\end{tikzcd}$$
			
			\item The unit given as a whiskering of $h_{0}$, and has components on $(X, Y) \in \mathcal{A} \otimes \mathcal{B}$ as depicted below.
			
			$$\begin{tikzcd}[font=\fontsize{9}{6}, column sep = 40]
				(X{,}HY)
				\arrow[dd, "(\eta_{X}{,}HY)"']
				\\
				&&{}\arrow[d, Rightarrow, shorten = 5, "(S'X{,}h_{0{,}Y})"]
				\\
				(S'X{,}HY)\arrow[rr, "(S'X{,}H\eta_{Y})"']
				\arrow[rrrr, bend left = 30, "(SX{,}\eta_{HY})"]
				&{}&(SX{,}T'HY)
				\arrow[rr, "(S'X{,}h_Y)"']
				&{}& (S'X{,}HTY)
			\end{tikzcd}$$
			
			\item The multiplication given as a pasting of $h_{2}$ with various interchangers in $\mathbf{Gray}$, and has components on $(X, Y) \in \mathcal{A} \otimes \mathcal{B}$ as depicted below. 
			
			$$\begin{tikzcd}[font=\fontsize{9}{6}]
				&({S}^{2}X{,}{T'}^2HY)\arrow[rr,"({S^2}X{,}T'hY)"]\arrow[dd, "(S^2X{,}\mu_{HY})"]
				\arrow[ldd, bend right = 20, "(\mu_{X}{,}{T'}^{2}Y)"']
				&&({S}^2X{,}T'HTY)\arrow[rr, "({S}^2X{,}h_{TY})"]\arrow[dd, Rightarrow, shorten = 15, "(S^2X{,}h_{2{,}Y})"]
				&&(S^2X{,}HT^2Y)\arrow[dd, "(S^2X{,}H\mu_Y)"']
				\arrow[rdd, bend left = 20, "(\mu_{X}{,}HT^2Y)"]
				\\
				\\
				(SX{.}{T'}^2HY)\arrow[rdd, bend right = 20, "(SX{,}h_{Y})"']
				&({S}^2X{,}{T'}HY)
				\arrow[l, Rightarrow, shift right = 4, shorten = 5, "{(\mu_{X}{,}\mu_{HY})}^{-1}"']
				\arrow[rrrr, "(S^2X{,}\mu_{HY})"']
				\arrow[dd, "(\mu_X{,}T'HY)"]
				&&{}\arrow[dd, Rightarrow, shorten = 15, "{(\mu_{X}{,}h_Y)}^{-1}"]
				&&(S^2X{,}HTY)\arrow[dd, "(\mu_{X}{,}HTY)"']
				&(SX{,}HT'^2Y)
				\arrow[l, Rightarrow, shorten =5,shift right = 4, "{(\mu_{X}{,}H\mu_{Y})}"']
				\arrow[ddl, bend left = 20, "(SX{,}H\mu_{Y})"]
				\\
				\\
				&({S}X{,}T'HY)\arrow[rrrr, "(SX{,}h_Y)"']&&{}&&(SX{,}HTY)
			\end{tikzcd}$$
		\end{itemize} 
		\item Let $(\A\otimes H', \A\otimes h'): (\A\otimes \B, S \otimes T) \to (\A\otimes \B')$ be defined similarly to part (1). Then there is a pseudomonad transformation $(\A, \tilde{\psi}):(\A\otimes H, \A\otimes h) \to (\B\otimes H', \A\otimes h')$, which is ambistrict if $\left(\psi, \tilde{\psi}\right)$ is.
		\item This assignment extends to a $2$-functor $(\A, S)\otimes -: \Past \to \Past$, which is the identity if $\A$ is terminal.
	\end{enumerate}
\end{lemma}

\begin{proof}
	The proof strategy is analogous to that for Lemma \ref{Lemma A S tensor - 2-functor}. Part (1) uses the axioms for $(H, h, h_{0}, h_{2})$, and the interchange laws for the $2$-cell components of the pseudomonads and of $h_{0}$ and $h_{2}$. Part (2) uses the axioms for the pseudomonad transformation $(1_{H}, \tilde{\psi})$. Part (3) follows from the functoriality aspects in the presentation of the $\mathbf{Gray}$-tensor product. 
\end{proof}

\noindent The $2$-functor $\PsG_\text{ast} \otimes \PsG_\text{ast} \to \PsG_\text{ast}$ also has images on interchangers. We describe these in Lemma \ref{interchanger between morphisms of pseudomonads in Gray monoid structure on Psmnd Gray}, to follow.

\begin{lemma}\label{interchanger between morphisms of pseudomonads in Gray monoid structure on Psmnd Gray}
	Let $\left(G, g, g_{0}, g_{2}\right): \left(\mathcal{A}, S\right) \rightarrow \left(\mathcal{A}', S'\right)$ and $\left(H, h, h_{0}, h_{2}\right): \left(\mathcal{B}, T\right) \rightarrow \left(\mathcal{B}', T'\right)$ be morphisms of pseudomonads. Then there is an ambistrict pseudomonad transformation as depicted below left, whose component on $(X, Y) \in\A'\otimes\B'$ is given by the $2$-cell depicted below right. 
	
	$$\begin{tikzcd}[font=\fontsize{9}{6}]
	(\A\otimes \B{,}S\otimes T) \arrow[rr, "(G\otimes \B{,}g\otimes \B)"]\arrow[dd, "(\A\otimes H{,}\A\otimes h)"'] &{}\arrow[dd, Rightarrow, shorten = 10, "\left(1_{HG}{,} h_{g}\right)"]
	& (\A'\otimes \B{,}S'\otimes T)\arrow[dd, "(\A'\otimes H{,}\A'\otimes h)"]
	\\
	\\
	(\A\otimes \B'{,}S\otimes T')\arrow[rr, "(G\otimes \B'{,}g\otimes \B')"'] &{}& (\A'\otimes \B'{,}S'\otimes T')&{}
\end{tikzcd}\begin{tikzcd}[font=\fontsize{9}{6}]
(S'GX{,}T'HY)\arrow[dd, "(S'GX{,}h_{Y})"']\arrow[rr, "(g_{X}{,}T'HY)"]&{}\arrow[dd, Rightarrow, shorten = 10, "{(g_{X}{,}h_{Y})}^{-1}"]&(GSX{,}T'HY)\arrow[dd, "(GSX{,}h_{Y})"]
\\
\\
(S'GX{,}HTY)\arrow[rr, "(g_{X}{,}HTY)"']&{}&(GSX{,}HTY)
\end{tikzcd}$$
	
\end{lemma}

\begin{proof}
	For the unit law, the right side of the equation given in \cite{Formal Theory of Pseudomonads} corresponds to the pasting displayed below.
	
	$$\begin{tikzcd}[font=\fontsize{9}{6}]
		&&&&(S'GX{,}T'HY) \arrow[rr, "(g_{X}{,}1)"]\arrow[dd, "(1{,}h_{Y})"] &{}\arrow[dd, Leftarrow, shorten = 10, "{(g_{X}{,}h_{Y})}^{-1}"]& (GSX{,}T'HY)\arrow[dd, "(1{,}h_{Y})"]
		\\
		&&&{}\arrow[d, Leftarrow, shorten = 2, "(1{,}h_{0{,}Y})"]
		\\
		(GX{,}HY)\arrow[rr, "(\eta_{GX}{,}1)"]\arrow[rrdd, bend right = 10, "(G\eta_{X}{,}1)"']
		&&(S'GX{,}HY)\arrow[dd, Leftarrow, shift right = 5,  shorten = 10, "(g_{0{,}X}{,}1)"']
		\arrow[rr, "(1{,}H\eta_{Y})"']
		\arrow[rruu, "(1{,}\eta_{HY})"]
		\arrow[dd, "(g_{X}{,}1)"]
		&{}\arrow[dd, Leftarrow, shorten = 10,"{(g_{X}{,}H\eta_{Y})}"]
		&(S'GX{,}HTY)
		\arrow[rr, "(g_{X}{,}1)"']
		&{}& (GSX{,}HTY)
		\\
		\\
		&&(GSX{,}HY)\arrow[rrrruu, bend right = 10, "(1{,}H\eta_{Y})"']
		&{}
	\end{tikzcd}$$
	
	\noindent By the interchange law corresponding to pseudonaturality of $(g_{X}, -)$ on the $2$-cell $h_{0, Y}$, this is equal to the pasting diagram below. This is the required left side of the equation in the unit law.
	
		$$\begin{tikzcd}[font=\fontsize{9}{6}]
		&&&&(S'GX{,}T'HY) \arrow[rr, "(g_{X}{,}1)"]\arrow[dd, Leftarrow, shorten = 15, red,"{(g_{X}{,}\eta_{HY})}"']
		&& (GSX{,}T'HY)\arrow[dd, "(1{,}h_{Y})"]
		\\
		&&&{}
		\\
		(GX{,}HY)\arrow[rr, "(\eta_{GX}{,}1)"]\arrow[rrdd, bend right = 10, "(G\eta_{X}{,}1)"']
		&&(S'GX{,}HY)\arrow[dd, Leftarrow, shift right = 5,  shorten = 10, "{(g_{0{,}X}{,}1)}"']
		\arrow[rruu, "(1{,}\eta_{HY})"]
		\arrow[dd, "(g_{X}{,}1)"]
		&{}&{}\arrow[dd, Leftarrow, shorten = 15,red, "(1{,}h_{0{,}Y})"']
		&{}& (GSX{,}HTY)
		\\
		\\
		&&(GSX{,}HY)\arrow[rrrruu, bend right = 10, "(1{,}H\eta_{Y})"']\arrow[rrrruuuu,blue, "(1{,}\eta_{HY})"description]
		&&{}
	\end{tikzcd}$$

	\noindent The multiplication law follows similarly, using the interchange law in the $\mathbf{Gray}$-tensor product $\A'\otimes\B'$ for the $2$-cell $h_{2, Y}$.
\end{proof}

\noindent Lemma \ref{symmetry for Gray tensor product of pseudomonads}, to follow, describes the symmetry for the $\Gr$-monoid structure on $\PsG_\text{ast}$.

\begin{lemma}\label{symmetry for Gray tensor product of pseudomonads}
	Let $(\mathcal{A}, S)$ and $(\mathcal{B}, T)$ be pseudomonads in $\mathbf{Gray}$ and let $\tau_{\A, \B}: \mathcal{A} \otimes \mathcal{B}\to \B \otimes \A$ denote the symmetry for the $\mathbf{Gray}$-tensor product. Then $\tau$ underlines an invertible morphism of pseudomonads whose $2$-cell component is the identity and whose unit and multiplication constraints on $(X, Y)$ are given by the interchangers ${(\eta_{X})}_{(\eta_{Y})}$ and ${(\mu_{X})}_{(\mu_{Y})}$ respectively.
\end{lemma}

\begin{proof}
	This follows from $\mathbf{Gray}$-naturality of $(\A, \B) \mapsto \tau_{\A, \B}$.
\end{proof}

\noindent We are now ready to describe the symmetric $\Gr$-monoid structure on $\PsG_\text{ast}$.

\begin{theorem}\label{Theorem Gray tensor of pseudomonads}
	\hspace{1mm}
	\begin{enumerate}
		\item There is a $\mathbf{Gray}$-monoid structure on ${\mathbf{Psmnd}\left(\mathbf{Gray}\right)}_{\text{ast}}$, with tensor given at the level of objects as in Lemma \ref{Gray tensor of two pseudomonads}, tensor on the left given as per Lemma \ref{Lemma - tensor B T 2-functor}, tensor on the right given as per Lemma \ref{Lemma A S tensor - 2-functor}, and interchangers given as described in Lemma \ref{interchanger between morphisms of pseudomonads in Gray monoid structure on Psmnd Gray}. 
		\item The pseudomonad morphism described in Lemma \ref{symmetry for Gray tensor product of pseudomonads} equips the underlying monoidal category with a symmetry, in such a way that the $\mathbf{Gray}$-monoid ${\mathbf{Psmnd}\left(\mathbf{Gray}\right)}_{\text{ast}}$ becomes symmetric, with a $2$-natural braiding and identity modification components.
	\end{enumerate}
\end{theorem}

\begin{proof}
	For part (1), by Lemmas \ref{Lemma - tensor B T 2-functor} and \ref{Lemma A S tensor - 2-functor}, it remains to show that the interchangers are functorial and natural in both arguments. But this follows from the analogous properties for the $\mathbf{Gray}$-tensor product. $2$-naturality of the braiding uses naturality of the braiding for the $\mathbf{Gray}$-tensor product. The braid equations and symmetry use these properties for the $\mathbf{Gray}$-tensor product on $\mathbf{Gray}$.
\end{proof}
\subsection{Semi-strict opmonoidal pseudomonads as pseudomonoids in the $\mathbf{Gray}$-monoid of pseudomonads}\label{Section monoidal pseudomonads as pseudomonoids in pseudomonads}

We explain how a pseudomonoid in the $\mathbf{Gray}$-monoid of Theorem \ref{Theorem Gray tensor of pseudomonads} is precisely a semi-strict opmonoidal pseudomonad. Indeed, the correspondence extends what has been described in Subsection \ref{Subsection opmonoidal monads as monoids in Mnd Cat}. We compare the data involved first, and then discuss the correspondence between conditions on these data. This correspondence is entirely an exercise in unravelling definitions, and does not require any further calculation beyond an inspection of the axioms for both sets of data. Again due to space constraints, we give only selected examples of the equalities of pasting diagrams corresponding to these axioms.
\\
\\
\noindent DATA
\\

\begin{itemize}
	\item The monoidal product consists of a $2$-functor $\bigoplus: \mathcal{A} \otimes \mathcal{A} \rightarrow \mathcal{A}$, a pseudonatural transformation $\chi: T.\bigoplus \rightarrow \bigoplus.\left(T \bigoplus T\right)$, and two invertible modifications $\eta_{-, ?}$ and $\mu_{-, ?}$ as depicted below.
	
	$$\begin{tikzcd}[font = \fontsize{9}{6}, row sep = 12]
		\mathcal{A} \otimes \mathcal{A}
		\arrow[rr, "\bigoplus"]
		\arrow[dd, bend right = 30, "T\otimes T"']
		&{}\arrow[dd, Rightarrow, shorten = 5, "\chi"]
		& \mathcal{A}
		\arrow[dd, bend left = 30, "1_{\mathcal{A}}"name=A]
		\arrow[dd, bend right = 30, "T"'name=B]
		\arrow[from=A, to=B, Rightarrow, shorten =5,"\eta"]
		\\
		&&&&\cong_{\eta_{-{,}?}}&{}
		\\
		\mathcal{A} \otimes \mathcal{A} \arrow[rr, "\bigoplus"'] &
		{}& \mathcal{A}
	\end{tikzcd}\begin{tikzcd}[row sep = 15,font = \fontsize{9}{6}]
		\mathcal{A} \otimes \mathcal{A}
		\arrow[rr, "\bigoplus"]
		\arrow[dd, bend left = 30, "1_{\mathcal{A}}\otimes 1_{\mathcal{A}}"name=A]
		\arrow[dd, bend right = 30, "T\otimes T"'name=B]
		&{}
		& \mathcal{A}
		\arrow[dd, bend left = 30, "1_{\mathcal{A}}"]
		\\
		&=
		\\
		\mathcal{A} \otimes \mathcal{A} \arrow[rr, "\bigoplus"'] &
		{}& \mathcal{A}
		\arrow[from=A, to=B,Rightarrow, shorten=5, "\eta \otimes \eta"]
	\end{tikzcd}$$
	
	$$\begin{tikzcd}[row sep = 15, font = \fontsize{9}{6}]
		\mathcal{A} \otimes \mathcal{A} \arrow[rr, "\bigoplus"] \arrow[dddd, "T\otimes T"']\arrow[ddr, "T\otimes T"]
		&{}\arrow[dd, Rightarrow, shorten = 5, shift left = 3, "\chi"]
		&\mathcal{A}\arrow[ddr, "T"]
		\\
		{}\arrow[dd, Rightarrow, shift left = 5, "\mu\otimes\mu"]
		\\
		&\mathcal{A} \otimes \mathcal{A} \arrow[rr, "\bigoplus"] \arrow[ddl, "T\otimes T"]\arrow[dd, Rightarrow, shorten = 5, shift left = 3, "\chi"]
		&&\mathcal{A}\arrow[ddl, "T"]&\cong_{\mu_{-{,}?}}&{}
		\\
		{}
		\\
		\mathcal{A} \otimes \mathcal{A} \arrow[rr, "\bigoplus"'] &{}&\mathcal{A}
	\end{tikzcd}\begin{tikzcd}[row sep = 15, font = \fontsize{9}{6}]
		\mathcal{A} \otimes \mathcal{A} \arrow[rr, "\bigoplus"] \arrow[dddd, "T"']&{}\arrow[dddd, Rightarrow, shorten = 5, "\chi"]&\mathcal{A}\arrow[dddd, "T"']\arrow[ddr, "T"]
		\\
		&&{}\arrow[dd, Rightarrow, shift left = 5, "\mu"]
		\\
		&&&\mathcal{A}\arrow[ddl, "T"]
		\\
		&&{}
		\\
		\mathcal{A} \otimes \mathcal{A} \arrow[rr, "\bigoplus"'] &{}&\mathcal{A}
	\end{tikzcd}$$
	\item The monoidal unit consists of a $2$-functor $I: \mathbf{1} \rightarrow \mathcal{A}$, a pseudonatural transformation $\iota: TI \rightarrow I$ and two invertible modifications $\eta_{0}$ and $\mu_{0}$ as depicted below.
	
	$$\begin{tikzcd}[row sep = 15,font = \fontsize{9}{6}]
		\mathbf{1}\arrow[rr, "I"] \arrow[dd, bend right=30, "!"']
		&{}\arrow[dd, Rightarrow, shorten = 5, "\iota"]
		& \mathcal{A}
		\arrow[dd, bend right =30, "T"'name=B]
		\arrow[dd, bend left = 30, "1_{\mathcal{A}}"name=A]
		\arrow[from=A, to=B, Rightarrow, shorten =5, "\eta"]
		\\
		&&&&\cong_{\eta_{0}}&&{}
		\\
		\mathbf{1}\arrow[rr, "I"'] &{}& \mathcal{A}
	\end{tikzcd}\begin{tikzcd}[row sep = 15, font = \fontsize{9}{6}]
		\mathbf{1}\arrow[rr, "I"] \arrow[dd, bend right=30, "!"']
		\arrow[dd, bend left = 30, "!"]
		&& \mathcal{A}
		\arrow[dd, bend left = 30, "1_{\mathcal{A}}"]
		\\
		=&=
		\\
		\mathbf{1}\arrow[rr, "I"'] &{}& \mathcal{A}
	\end{tikzcd}$$

	$$\begin{tikzcd}[row sep = 15, font = \fontsize{9}{6}]
		\mathbf{1} \arrow[rr, "I"] \arrow[dddd, "!"']\arrow[ddr, "!"]&{}\arrow[dd, Rightarrow, shorten = 5, shift left = 3, "\iota"]&\mathcal{A}\arrow[ddr, "T"]
		\\
		\\
		{}\arrow[r, equal, shorten = 5]&\mathbf{1} \arrow[rr, "I"] \arrow[ddl, "!"]\arrow[dd, Rightarrow, shorten = 5, shift left = 3, "\iota"]&&\mathcal{A}\arrow[ddl, "T"]&\cong_{\mu_{0}}&{}
		\\
		{}
		\\
		\mathbf{1} \arrow[rr, "I"'] &{}&\mathcal{A}
	\end{tikzcd}\begin{tikzcd}[row sep = 15, font = \fontsize{9}{6}]
		\mathbf{1} \arrow[rr, "I"] \arrow[dddd, "!"']&{}\arrow[dddd, Rightarrow, shorten = 5, "\iota"]&\mathcal{A}\arrow[dddd, "T"']\arrow[ddr, "T"]
		\\
		&&{}\arrow[dd, Rightarrow, shift left = 5, "\mu"]
		\\
		&&&\mathcal{A}\arrow[ddl, "T"]
		\\
		&&{}
		\\
		\mathbf{1}\arrow[rr, "I"'] &{}&\mathcal{A}
	\end{tikzcd}$$
	
	\item As an ambistrict $2$-cell, the associator consists of \begin{itemize}
		\item the assertion that the $2$-functor $\bigoplus: \mathcal{A} \otimes \mathcal{A} \rightarrow \mathcal{A}$ is a strictly associative product.
		\item an invertible modification $\omega$ as depicted below.
	\end{itemize}
\end{itemize}

$$
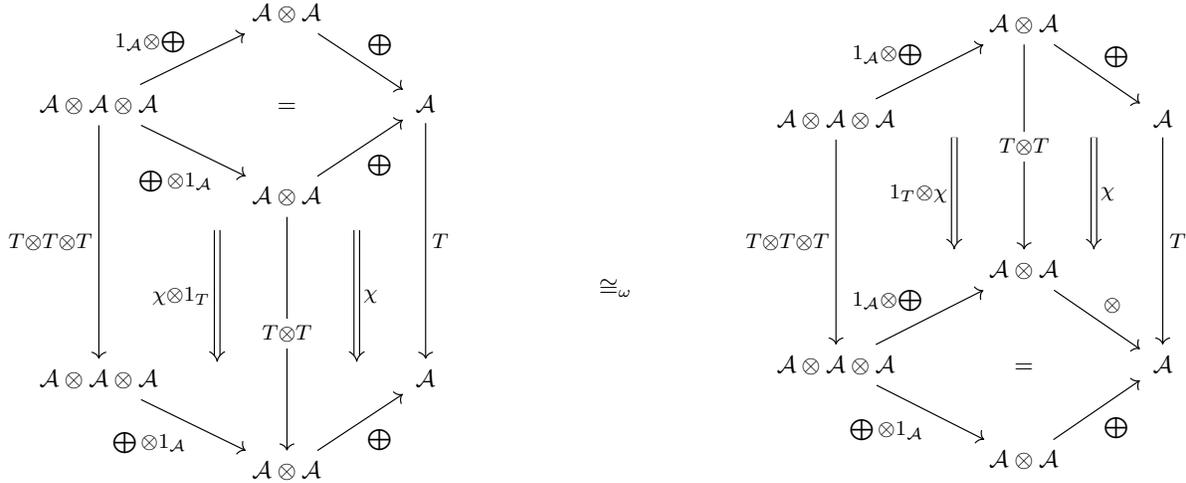
\begin{tikzcd}[font = \fontsize{9}{6}]
	&\mathcal{A} \otimes \mathcal{A}\arrow[rd, "\bigoplus"]
	\\
	\mathcal{A}\otimes\mathcal{A}\otimes\mathcal{A}
	\arrow[ddd, "T\otimes T\otimes T"']
	\arrow[ru, "1_\mathcal{A}\otimes\bigoplus"]
	\arrow[rd, near end, "\bigoplus\otimes 1_{\mathcal{A}}"']
	&=&\mathcal{A}\arrow[ddd, "T"]
	\\
	&\mathcal{A}\otimes \mathcal{A}\arrow[ru, "\bigoplus"']
	\arrow[dd, Rightarrow, shorten = 5, shift right = 10, "\chi \otimes 1_{T}"']
	\arrow[dd, Rightarrow, shorten = 5, shift left = 10, "\chi"]
	\arrow[ddd, "T\otimes T" description]
	\\
	&&&&\cong_{\omega}&{}
	\\
	\mathcal{A}\otimes \mathcal{A} \otimes \mathcal{A}\arrow[rd, "\bigoplus \otimes 1_\mathcal{A}"']&{}&\mathcal{A}
	\\
	&\mathcal{A}\otimes \mathcal{A}\arrow[ru, "\bigoplus"']
\end{tikzcd}\begin{tikzcd}[font = \fontsize{9}{6},row sep = 21]
	&\mathcal{A} \otimes \mathcal{A}\arrow[rd, "\bigoplus"]
	\arrow[ddd, "T\otimes T"description]
	\\
	\mathcal{A}\otimes\mathcal{A}\otimes\mathcal{A}
	\arrow[ddd, "T\otimes T\otimes T"']
	\arrow[ru, "1_\mathcal{A}\otimes\bigoplus"]
	&{}\arrow[dd, Rightarrow, shift right = 10, "1_{T}\otimes \chi"']\arrow[dd, Rightarrow, shift left = 10, "\chi"]&\mathcal{A}\arrow[ddd, "T"]
	\\
	\\
	&\mathcal{A} \otimes \mathcal{A}\arrow[rd, "\otimes"]
	\\
	\mathcal{A}\otimes \mathcal{A} \otimes \mathcal{A}
	\arrow[rd, "\bigoplus \otimes 1_\mathcal{A}"']
	\arrow[ru, "1_\mathcal{A}\otimes \bigoplus"]
	&=&\mathcal{A}
	\\
	&\mathcal{A}\otimes \mathcal{A}\arrow[ru, "\bigoplus"']
\end{tikzcd}$$

\begin{itemize}
	\item As an ambistrict $2$-cell, the left unitor consists of \begin{itemize}
		\item the assertion that the left unit law $\bigoplus.\left(I\otimes 1_{\mathcal{A}}\right) =1_{\mathcal{A}}$ holds on the nose.
		\item an invertible modification $\gamma$ as depicted below left.
	\end{itemize}
\end{itemize}

$$
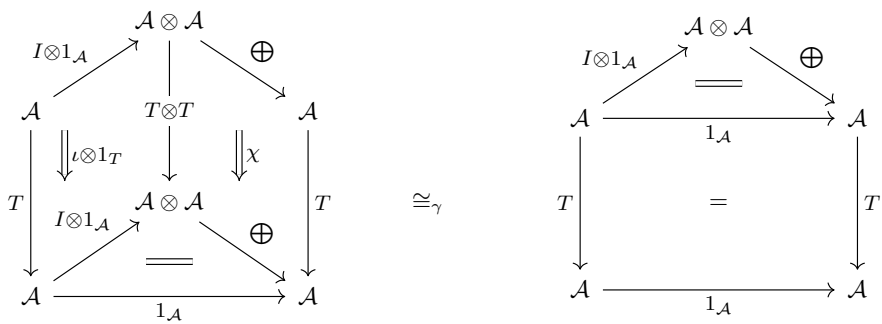
\begin{tikzcd}[font = \fontsize{9}{6}]
	&\mathcal{A}\otimes \mathcal{A}
	\arrow[rd, "\bigoplus"]
	\arrow[dd, "T\otimes T"description]
	\\
	\mathcal{A} \arrow[dd, "T"']\arrow[ru, "I \otimes 1_\mathcal{A}"]
	&{}\arrow[d, shift left = 10, Rightarrow, "\chi"]\arrow[d, shift right = 15, Rightarrow, "\iota \otimes 1_{T}"]
	& \mathcal{A}
	\arrow[dd, "T"]
	\\
	&\mathcal{A}\otimes \mathcal{A}\arrow[rd, "\bigoplus"]
	&&\cong_{\gamma}&{}
	\\
	\mathcal{A}\arrow[ru, near end, "I\otimes 1_{\mathcal{A}}"]\arrow[rr, "1_\mathcal{A}"']
	\arrow[rr, equal, shorten = 35, shift left = 5]
	&&\mathcal{A}
\end{tikzcd}\begin{tikzcd}[font = \fontsize{9}{6}]
	&\mathcal{A}\otimes \mathcal{A}
	\arrow[rd, "\bigoplus"]
	\\
	\mathcal{A} \arrow[dd, "T"']\arrow[ru, "I \otimes 1_\mathcal{A}"]
	\arrow[rr, "1_{\mathcal{A}}"']
	\arrow[rr, equal, shorten = 35, shift left = 5]
	&{}& \mathcal{A}
	\arrow[dd, "T"]
	\\
	&=
	\\
	\mathcal{A}\arrow[rr, "1_\mathcal{A}"']
	&&\mathcal{A}
\end{tikzcd}$$

\begin{itemize}
	\item As an ambistrict $2$-cell, the right unitor consists of \begin{itemize}
		\item the assertion that the right unit law $\bigoplus.\left(1_{\mathcal{A}}\otimes I\right) =1_{\mathcal{A}}$ holds on the nose.
		\item an invertible modification $\delta$ as depicted below right.
	\end{itemize}
\end{itemize}

$$\begin{tikzcd}[font = \fontsize{9}{6}]
	&\mathcal{A}\otimes \mathcal{A}
	\arrow[rd, "\bigoplus"]
	\arrow[dd, "T\otimes T"description]
	\\
	\mathcal{A} \arrow[dd, "T"']\arrow[ru, "1_\mathcal{A} \otimes I"]
	&{}\arrow[d, shift left = 10, Rightarrow, "\chi"]\arrow[d, shift right = 15, Rightarrow, "1_{T} \otimes \iota"]
	& \mathcal{A}
	\arrow[dd, "T"]
	\\
	&\mathcal{A}\otimes \mathcal{A}\arrow[rd, "\bigoplus"]
	&&\cong_{\delta}&{}
	\\
	\mathcal{A}\arrow[ru, near end, "1_\mathcal{A}\otimes I"]\arrow[rr, "1_\mathcal{A}"']
	\arrow[rr, equal, shorten = 35, shift left = 5]
	&&\mathcal{A}
\end{tikzcd}\begin{tikzcd}[font = \fontsize{9}{6}]
	&\mathcal{A}\otimes \mathcal{A}
	\arrow[rd, "\bigoplus"]
	\\
	\mathcal{A} \arrow[dd, "T"']\arrow[ru, "1_\mathcal{A} \otimes I"]
	\arrow[rr, "1_{\mathcal{A}}"']
	\arrow[rr, equal, shorten = 35, shift left = 5]
	&{}& \mathcal{A}
	\arrow[dd, "T"]
	\\
	&=
	\\
	\mathcal{A}\arrow[rr, "1_\mathcal{A}"']
	&&\mathcal{A}
\end{tikzcd}$$

\noindent AXIOMS
\\
\begin{itemize}
	\item The pentagon and unit laws for the pseudomonoid in $\left(\mathbf{Psmnd}\left(\mathbf{Gray}\right)_\text{ast}, \otimes\right)$ involve equations between the invertible $3$-cells in $\mathbf{Gray}$ which are built out of ambistrict $2$-cells in $\mathbf{Psmnd}\left(\mathbf{Gray}\right)$. These equations correspond to the pentagon and unit laws for $\left(T, \chi, \omega, \gamma, \delta\right)$ as an opmonoidal $2$-functor.
\end{itemize}

\begin{itemize}
	\item The unit and multiplication laws for the ambistrict $2$-cell $\omega$ correspond to the associativity laws for the monoidal pseudonatural transformations $\left(\eta, \eta_{0}, \eta_{-, ?}\right)$ and $\left(\mu, \mu_{0}, \mu_{-, ?}\right)$, respectively.
\end{itemize}

\begin{itemize}
	\item The unit and multiplication laws for the ambistrict $2$-cell $\gamma$ correspond to the left unit laws for the monoidal pseudonatural transformations $\left(\eta, \eta_{0}, \eta_{-, ?}\right)$ and $\left(\mu, \mu_{0}, \mu_{-, ?}\right)$ respectively.
	\item The unit and multiplication laws for the ambistrict $2$-cell $\delta$ correspond to the right unit laws for the monoidal pseudonatural transformations $\left(\eta, \eta_{0}, \eta_{-, ?}\right)$ and $\left(\mu, \mu_{0}, \mu_{-, ?}\right)$ respectively.
	\item As a morphism of pseudomonads, $\left(\bigoplus, \chi, \eta_{-, ?}, \mu_{-, ?}\right)$ must satisfy an associativity law, a left unit law, and a right unit law. These respectively correspond to the composition conditions needed for $\mathbf{a}$, $\mathbf{l}$ and $\mathbf{r}$ to be monoidal modifications.
	\item As a morphism of pseudomonads, $\left(I, \iota, \eta_{0}, \mu_{0}\right)$ must satisfy an associativity law, a left unit law, and a right unit law. These respectively correspond to the unit conditions needed for $\mathbf{a}$, $\mathbf{l}$ and $\mathbf{r}$ to be monoidal modifications.
\end{itemize}

\begin{theorem}\label{theorem semi-strict opmonoidal pseudomonads are pseudomonoids in the Gray monoid of pseudomonads}
	The correspondence described above is a bijection between the class of semi-strict opmonoidal pseudomonads, and the class of pseudomonoids in $(\mathbf{Psmnd}(\Gr)_\text{ast}, \otimes. \mathbf{1})$.
\end{theorem}

\begin{proof}
\noindent For brevity we give only a few of the equalities of pasting diagrams involved in the identification between pseudomonoids in the $\mathbf{Gray}$-monoid $\left(\mathbf{Psmnd}\left(\mathbf{Gray}\right)_\text{ast}, \otimes \right)$ and pseudomonads in the $\mathbf{Gray}$-category $\mathbf{opGray}$-$\mathbf{Monoids}$. The most complex part of the pseudomonoid structure is the pentagon law. This corresponds to the equation of pasting diagrams displayed below, which can in turn be seen as the pentagon law for the opmonoidal $2$-functor $(S, \iota, \chi, \gamma, \omega, \delta): (\A, \otimes, I) \to (\A, \otimes, I)$.

\begin{tikzcd}[font=\fontsize{6}{3}]
	T\left(W\bigoplus X \bigoplus Y \bigoplus Z\right)
	\arrow[rr, "\chi_{{}_{W\bigoplus X\bigoplus Y{,}Z}}"]
	\arrow[dddd, "\chi_{{}_{W{,}X\bigoplus Y \bigoplus Z}}"']
	&{}\arrow[dddd, Rightarrow, shorten = 30, "\omega_{W{,}X\bigoplus Y{,}Z}"]
	& T\left(W\bigoplus X \bigoplus Y\right)\bigoplus TZ
	\arrow[rdd, "\chi_{{}_{W\bigoplus X{,}Y}}\bigoplus 1"]
	\arrow[dddd, "\chi_{{}_{W{,}X\bigoplus Y}\bigoplus 1}"description]
	\\
	\\
	&&& T\left(W\bigoplus X\right)\bigoplus TY \bigoplus TZ
	\arrow[dddd, "\chi_{{}_{W{,}X}}\bigoplus 1\bigoplus 1"]
	\arrow[dd, Rightarrow, shorten = 12, shift right = 15, "\omega_{W{,}X{,}Y}\bigoplus 1"']
	\\
	\\
	TW\bigoplus T\left(X\bigoplus Y \bigoplus Z\right)
	\arrow[rdd, "1\bigoplus \chi_{{}_{X{,}Y\bigoplus Z}}"']
	\arrow[rr, "1\bigoplus \chi_{{}_{X\bigoplus Y{,}Z}}"]
	&{}\arrow[dd, Rightarrow, shorten =10, "1\bigoplus \omega_{X{,}Y{,}Z}"]
	&TW\bigoplus T\left(X\bigoplus Y\right) \bigoplus TZ
	\arrow[rdd, "1\bigoplus \chi_{{}_{X{,}Y}}\bigoplus 1"description]
	&{}
	\\
	\\
	&TW\bigoplus TX\bigoplus T\left(Y\bigoplus Z\right)
	\arrow[rr, "1\bigoplus 1\bigoplus \chi_{{}_{Y{,}Z}}"']
	&& TW\bigoplus TX\bigoplus TY\bigoplus TZ 
\end{tikzcd}

=\begin{tikzcd}[font=\fontsize{6}{3}]
	T\left(W\bigoplus X \bigoplus Y \bigoplus Z\right)
	\arrow[dddd, "\chi_{{}_{W{,}X\bigoplus Y \bigoplus Z}}"']
	\arrow[rdd, "\chi_{{}_{W\bigoplus X{,}Y\bigoplus Z}}"description] \arrow[rr, "\chi_{{}_{W\bigoplus X \bigoplus Y{,}Z}}"]
	&{}\arrow[dd, Rightarrow, shorten = 10, "\omega_{W\bigoplus X{,}Y{,}Z}"]
	& T\left(W\bigoplus X \bigoplus Y\right)\bigoplus TZ
	\arrow[rdd, "\chi_{{}_{W\bigoplus X{,}Y}}\bigoplus 1"]
	\\
	\\
	&T\left(W\bigoplus X\right)\bigoplus T\left(Y \bigoplus Z\right)
	\arrow[dd, Rightarrow, shorten = 12, shift right = 20, "\omega_{W{,}X{,}Y\bigoplus Z}"']
	\arrow[dddd, "\chi_{{}_{W{,}X}}\bigoplus 1"description]
	\arrow[rr, "1\bigoplus \chi_{{}_{Y{,}Z}}"]
	&{}\arrow[dddd, Rightarrow, shorten = 30, "{\left(\chi_{{}_{W{,}X}}\right)}_{\left(\chi_{{}_{Y{,}Z}}\right)}"]
	& T\left(W\bigoplus X\right)\bigoplus TY \bigoplus TZ
	\arrow[dddd, "\chi_{{}_{W{,}X}}\bigoplus 1\bigoplus 1"]
	\\
	\\
	TW\bigoplus T\left(X\bigoplus Y \bigoplus Z\right)
	\arrow[rdd, "1\bigoplus \chi_{{}_{X{,}Y\bigoplus Z}}"']
	&{}
	\\
	\\
	&TW\bigoplus TX\bigoplus T\left(Y\bigoplus Z\right)
	\arrow[rr, "1\bigoplus 1\bigoplus \chi_{{}_{Y{,}Z}}"']
	&{}& TW\bigoplus TX\bigoplus TY\bigoplus TZ 
\end{tikzcd}

\noindent We also give as an example the axioms needed for $(\eta, \eta_{0}, \eta_{2}): 1_\A \Rightarrow (S, \chi, \iota, \gamma, \omega, \delta)$ to be well-defined as an opmonoidal pseudonatural transformation. Note that these are indeed also respectively the unit axioms required for $\gamma$, $\omega$ and $\delta$ to be ambistrict pseudomonad transformations.
\\
\\
\emph{Left unit}

$$\begin{tikzcd}[column sep = 18,font = \fontsize{7}{4}]
	I\bigoplus X
	\arrow[rdd,"\eta_{I}\bigoplus 1"]
	\arrow[rdddd,shift right = 1, bend right = 30, "1_{I\bigoplus X}"']
	\arrow[rr, "\eta_{I\bigoplus X}"]
	&{}&S(I\bigoplus X)\arrow[rdd, "\chi_{I{,}X}"]
	\\
	&{}
	\\
	&SI\bigoplus X
	\arrow[uu, Rightarrow,shift right = 8, shorten = 15, "\eta_{2{,}I{,}X}"']
	\arrow[rr, "1\bigoplus \eta_{X}"]\arrow[dd, "\iota\bigoplus 1"]
	&{}&SI\bigoplus SX
	\arrow[dd, "\iota \bigoplus 1"]
	&=
	\\
	&{}\arrow[uu, shorten = 15, Rightarrow,shift left = 8,  "\eta_{0}\bigoplus 1"]
	\\
	&I\bigoplus X
	\arrow[rr, "1\bigoplus \eta_{X}"']
	&{}\arrow[uu, Rightarrow, shorten = 15, "{(\iota)}_{(\eta_{X})}"]
	&I\bigoplus SX
\end{tikzcd}\begin{tikzcd}[font = \fontsize{9}{6}, column sep = 15]
	I\bigoplus X
	\arrow[rdddd, bend right = 30, "1_{I\bigoplus X}"']
	\arrow[rr, "\eta_{I\bigoplus X}"]
	&{}&S(I\bigoplus X)\arrow[rdd, "\chi_{I{,}X}"]
	\\
	&&&{}
	\\
	&=
	&&
	SI\bigoplus SX\arrow[dd, "\iota \bigoplus 1"]
	\\
	&&&{}\arrow[uu, Rightarrow, shorten = 15, shift left = 10, "\gamma_{X}^{-1}"]
	\\
	&I\bigoplus X
	\arrow[rr, "1\bigoplus \eta_{X}"']
	&{}
	&I\bigoplus SX
	\arrow[luuuu, bend left = 30,leftarrow, "1_{SX}"]
\end{tikzcd}$$

\noindent \emph{Associativity}

$$\begin{tikzcd}[font = \fontsize{9}{6}]
	&X\bigoplus Y\bigoplus Z
	\arrow[rr, "\eta_{X\bigoplus Y \bigoplus Z}"]
	\arrow[ldd, "\eta_{X}\bigoplus 1\bigoplus 1"']
	&&S(X\bigoplus Y\bigoplus Z)
	\arrow[rdd, "\chi_{X\bigoplus Y{,}Z}"]
	\arrow[ldd, "\chi_{X{,}Y\bigoplus Z}"']
	\\
	\\
	SX\bigoplus Y\bigoplus Z\arrow[rr, "1\bigoplus \eta_{Y\bigoplus Z}"]
	\arrow[rdd, "1\bigoplus \eta_{Y}\bigoplus 1"']
	&{}\arrow[uu, Rightarrow, shorten = 15, "\eta_{2{,}X{,}Y\bigoplus Z}"']
	&SX\bigoplus S(Y\bigoplus Z)
	\arrow[rdd, "1\bigoplus \chi_{Y{,}Z}"']
	&&S(X\bigoplus Y)\bigoplus SZ
	\arrow[ldd, "\chi_{X{,}Y}\bigoplus 1"]
	\\
	\\
	&SX\bigoplus SY \bigoplus Z
	\arrow[uu, Rightarrow, shorten = 15,"1\bigoplus \eta_{2{,}Y{,}Z}"']
	\arrow[rr, "1\bigoplus 1\bigoplus \eta_{Z}"']
	&&SX\bigoplus SY \bigoplus SZ
	\arrow[uuuu, Rightarrow, shorten = 25, "\omega_{X{,}Y{,}Z}"]
\end{tikzcd}$$

=$$\begin{tikzcd}[font = \fontsize{9}{6}]
	&X\bigoplus Y\bigoplus Z
	\arrow[dddd, Leftarrow, shorten = 25, "\eta_{2{,}X{,}Y}\bigoplus 1"]
	\arrow[rdd, "\eta_{X\bigoplus Y}\bigoplus 1"]
	\arrow[rr, "\eta_{X\bigoplus Y \bigoplus Z}"]
	\arrow[ldd, "\eta_{X}\bigoplus 1\bigoplus 1"']
	&&S(X\bigoplus Y\bigoplus Z)
	\arrow[dd, Leftarrow, shorten = 15, "\eta_{2{,}X\bigoplus Y{,}Z}"]
	\arrow[rdd, "\chi_{X\bigoplus Y{,}Z}"]
	\\
	\\
	SX\bigoplus Y\bigoplus Z
	\arrow[rdd, "1\bigoplus \eta_{Y}\bigoplus 1"']
	&{}
	& S(X\bigoplus Y)\bigoplus Z
	\arrow[dd, Leftarrow, shorten = 15, shift left = 5, "{(\chi_{X{,}Y})}_{(\eta_{Z})}"]
	\arrow[ldd, "\chi_{X{,}Y}\bigoplus 1"]
	\arrow[rr, "1\bigoplus \eta_{Z}"']
	&{}&S(X\bigoplus Y)\bigoplus SZ
	\arrow[ldd, "\chi_{X{,}Y}\bigoplus 1"]
	\\
	\\
	&SX\bigoplus SY \bigoplus Z
	\arrow[rr, "1\bigoplus 1\bigoplus \eta_{Z}"']
	&{}&SX\bigoplus SY \bigoplus SZ
\end{tikzcd}$$

\noindent \emph{Right unit}

$$\begin{tikzcd}[font = \fontsize{9}{6}]
	X\bigoplus I \arrow[rr, "\eta_{X\bigoplus I}"]
	\arrow[dd, "\eta_{X}\bigoplus 1"']
	&{}\arrow[dd, Rightarrow, shorten = 15, "\eta_{2{,}X{,}I}"]
	& S(X\bigoplus I)\arrow[dd, "\chi_{X{,}I}"]
	\arrow[rdddd, bend left = 30, "1_{I}"]
	\\
	&&&{}\arrow[dd,shift right = 8,  Rightarrow, shorten = 15, "\delta_{X}"']
	\\
	SX\bigoplus I\arrow[rrrdd, bend right = 20, "1_{SX}"']
	\arrow[rr, "1\bigoplus \eta_{I}"]
	&{}& SX\bigoplus SI
	\arrow[dd, Rightarrow, shorten = 15,shift right = 10, "1\bigoplus \eta_{0}"]
	\arrow[ddr, "1\bigoplus \iota"']&&&=&1_{\eta_{X}}
	\\
	&&&{}
	\\
	&
	&{}&SX \bigoplus I
\end{tikzcd}$$
	
\end{proof}

\noindent Later, in Example \ref{Example pseudoalgebra structure on the unit}, we will also give the equalities of pasting diagrams needed for $(I, \iota, \eta_{0}, \mu_{0}): \mathbf{1} \to (\A, S)$ to be well-defined as a morphism of pseudomonads, which also correspond to the unit laws for the monoidal pseudonatural transformations $(\eta, \eta_{0}, \eta_{2})$ and $(\mu, \mu_{0}, \mu_{2})$.

\section{Interaction between the Eilenberg-Moore construction and $\mathbf{Gray}$-tensor products}\label{Section Interaction between the Eilenberg-Moore construction and Gray tensor products}

\noindent There is a $\mathbf{Gray}$-functor $\mathbf{EM}: \mathbf{Psmnd}\left(\mathbf{Gray}\right) \rightarrow \mathbf{Gray}$ which sends a pseudomonad $\left(\mathcal{A}, S\right)$ to the forgetful $2$-functor from the $2$-categories of pseudoalgebras, pseudomorphisms and pseudoalgebra $2$-cells, which we denote by $\mathcal{A}^S$. The aim of this section is to study the interaction of this $\Gr$-functor with $\Gr$ tensor products. We find that ambistrict pseudomonad transformations are sent to $2$-natural transformations, whose components are pseudomorphisms of a particularly simple form to be described in Corollary \ref{EM of ambistrict is 2-natural} part (4). We also examine how $\mathbf{EM}$ respects the $\mathbf{Gray}$-monoid structure of Theorem \ref{Theorem Gray tensor of pseudomonads}. In particular, respect for the tensor product is mediated by a family of $2$-functors which we describe in detail in Remark \ref{Remark explicit description of compositor for EM}. This family of $2$-functors is natural as the pseudomonads vary, up to invertible $2$-natural transformations that will be described in Proposition \ref{Proposition 2-cell component of compositor as the second pseudomonad varies}. Moreover, we find in Proposition \ref{associativity and unitality of EM} that this family of $2$-functors satisfies associativity and unit conditions on the nose.
\\
\\
\noindent Most of the proofs in this Section are via inspection of Remarks \ref{Remark explicit EM} and \ref{Remark explicit description of compositor for EM}. The necessary calculations use relations in the presentation of the $\mathbf{Gray}$-tensor product, as described in Remark \ref{Remark Gray tensor product}, rather than axioms for pseudomonads or their higher dimensional cells. There are subtleties involved in packaging the calculations in this Section into a rigorous statement about normal lax monoidal $2$-functoriality of $\mathbf{EM}$. These subtleties will be addressed in Section \ref{Section lifting monoidal structures}.

\subsection{Ambistrict pseudomonad transformations and $2$-natural transformations with tight components}\label{Subsection the explicit description of EM}

\noindent Remark \ref{Remark explicit EM} recalls the explicit description of the $\Gr$-functor $\mathbf{EM}: \mathbf{Psmnd}\left(\mathbf{Gray}\right) \rightarrow \mathbf{Gray}$ constructed in \cite{Formal Theory of Pseudomonads}, and Corollary \ref{EM of ambistrict is 2-natural} relates special properties of pseudomonad transformations with properties of their images under $\mathbf{EM}$.

\begin{remark}\label{Remark explicit EM}
	On $1$-cells, $\mathbf{EM}$ sends a morphism of pseudomonads $\left(G, g, g_{0}, g_{2}\right): \left(\mathcal{A}, S\right) \rightarrow \left(\mathcal{B}, T\right)$ to the $2$-functor $\overline{G}: \mathcal{A}^{S} \rightarrow \mathcal{B}^{T}$ which 
	
	\begin{itemize}
		\item maps a pseudoalgebra $\left(X, x: SX \rightarrow X, x_{0}: x.\eta_{X} \Rightarrow 1_{X}, x_{2}: x.Tx \Rightarrow x.\mu_{X}\right)$ to the pseudoalgebra whose structure map is as displayed below
		
		$$\begin{tikzcd}[font=\fontsize{9}{6}]
			TGX \arrow[rr, "g_{X}"] && GSX \arrow[rr, "Gx"] && GX{,}
		\end{tikzcd}$$
		\noindent and whose unit and multiplication constraints are given in Equation \ref{Equation unit and multiplication constraints of G X}.
		\begin{equation}\label{Equation unit and multiplication constraints of G X}
			\begin{tikzcd}[font=\fontsize{9}{6}]
				&{}\arrow[dd, Rightarrow, shorten = 10, "g_{0{,}X}"]
				&TGX\arrow[dd, "g_{X}"]
				\\
				\\
				GX\arrow[rr, "G\eta_{X}"]\arrow[rruu, bend left = 30, "\eta_{GX}"] \arrow[rrdd, bend right = 30, "1_{GX}"']
				&{}\arrow[dd, Rightarrow, shorten = 10, "Gx_{0}"]
				& GSX\arrow[dd, "Gx"]
				\\
				\\
				&{}&GX&{}
			\end{tikzcd}\begin{tikzcd}[font=\fontsize{9}{6}]
				T^{2}GX \arrow[rr, "Tg_{X}"] \arrow[dddd, "\mu_{GX}"']
				&{}\arrow[dddd, Rightarrow, shorten = 30, "g_{2{,}X}"]
				& TGSX \arrow[rr, "TGx"]\arrow[dd, "g_{SX}"description] &{}\arrow[dd, Rightarrow, shorten = 15, "g_{x}"]& TGX\arrow[dd, "g_{X}"]
				\\
				\\
				&& GS^{2}X\arrow[rr, "GSx"]\arrow[dd, "G\mu_{X}"description] &{}\arrow[dd, Rightarrow, shorten = 15, "Gx_{2}"]
				& GSX\arrow[dd, "Gx"]
				\\
				\\
				TGX \arrow[rr, "g_{X}"'] &{}& GSX \arrow[rr, "Gx"'] &{}& GX
			\end{tikzcd}
		\end{equation}
		\item maps a pseudomorphism $\left(p: X \rightarrow Y, \overline{p}: y.Sp \Rightarrow p.x\right)$ for $\left(\mathcal{A}, S\right)$ to the pseudomorphism whose underlying $1$-cell is $Gp: GX \rightarrow GY$ and whose $2$-cell component is given by the pasting depicted to the left of Equation \ref{Equation 2-cell components of pseudomorphisms G p and phi X}.
		\item Sends the $2$-cell of $S$-pseudoalgebras $\phi: \left(p, \overline{p}\right) \rightarrow \left(q, \overline{q}\right)$ to the $2$-cell of $T$-pseudoalgebras determined by $G\phi: Gp \Rightarrow Gq$.
	\end{itemize}
	
	\noindent On $2$-cells, $\mathbf{EM}$ maps the $2$-cell of pseudomonads $\left(\phi, \tilde{\phi}\right): \left(G, g, g_{0}, g_{2}\right) \Rightarrow \left(H, h, h_{0}, h_{2}\right)$ to the pseudonatural transformation whose $1$-cell component on $\left(X, x, x_{0}, x_{2}\right)$ is the pseudomorphism $\left(GX, G\phi.g_{X}\right) \rightarrow \left(HX, H\phi.g_{X}\right)$ which has $1$-cell component given by $\phi_{X}: GX \rightarrow HX$ and $2$-cell component given by the pasting in $\mathcal{B}$ depicted to the right of Equation \ref{Equation 2-cell components of pseudomorphisms G p and phi X}. Meanwhile, the $2$-cell component of $\mathbf{EM}\left(\phi, \tilde{\phi}\right)$ on a pseudomorphism $\left(p, \overline{p}\right): \left(X, x, x_{0}, x_{2}\right) \rightarrow \left(Y, y, y_{0}, y_{2}\right)$ is given by the $2$-cell component in $\B$ of the pseudonatural transformation $\phi: G \Rightarrow H$ on the morphism $p: X \to Y$ in $\A$.
	
\begin{equation}\label{Equation 2-cell components of pseudomorphisms G p and phi X}
		\begin{tikzcd}[font=\fontsize{9}{6}]
		TGX \arrow[rr, "TGp"]\arrow[dd, "g_{X}"'] &{}\arrow[dd, Rightarrow, shorten = 10, "g_{p}"]& TGY\arrow[dd, "g_{Y}"]
		\\
		\\
		GSX \arrow[rr, "GSp"]
		\arrow[dd, "Gx"']
		&{}\arrow[dd, Rightarrow, shorten = 10, "G\overline{p}"]
		& GSY\arrow[dd, "Gy"]
		\\
		\\
		GX \arrow[rr, "Gp"'] &{}& GY&{}
	\end{tikzcd}\begin{tikzcd}[font=\fontsize{9}{6}]
		TGX \arrow[rr, "T\phi_{X}"]\arrow[dd, "g_{X}"'] &{}\arrow[dd, Rightarrow, shorten = 10, "\tilde{\phi}_{X}"]
		& THX
		\arrow[dd, "h_{X}"]
		\\
		\\
		GSX \arrow[rr, "\phi_{SX}"]
		\arrow[dd, "Gx"']
		&{}\arrow[dd, Rightarrow, shorten = 10, "\phi_{x}"]
		& HSX\arrow[dd, "Hx"]
		\\
		\\
		GX \arrow[rr, "\phi_{X}"']
		&{}& HX
	\end{tikzcd}
\end{equation}
\end{remark}

\noindent In Definition \ref{definition tight pseudomorphism}, we isolate a property shared by many of the pseudomorphisms we will consider.

\begin{definition}\label{definition tight pseudomorphism}
	A pseudomorphism will be called \emph{tight} if its $1$-cell component is the identity.
\end{definition}

\begin{corollary}\label{EM of ambistrict is 2-natural}
	Let $(\phi, \tilde{\phi}): (G, g, g_{0}, g_{2}) \to (H, h, h_{0}, h_{2})$ be a pseudomonad transformation.
	
	\begin{enumerate}
		\item The $2$-cell components of the pseudonatural transformation $\mathbf{EM}(\phi, \tilde{\phi})$ on tight pseudomorphisms, are given by identities.
		\item The pseudonatural transformation $\mathbf{EM}(\phi, \tilde{\phi})$ is $2$-natural if and only if $\phi$ is $2$-natural.
		\item The $1$-cell components of $\mathbf{EM}(\phi, \tilde{\phi})$ are tight if and only if $\phi$ is an invertible icon. 
		\item The pseudomonad transformation $(\phi, \tilde{\phi})$ is ambistrict if and only if $\mathbf{EM}(\phi, \tilde{\phi})$ is $2$-natural and has tight components.
	\end{enumerate}
\end{corollary}

\begin{proof}
	Part (1) follows from the unit law for the pseudonatural transformation $\phi: G \Rightarrow H$. Part (4) is a consequence of parts (2) and (3), since by Proposition \ref{Proposition characterising ambistrict pseudomonad transformations}, $\left(\phi, \tilde{\phi}\right)$ is ambistrict precisely when $\phi: G \rightarrow H$ is an identity pseudonatural transformation. But parts (2) and (3) follow by specialising the description given in Remark \ref{Remark explicit EM}.
\end{proof}

\noindent In Example \ref{Example pseudoalgebra structure on the unit}, to follow, we explicitly describe the image under $\mathbf{EM}$ of a particularly simple morphism of pseudomonads. The result is a $2$-functor of the form $\mathbf{1} \to \A^S$. As we will discuss in Subsection \ref{subsection monoidal structure on pseudoalgebras}, this pseudoalgebra will be the unit in a monoidal structure on the $2$-category of pseudoalgebras $\mathcal{A}^S$. An explicit description of the rest of the monoidal structure will be given in Theorem \ref{monoidal bicategory structure on pseudoalgebras}.

\begin{example}\label{Example pseudoalgebra structure on the unit}
	If $(\mathcal{A}, \bigoplus, I)$ is a $\mathbf{Gray}$-monoid equipped with a semi-strict opmonoidal pseudomonad $(S, \eta, \mu)$, then as observed in Subsection \ref{Section monoidal pseudomonads as pseudomonoids in pseudomonads}, there is a morphism of pseudomonads of the form $(I, \iota, \eta_{0}, \mu_{0}): (\mathbf{1}, 1_\mathbf{1}) \to (\A, S)$. Under $\mathbf{EM}$, this gives rise to a $2$-functor $\mathbf{1} \to \A^{S}$. This $2$-functor selects a pseudoalgebra for $(\A, S)$ on the unit $I$ for the $\mathbf{Gray}$-monoid structure. The structure map for this pseudoalgebra is given by $\iota: SI \to I$, while its unit and multiplication constraints are given as displayed below. 
	
	$$\begin{tikzcd}[font=\fontsize{9}{6}]
		I \arrow[rrdd, bend right = 30, "1_{I}"']
		\arrow[rr, "\eta_{I}"]
		&{}\arrow[dd, Leftarrow, shorten = 10, "\eta_{0}"]
		& SI\arrow[dd, "\iota"]
		\\
		\\
		&{}&I&{}
	\end{tikzcd}\begin{tikzcd}[font=\fontsize{9}{6}]
	S^{2}I\arrow[dd, "\mu_{I}"']
	\arrow[rr, "S\iota"]&{}\arrow[dd, Rightarrow, shorten = 10, "\mu_{0}"]
	&SI\arrow[dd, "\iota"]
	\\
	\\
	SI\arrow[rr, "\iota"']&{}&I
\end{tikzcd}$$

	\noindent The axioms required for these data to comprise a well-defined pseudoalgebra are as depicted below. Observe that these are precisely the axioms required for $(I, \iota, \eta_{0}, \mu_{0})$ to be well-defined as a morphism of pseudomonads, and also exactly the unit axioms needed for $\mathbf{a}$, $\mathbf{l}$ and $\mathbf{r}$ to be well-defined as monoidal modifications.

	$$\begin{tikzcd}[font=\fontsize{9}{6}]
		S^{3}I \arrow[dd, "\mu_{SI}"']
		\arrow[rr, "S^{2}\iota"]
		\arrow[rd, "S\mu_{I}"description]
		&{}\arrow[d, Rightarrow, shorten = 5, "S\mu_{0}"]& S^{2}X\arrow[rd, "S\iota"]
		\\
		&S^{2}I\arrow[rr, "S\iota"]\arrow[dd, "\mu_{I}"]
		\arrow[d, Rightarrow, shorten = 5, shift right = 12, "\mathbf{a}_{I}"]
		&{}\arrow[dd, Rightarrow, shorten = 15, "\mu_{0}"]
		&SX\arrow[dd, "\iota"]
		\\
		S^{2}I\arrow[rd, "\mu_{I}"']&{}&&&=
		\\
		&SI \arrow[rr, "\iota"']
		&{}& I
	\end{tikzcd}\begin{tikzcd}[font=\fontsize{9}{6}]
		S^{3}I \arrow[dd, "\mu_{SI}"']
		\arrow[rr, "S^{2}\iota"]
		&{}\arrow[dd, Rightarrow, shorten = 20, Rightarrow, "\mu_{\iota}"]
		& S^{2}I\arrow[rd, "S\iota"]\arrow[dd, "\mu_{I}"description]
		\\
		&&&SX\arrow[dd, "\iota"]
		\arrow[d, Rightarrow, shorten = 5, shift right = 12, "\mu_{0}"]
		\\
		S^{2}I\arrow[rd, "\mu_{I}"']\arrow[rr, "S\iota"]
		&{}\arrow[d, Rightarrow,"\mu_{0}"]
		&SI\arrow[rd, "\iota"]
		&{}
		\\
		&SI \arrow[rr, "\iota"']
		&{}& I
	\end{tikzcd}$$
	
	$$\begin{tikzcd}[font=\fontsize{9}{6}]
		&&{}\arrow[dd, Rightarrow, shorten = 10, "S\eta_{0}"]
		&&SI\arrow[rrdd, "\iota"]\arrow[dddd, Rightarrow, shorten = 15, "\mu_{0}"]
		\\
		\\
		SI \arrow[rrrruu, bend left = 30, "1_{SI}"]
		\arrow[rrrrdd, bend right = 30, "1_{SI}"']
		\arrow[rr, "S\eta_{I}"]
		&& S^{2}I \arrow[rruu, "S\iota"]\arrow[rrdd, "\mu_{I}"]\arrow[dd, Rightarrow, shorten = 10, "\mathbf{r}_{I}"]
		&&&&I
		&=
		&SI\arrow[rr, bend left = 20, "\iota"name=A]\arrow[rr, bend right = 20, "\iota"'name=B] && I
		\\
		\\
		&&{}&&SI\arrow[rruu, "\iota"]
		\arrow[from=A, to=B, Rightarrow, shorten = 5, "1_{\iota}"]
	\end{tikzcd}$$

	$$\begin{tikzcd}[font = \fontsize{9}{6}]
		SI \arrow[rr, "\iota"]\arrow[dd, "\eta_{SI}"'] &{}\arrow[dd, Rightarrow, shorten = 10, "\eta_{\iota}"]
		& I\arrow[dd, "\eta_{I}"]\arrow[rdddd, bend left = 30, "1_{I}"]
		\\
		&&&{}\arrow[dd,shift right = 8,  Rightarrow, shorten = 15, "\eta_{0}"]
		\\
		S^2I\arrow[rdd, "\mu_{I}"']\arrow[rr, "S\iota"]
		&{}\arrow[dd, Rightarrow, shorten = 10, "\mu_{0}"]
		& SI\arrow[ddr, "\iota"']&&&=&{}
		\\
		&&&{}
		\\
		&SI\arrow[rr, "\iota"']
		&&I
	\end{tikzcd}\begin{tikzcd}[font = \fontsize{9}{6}]
	SI\arrow[dddd,shift left = 8, Rightarrow, shorten = 50, "\mathbf{l}_{I}"]
	\arrow[rdddd, bend left = 30, "1_{SI}"] \arrow[rr, "\iota"]\arrow[dd, "\eta_{SI}"'] &{}
	& I\arrow[rdddd, bend left = 30, "1_{I}"]
	\\
	\\
	S^2I\arrow[rdd, "\mu_{I}"']&&=
	&{}&
	\\
	&&&{}
	\\
	{}&SI\arrow[rr, "\iota"']
	&&I
\end{tikzcd}$$
	
\end{example}

\begin{notation}
	Whenever it is possible to do so without ambiguity, the image of a morphism of pseudomonads $\left(G, g, g_{0}, g_{2}\right): \left(\mathcal{A}, S\right) \rightarrow \left(\mathcal{A}', S'\right)$ under $\mathbf{EM}$ will be denoted as $\overline{G}: \mathcal{A}^{S} \rightarrow \mathcal{A'}^{S'}$.
\end{notation}

\subsection{Towards normal lax monoidality of $\mathbf{EM}$}\label{Subsection normal lax monoidality of EM}

\begin{remark}\label{Remark explicit description of compositor for EM}
	Let $\left(\mathcal{A}, S\right)$ and $\left(\mathcal{B}, T\right)$ be pseudomonads in $\mathbf{Gray}$. There is a canonical $2$-functor 
	
	$$M_{S, T}: \mathcal{A}^{S}\otimes \mathcal{B}^{T} \rightarrow {\left(\mathcal{A}\otimes \mathcal{B}\right)}^{\left(S\otimes T\right)}$$ 
	
	\noindent induced by the universal property of the Eilenberg-Moore object by the pseudoalgebra 
	
	$$\left(U^{S} \otimes U^{T}, \varepsilon^{S}\otimes \varepsilon^{T}, u^{S} \otimes u^{T}, m^{S} \otimes m^{T}\right)$$
	
	\noindent for the pseudomonad $\mathbf{Gray}\left(\mathcal{A}^{S}\otimes \mathcal{B}^{T}, S\otimes T\right)$ on the $2$-category $\mathbf{Gray}\left(\mathcal{A}^{S}\otimes \mathcal{B}^{T}, \mathcal{A}\otimes \mathcal{B}\right)$. We omit proof that this pseudoalgebra is well-defined as it is similar to the proof of the pseudomonad axioms in Lemma \ref{Gray tensor of two pseudomonads}. We give a more elementary description of this $2$-functor. On objects, it sends a pair of pseudoalgebras $\left(\left(X, x, x_{0}, x_{2}\right), \left(Y, y, y_{0}, y_{2}\right)\right)$ to the pseudoalgebra with structure map given as displayed below
	
	$$\begin{tikzcd}[font = \fontsize{9}{6}]
		\left(SX{,}TY\right) \arrow[rr, "\left(x{,}1\right)"] && \left(X{,}TY\right) \arrow[rr, "\left(1{,}y\right)"] && \left(X{,}Y\right) 
	\end{tikzcd}$$

\noindent and with unit and composition constraints given as depicted in Equation \ref{Equation unit and composition constraints of M X Y}.
	\begin{equation}\label{Equation unit and composition constraints of M X Y}
	\begin{tikzcd}[font = \fontsize{9}{6}, column sep = 16]
		\left(X{,}Y\right)\arrow[rr, "\left(\eta_{X}{,}1\right)"]
		\arrow[rrdd, bend right = 30, "1"']
		&{}\arrow[dd, shorten = 15, Rightarrow, "\left(x_{0}{,}1\right)"]
		& \left(SX{,}Y\right)\arrow[rr, "\left(1{,}\eta_{Y}\right)"]
		\arrow[dd, "\left(x{,}1\right)"]
		&{}\arrow[dd, Rightarrow, shorten = 15, "{(x)}_{\left(\eta_{Y}\right)}"]& \left(SX{,}TY\right)
		\arrow[dd, "\left(x{,}1\right)"]
		\\
		\\
		&{}& \left(X{,}Y\right)
		\arrow[rr, "\left(1{,}\eta_{Y}\right)"]
		\arrow[rrdd, bend right =30, "1"']
		&{}\arrow[dd, Rightarrow, shorten = 15, "\left(1{,}y_{0}\right)"]
		&\left(X{,}TY\right)
		\arrow[dd, "\left(1{,}y\right)"]
		\\
		\\
		&&&{}&\left(X{,}Y\right) &{}
	\end{tikzcd}\begin{tikzcd}[font=\fontsize{9}{6}, column sep = 16]
		\left(S^{2}X{,}T^{2}Y\right)\arrow[rr, "\left(Sx{,}1\right)"]
		\arrow[dd, "\left(\mu_{X}{,}1\right)"']
		&{}\arrow[dd, Rightarrow,shift right = 7,  shorten = 15, "\left(x_{2}{,}1\right)"]
		&\left(SX{,}T^{2}Y\right)
		\arrow[rr, "\left(1{,}Ty\right)"]
		\arrow[dd, "\left(x{,}1\right)"description]
		&{}\arrow[dd, Rightarrow, shorten = 15, "{\left(x\right)}_{(Ty)}"]
		&\left(SX{,}TY\right)
		\arrow[dd, "\left(x{,}1\right)"]
		\\
		\\
		\left(SX{,}T^{2}Y\right)
		\arrow[dd, "\left(1{,}\mu_{Y}\right)"']
		\arrow[rr, "\left(x{,}1\right)"]
		&{}\arrow[dd, Rightarrow, shorten = 15, "{\left(x\right)}_{(\mu_{Y})}"]
		&\left(X{,}T^{2}Y\right)
		\arrow[dd, "\left(1{,}\mu_{Y}\right)" description]
		\arrow[rr, "\left(1{,}Ty\right)"]
		&{}\arrow[dd, Rightarrow, shorten = 15, "\left(1{,}y_{2}\right)"]
		&\left(X{,}TY\right)
		\arrow[dd, "\left(1{,}y\right)"]
		\\
		\\
		\left(SX{,}TY\right)
		\arrow[rr, "\left(x{,}1\right)"']
		&{}&\left(X{,}TY\right)
		\arrow[rr, "\left(1{,}y\right)"']
		&{}&\left(X{,}Y\right)
	\end{tikzcd}
\end{equation}
	\noindent We now describe $M_{S, T}: \mathcal{A}^{S}\otimes\mathcal{B}^{T} \rightarrow \left(\mathcal{A}\otimes \mathcal{B}\right)^{S \otimes T}$ on morphisms. Let $\left(\mathbf{p}: \mathbf{X} \rightarrow \mathbf{X'}\right) := \left(p, \overline{p}\right): \left(X, x, x_{0}, x_{2}\right) \rightarrow \left(X', x', x_{0}', x_{2}'\right)$ be a pseudomorphism for $\left(\mathcal{A}, {S}\right)$, and let $\mathbf{Y}:= \left(Y, y, y_{0}, y_{2}\right)$ be a pseudoalgebra for $\left(\mathcal{B}, T\right)$. Then $M$ sends the generating morphism $\left(\mathbf{p}, \mathbf{Y}\right): \left(\mathbf{X}, \mathbf{Y}\right) \rightarrow \left(\mathbf{X}', \mathbf{Y}\right)$ to the pseudomorphism for $\left(\mathcal{A}\otimes\mathcal{B}, S\otimes T\right)$ whose underlying $1$-cell is $\left(p, Y\right): \left(X, Y\right) \rightarrow \left(X', Y\right)$ and whose $2$-cell component is given by the pasting depicted to the left of Equation \ref{Equation M of p and q}. Similarly, if $\left(\mathbf{q}: \mathbf{Y} \rightarrow \mathbf{Y'}\right):= \left(q, \overline{q}\right): \left(Y, y, y_{0}, y_{2}\right) \rightarrow \left(Y', y', y_{0}', y_{2}'\right)$ is a pseudomorphism for $\left(\mathcal{B}, T\right)$, then $M: \mathcal{A}^{S}\otimes\mathcal{B}^{T} \rightarrow \left(\mathcal{A}\otimes \mathcal{B}\right)^{S \otimes T}$ sends the generating morphism $\left(\mathbf{X}, \mathbf{q}\right): \left(\mathbf{X}, \mathbf{Y}\right) \rightarrow \left(\mathbf{X}, \mathbf{Y'}\right)$ to the pseudomorphism for $\left(\mathcal{A}\otimes\mathcal{B}, S\otimes T\right)$ whose underlying $1$-cell is $\left(X, q\right): \left(X,Y\right) \rightarrow \left(X, Y'\right)$ and whose $2$-cell component is given by the pasting depicted to the right of Equation \ref{Equation M of p and q}. 

	\begin{equation}\label{Equation M of p and q}
		\begin{tikzcd}[font=\fontsize{9}{6}]
			\left(SX{,}TY\right)\arrow[rr, "\left(Sp{,}1\right)"]
			\arrow[dd, "\left(x{,}1\right)"']
			&{}\arrow[dd, Rightarrow, shorten = 15, "\left(\overline{p}{,}1\right)"]& \left(SX'{,}TY\right)
			\arrow[dd, "\left(x'{,}1\right)"]
			\\
			\\
			\left(X{,}TY\right)\arrow[rr, "\left(p{,}1\right)"] \arrow[dd, "\left(1{,}y\right)"']
			&{}\arrow[dd, Rightarrow, shorten = 15,"{(p)}_{(y)}"]& \left(X'{,}TY\right)\arrow[dd, "\left(1{,}y\right)"]
			\\
			\\
			\left(X{,}Y\right)
			\arrow[rr, "\left(p{,}1\right)"'] &{}& \left(X'{,}Y\right)&{}
		\end{tikzcd}\begin{tikzcd}[font=\fontsize{9}{6}]
			\left(SX{,}TY\right)\arrow[rr, "\left(1{,}Tq\right)"]
			\arrow[dd, "\left(x{,}1\right)"'] &{}\arrow[dd, Rightarrow, shorten = 15, "{\left(x\right)}_{(Tq)}"]& \left(SX{,}TY'\right)\arrow[dd, "\left(x{,}1\right)"]
			\\
			\\
			\left(X{,}TY\right)\arrow[rr, "\left(1{,}Tq\right)"] \arrow[dd, "\left(1{,}y\right)"']&{}\arrow[dd, Rightarrow, shorten = 15, "\left(1{,}\overline{q}\right)"]& \left(X{,}Y'\right)\arrow[dd, "\left(1{,}y'\right)"]
			\\
			\\
			\left(X{,}Y\right)\arrow[rr, "\left(1{,}q\right)"'] &{}& \left(X{,}Y'\right)
		\end{tikzcd}
	\end{equation}
	\noindent Finally, $M: \mathcal{A}^{S}\otimes\mathcal{B}^{T} \rightarrow \left(\mathcal{A}\otimes \mathcal{B}\right)^{S \otimes T}$ extends to $2$-cells from each component in the evident way, and sends interchangers in $\mathcal{A}^{S}\otimes \mathcal{B}^T$ of $\mathbf{p}: \mathbf{X} \to \mathbf{X}'$ in $\A^S$ and $\mathbf{q}: \mathbf{Y} \to \mathbf{Y}'$ in $\B^T$ to the $2$-cell in ${(\mathcal{A}\otimes \B)}^{(S\otimes T)}$ determined by the interchanger in $\A \otimes \B$ of $p\in \A(X, X')$ and $q \in \B(Y , Y')$. 
\end{remark}

\begin{notation}\label{notation pseudoalgebras and pseudomorphisms}
	We will use boldface to abbreviate notation for pseudoalgebras, and for pseudomorphisms, as in Remark \ref{Remark explicit description of compositor for EM}.
\end{notation}

\begin{lemma}\label{lemma tensors and interchangers of tight pseudomorphisms}
	Let $\mathbf{p}: \mathbf{X} \to \mathbf{X}'$ and $\mathbf{q}: \mathbf{Y} \to \mathbf{Y}'$ be pseudomorphisms.
	
	\begin{enumerate}
		\item The pseudomorphism $M_{S, T}(\mathbf{X}, \mathbf{q}): M_{S, T}(\mathbf{X}, \mathbf{Y}) \to M_{S, T}(\mathbf{X}, \mathbf{Y}')$ is tight if and only if $\mathbf{q}: \mathbf{Y} \to \mathbf{Y}'$ is tight.
		\item The pseudomorphism $M_{S, T}(\mathbf{p}, \mathbf{Y}): M_{S, T}(\mathbf{X}, \mathbf{Y}) \to M_{S, T}(\mathbf{X}', \mathbf{Y})$ is tight if and only if $\mathbf{p}: \mathbf{X} \to \mathbf{X}'$ is tight. 
		\item If either $\mathbf{p}: \mathbf{X} \to \mathbf{X}'$ or $\mathbf{q}: \mathbf{Y} \to \mathbf{Y'}$ are tight, then the image of their interchanger in $\A^S \otimes \B^T$ under the $2$-functor $M_{S, T}: \A^S \otimes \B^T \to {(\A \otimes \B)}^{(S\otimes T)}$ is an identity.
	\end{enumerate}
\end{lemma}

\begin{proof}
	This is an easy inspection of the final comments in Remark \ref{Remark explicit description of compositor for EM}.
\end{proof}

\begin{proposition}\label{Proposition M is natural in the first component}
	Let $\left(G, g, g_{0}, g_{2}\right):\left(\mathcal{A}, S\right) \rightarrow \left(\mathcal{A}', S'\right)$ be a morphism of pseudomonads. Then there is a commutative diagram of $2$-functors as depicted below.
	
	$$\begin{tikzcd}[font=\fontsize{9}{6}]
		\mathcal{A}^{S} \otimes \mathcal{B}^{T} \arrow[rr, "M_{S{,}T}"]
		\arrow[dd, "\overline{G}\otimes \mathcal{B}^{T}"'] && \left(\mathcal{A} \otimes \mathcal{B}\right)^{S\otimes T}
		\arrow[dd, "\overline{G\otimes \mathcal{B}}"]
		\\
		\\
		\mathcal{A'}^{S'} \otimes \mathcal{B}^{T} \arrow[rr, "M_{S'{,}T}"'] && \left(\mathcal{A}' \otimes \mathcal{B}\right)^{S'\otimes T}
	\end{tikzcd}$$
\end{proposition}

\begin{proof}
	\noindent The proof is via an inspection of Remark \ref{Remark explicit EM}, in which we gave an explicit description of $1$-cells and $2$-cells in the image of $\mathbf{EM}$, and Remark \ref{Remark explicit description of compositor for EM}, in which we gave an explicit description of the $2$-functors $M_{S, T}$. Specifically, when one unwinds the definitions of each pair of cells which need to be proved to be equal, one finds that they are already equal just by the pasting theorem for $2$-categories; the proof uses neither the axioms for the components involved, nor the relations in the presentation of the $\mathbf{Gray}$-tensor product. We give more detail below. 
	\\
	\\
	\noindent A pair of pseudoalgebras $\big(\left(X, x, x_{0},x_{2}\right), \left(Y, y, y_{0}, y_{2}\right)\big)$ in $\mathcal{A}^{S}\otimes \mathcal{B}^T$ is mapped to the pseudoalgebra for $\left(\mathcal{A}\otimes \mathcal{B}, S \otimes T\right)$ with structure map \begin{tikzcd}[font=\fontsize{9}{6}]
		\left(S'GX{,}Y\right) \arrow[rr, "\left(g_{X}{,}Y\right)"] && \left(GSX{,}Y\right) \arrow[rr, "\left(Gx{,}Y\right)"] && \left(GX{, Y}\right) 
	\end{tikzcd}. The unit and multiplication constraints for this pseudoalgebra are as depicted below (c.f. Equations \ref{Equation unit and multiplication constraints of G X} and \ref{Equation unit and composition constraints of M X Y}).
	
	$$\begin{tikzcd}[column sep = 16, font=\fontsize{9}{6}]
		&{}\arrow[dd, Rightarrow, shorten = 15, "\left(g_{0{,}X}{,}1\right)"]
		&\left(S'GX{,}Y\right) \arrow[rr, "\left(1{,}\eta_Y\right)"]\arrow[dd, "\left(g_{X}{,}1\right)"]
		&{}\arrow[dd, Rightarrow, shorten = 15, "{(g_{X})}_{\left(\eta_{Y}\right)}"]
		& \left(S'GX{, }TY\right)
		\arrow[dd, "\left(g_{X}{,}1\right)"]
		\\
		\\
		\left(GX{,}Y\right)
		\arrow[rruu, bend left = 30, "({\eta'}_{GX}{,}1)"]
		\arrow[rr, "\left(\eta_{GX}{,}1\right)"]
		\arrow[rrdd, bend right = 30, "1"']
		&{}\arrow[dd, shorten = 15, Rightarrow, "\left(Gx_{0}{,}1\right)"]& \left(GSX{,}Y\right)\arrow[rr, "\left(1{,}\eta_{Y}\right)"] \arrow[dd, "\left(Gx{,}1\right)"]
		&{}\arrow[dd, Rightarrow, shorten = 15, "{(Gx)}_{\left(\eta_{Y}\right)}"]& \left(GSX{,}TY\right)\arrow[dd, "\left(Gx{,}1\right)"]
		\\
		\\
		&{}& \left(GX{,}Y\right)
		\arrow[rr, "\left(1{,}\eta_{Y}\right)"]
		\arrow[rrdd, bend right =30, "1"']
		&{}\arrow[dd, Rightarrow, shorten = 15, "\left(1{,}y_{0}\right)"]
		&\left(GX{,}TY\right)
		\arrow[dd, "\left(1{,}y\right)"]
		\\
		\\
		&&&{}&\left(GX{,}Y\right)
	\end{tikzcd}$$
	
	$$\begin{tikzcd}[font=\fontsize{9}{6}]
		\left({S'}^{2}GX{,}T^{2}Y\right)\arrow[dddd, "\left({\mu'}_{GX}{,}1\right)"']
		\arrow[rr, "\left(S'g_{X}{,}1\right)"]
		&&\left(S'GSX{,}T^{2}Y\right)
		\arrow[dd, "\left(g_{SX}{,}1\right)"description]
		\arrow[rr, "\left(S'Gx{,}1\right)"]
		&{}\arrow[dd, Rightarrow, shorten = 15, shift right = 5, "\left(g_{x}{,}1\right)"]
		& \left(S'GX{,}T^{2}Y\right)\arrow[dd, "\left(g_{X}{,}1\right)"description]
		\arrow[rr, "\left(1{,}Ty\right)"]
		&{}\arrow[dd, Rightarrow, shorten = 15, "{(g_{X})}_{\left(Ty\right)}"]
		&\left(S'GX{,}TY\right)
		\arrow[dd, "\left(g_{X}{,}1\right)"]
		\\
		&{}\arrow[dd, Rightarrow, shorten = 5, "\left(g_{2{,}X}{,}1\right)"']
		\\
		{}&&\left(GS^{2}X{,}T^{2}Y\right)\arrow[rr, "\left(GSx{,}1\right)"]\arrow[dd, "\left(G\mu_{X}{,}1\right)"description]
		&{}\arrow[dd, Rightarrow,shift right = 7,  shorten = 15, "\left(Gx_{2}{,}1\right)"]
		&\left(GSX{,}T^{2}Y\right)\arrow[rr, "\left(1{,}Ty\right)"]\arrow[dd, "\left(Gx{,}1\right)"description]
		&{}\arrow[dd, Rightarrow, shorten = 15, "{\left(Gx\right)}_{(Ty)}"]
		&\left(GSX{,}TY\right)\arrow[dd, "\left(Gx{,}1\right)"]
		\\
		&{}
		\\
		\left(S'GX{,}T^{2}Y\right)\arrow[rr, "\left(g_{X}{,}1\right)"]\arrow[dd, "\left(1{,}\mu_{Y}\right)"']
		&{}\arrow[dd, Rightarrow, shorten = 15, "{(g_{X})}_{\left(\mu_{Y}\right)}"]
		&\left(GSX{,}T^{2}Y\right)
		\arrow[dd, "\left(1{,}\mu_{Y}\right)"description]\arrow[rr, "\left(Gx{,}1\right)"]
		&{}\arrow[dd, Rightarrow, shorten = 15, "{(g_{X})}_{\left(\mu_{Y}\right)}"]
		&\left(GX{,}T^{2}Y\right)\arrow[dd, "\left(1{,}\mu_{Y}\right)" description]\arrow[rr, "\left(1{,}Ty\right)"]
		&{}\arrow[dd, Rightarrow, shorten = 15, "\left(1{,}y_{2}\right)"]
		&\left(GX{,}TY\right)\arrow[dd, "\left(1{,}y\right)"]
		\\
		\\
		\left(S'GX{,}TY\right)\arrow[rr, "\left(g_{X}{,}1\right)"']
		&{}&\left(GSX{,}TY\right)\arrow[rr, "\left(Gx{,}1\right)"']
		&{}&\left(GX{,}TY\right)\arrow[rr, "\left(1{,}y\right)"']&{}&\left(GX{,}Y\right)
	\end{tikzcd}$$
	
	\noindent Let $\left(\mathbf{p}: \mathbf{X} \rightarrow \mathbf{X'}\right):= \left(p, \overline{p}\right): \left(X, x, x_{0}, x_{2}\right) \rightarrow \left(X', x', x_{0}, {x}_{2}'\right)$ and $\left(\mathbf{q}: \mathbf{Y} \rightarrow \mathbf{Y'}\right) :=\left(q, \overline{q}\right): \left(Y, y, y_{0}, y_{2}\right) \rightarrow \left(Y', y', y_{0}', y_{2}'\right)$ be pseudomorphisms in $\mathcal{A}^{S}$ and $\mathcal{B}^{T}$ respectively, and consider the generating morphisms $\left(\mathbf{p}, \mathbf{Y}\right): \left(\mathbf{X}, \mathbf{Y}\right) \rightarrow \left(\mathbf{X}', \mathbf{Y}\right)$ and $\left(\mathbf{X}, \mathbf{q}\right): \left(\mathbf{X}, \mathbf{Y}\right) \rightarrow \left(\mathbf{X}, \mathbf{Y}'\right)$ in the $\mathbf{Gray}$-tensor product $\mathcal{A}^{S} \otimes \mathcal{B}^{T}$. Observe that the image of $\left(\mathbf{p}, \mathbf{Y}\right)$ under both $\overline{G \otimes \mathcal{B}}\circ M_{S, T}$ and $M_{S', T} \circ \left(\overline{G}\otimes \mathcal{B}^{T}\right)$ is the pseudomorphism with $1$-cell component $\left(Gp, Y\right): \left(GX, Y\right) \rightarrow \left(GX', Y\right)$ and $2$-cell component given as depicted to the left of Equation \ref{Equation proof M is natural in ambistrict pseudomonad transformations}. Similarly, the image of $\left(\mathbf{X}, \mathbf{q}\right): \left(\mathbf{X}, \mathbf{Y}\right)$ under both $\overline{G \otimes \mathcal{B}}\circ M_{S, T}$ and $M_{S', T} \circ \left(\overline{G}\otimes \mathcal{B}^{T}\right)$ is the pseudomorphism with $1$-cell component $\left(GX, q\right): \left(GX, Y\right) \rightarrow \left(GX, Y'\right)$ and $2$-cell component given as depicted to the right of Equation \ref{Equation proof M is natural in ambistrict pseudomonad transformations} (c.f. Equations \ref{Equation 2-cell components of pseudomorphisms G p and phi X} and \ref{Equation M of p and q}).

	\begin{equation}\label{Equation proof M is natural in ambistrict pseudomonad transformations}
		\begin{tikzcd}[font=\fontsize{9}{6}]
		\left(S'GX{,}TY\right)\arrow[dd, "\left(g_{X}{,}1\right)"']
		\arrow[rr, "\left(S'Gp{,}1\right)"] &{}\arrow[dd, Rightarrow, shorten = 15, "\left(g_{p}{,}1\right)"]
		& \left(S'GX'{,}TY\right)
		\arrow[dd, "\left(g_{X}{,}1\right)"]
		\\
		\\
		\left(GSX{,}TY\right)
		\arrow[rr, "\left(GSp{,}1\right)"]
		\arrow[dd, "\left(Gx{,}1\right)"'] &{}\arrow[dd, Rightarrow, shorten = 15, "\left(G\overline{p}{,}1\right)"]& \left(GSX'{,}TY\right)
		\arrow[dd, "\left(Gx'{,}1\right)"]
		\\
		\\
		\left(GX{,}TY\right)
		\arrow[rr, "\left(Gp{,}1\right)"] \arrow[dd, "\left(1{,}y\right)"']&{}\arrow[dd, Rightarrow, shorten = 15, "{\left(Gp\right)}_{(y)}"]& \left(GX'{,}TY\right)\arrow[dd, "\left(1{,}y\right)"]
		\\
		\\
		\left(GX{,}Y\right)\arrow[rr, "\left(Gp{,}1\right)"'] &{}& \left(GX'{,}Y\right)&{}
	\end{tikzcd}\begin{tikzcd}[font=\fontsize{9}{6}]
		\left(S'GX{,}TY\right)
		\arrow[rr, "\left(1{,}Tq\right)"]
		\arrow[dd, "\left(g_{X}{,}1\right)"'] &{}\arrow[dd, Rightarrow, shorten = 15, "{\left(g_{X}\right)}_{(Tq)}"]
		& \left(S'GX{,}TY'\right)\arrow[dd, "\left(g_{X}{,}1\right)"]
		\\
		\\
		\left(GSX{,}TY\right)\arrow[rr, "\left(1{,}Tq\right)"]\arrow[dd, "\left(Gx{,}1\right)"'] &{}\arrow[dd, Rightarrow, shorten = 15, "{\left(Gx\right)}_{(Tq)}"]& \left(GSX{,}TY'\right)\arrow[dd, "\left(Gx{,}1\right)"]
		\\
		\\
		\left(GX{,}TY\right)\arrow[rr, "\left(1{,}Tq\right)"] \arrow[dd, "\left(1{,}y\right)"']&{}\arrow[dd, Rightarrow, shorten = 15, "\left(1{,}\overline{q}\right)"]& \left(GX{,}Y'\right)\arrow[dd, "\left(1{,}y'\right)"]
		\\
		\\
		\left(GX{,}Y\right)\arrow[rr, "\left(1{,}q\right)"'] &{}& \left(GX{,}Y'\right)
	\end{tikzcd}
\end{equation}
	
	\noindent Finally, commutativity on $2$-cells is an even easier inspection. This completes the proof.
\end{proof}

\begin{proposition}\label{Proposition 2-naturality of compositor for EM in the first variable}
	Let $\left(G, g, g_{0}, g_{2}\right):\left(\mathcal{A}, S\right) \rightarrow \left(\mathcal{A}', S'\right)$ be as above and let $\left(1_{G}, \tilde{\phi}\right): \left(G, g, g_{0}, g_{2}\right) \rightarrow \left(G', g', g_{0}', g_{2}'\right)$ be an ambistrict $2$-cell of pseudomonads. Then the whiskerings $M_{S', T}\circ\left(\overline{\phi}\otimes \B^T\right)$ and $\overline{\left((1_{G}, \tilde{\phi}) \otimes \mathcal{B}\right)}\circ M_{S, T}$ are $2$-natural transformations, and are moreover equal. 	
\end{proposition}

\begin{proof}
	We first explain why the whiskerings are $2$-natural. By Lemma \ref{Lemma - tensor B T 2-functor} part (2), since $\left(1_{G}, \tilde{\phi}\right)$ is an ambistrict pseudomonad transformation, so is $\left(\left(1_{G}, \tilde{\phi}\right) \otimes \mathcal{B}\right)$. Hence by Corollary \ref{EM of ambistrict is 2-natural}, its image under $\mathbf{EM}: \PsG \to \mathbf{Gray}$ is $2$-natural, and in particular its whiskering with $M_{S, T}$ is also $2$-natural. On the other hand, as discussed in Remark \ref{Remark subtleties in Gray tensor of 2-natural transformations}, despite the fact that $\overline{\phi}: \overline{G} \to \overline{G'}$ is $2$-natural, the tensor $\overline{\phi}\otimes \B$ is only pseudonatural, with $2$-cell components given by interchangers in $\A^S \otimes \B^T$ involving the $1$-cell components of $\overline{\phi}$ in $\A^S$. But by Corollary \ref{EM of ambistrict is 2-natural}, these components are tight pseudomorphisms. Thus by Lemma \ref{lemma tensors and interchangers of tight pseudomorphisms} part (3), the $2$-functor $M_{S', T}: \A^S \otimes \B^T \to {(\A \otimes \B)}^{(S\otimes T)}$ sends such interchangers to identities. Hence the whiskering $M_{S', T}\circ (\overline{\phi}\otimes \B^T)$ is also $2$-natural, as required.
	\\
	\\
	\noindent We now show that the two $2$-natural transformations are equal. By inspection, their component pseudomorphism on a pair of pseudoalgebras $\left(X, x, x_{0}, x_{2}\right)$ and $\left(Y, y, y_{0}, y_{2}\right)$ is the tight pseudomorphism for $\left(\mathcal{A}'\otimes \mathcal{B}, S'\otimes T\right)$ with $2$-cell component given as depicted below.
	
	$$\begin{tikzcd}[font =\fontsize{9}{6}]
		(S'GX{,}TY)\arrow[dd, "1"']\arrow[rr, "(g_{X}{,}1)"]
		&{}\arrow[dd, Rightarrow, shorten = 15, "(\tilde{\phi}_{X}{,}1)"]
		&(GSX{,}TY)\arrow[dd, "1"]\arrow[rr, "(Gx{,}1)"]
		&&(GX{,}TY)\arrow[dd, "1"]\arrow[rr, "(1{,}y)"]
		&&(GX{,}Y)\arrow[dd, "1"]
		\\
		\\
		(S'GX{,}TY)\arrow[rr, "(g_{X}'{,}1)"']
		&{}&(GSX{,}TY)\arrow[rr, "(Gx{,}1)"']
		&&(GX{,}TY)\arrow[rr, "(1{,}y)"']
		&&(GX{,}Y)
	\end{tikzcd}$$
\end{proof}

\noindent Propositions \ref{Proposition M is natural in the first component} and \ref{Proposition 2-naturality of compositor for EM in the first variable} will amount to $2$-naturality of the assignment $\left(\left(\mathcal{A}, S\right), \left(\mathcal{B}, T\right)\right) \mapsto M_{S, T}$ in the first variable. On the other hand, this assignment is only pseudonatural as the pseudomonad $\left(\mathcal{B}, T\right)$ is varied. The pseudonaturality constraint is a $2$-natural isomorphism whose components are also tight pseudomorphisms. It is described in Proposition \ref{Proposition 2-cell component of compositor as the second pseudomonad varies}, to follow.

\begin{proposition}\label{Proposition 2-cell component of compositor as the second pseudomonad varies}
	Let $\left(H, h, h_{0}, h_{2}\right): \left(\mathcal{B}, T\right) \rightarrow \left(\mathcal{B}', T'\right)$ be a morphism of pseudomonads and $(\A, S)$ be another pseudomonad.
	
	\begin{enumerate}
		\item There is a $2$-natural isomorphism $M_{1, H}$ as depicted below left, with component on an object $\left(\mathbf{X}, \mathbf{Y}\right)\in \mathcal{A}^{S}\otimes \mathcal{B}^T$ is given by the tight pseudomorphism whose $2$-cell component is depicted below right.
		
		$$\begin{tikzcd}[column sep = 18, font=\fontsize{9}{6}]
			\mathcal{A}^{S} \otimes \mathcal{B}^{T} \arrow[rr, "M_{S{,}T}"]
			\arrow[dd, "\mathcal{A}^{S}\otimes \overline{H}"'] &{}\arrow[dd, Rightarrow, shorten = 10, "M_{1{,}h}"]
			& \left(\mathcal{A} \otimes \mathcal{B}\right)^{S\otimes T}
			\arrow[dd, "\overline{\mathcal{A} \otimes H}"]
			\\
			\\
			\mathcal{A}^{S} \otimes \mathcal{B'}^{T'} \arrow[rr, "M_{S{,}T'}"'] &{}& \left(\mathcal{A} \otimes \mathcal{B}'\right)^{S\otimes T'}&{}
		\end{tikzcd}\begin{tikzcd}[column sep = 18, font=\fontsize{9}{6}]
			\left(SX{,}T'HY\right) \arrow[rr, "\left(1{,}h_{Y}\right)"]
			\arrow[dd, "\left(x{,}1\right)"']
			&{}\arrow[dd, Rightarrow, shorten = 15, "\left(x{,}h_{Y}\right)"]
			& \left(SX{,} HTY\right) \arrow[dd, "\left(x{,}1\right)"]
			\\
			\\
			\left(X{,}T'HY\right) \arrow[rr, "\left(1{,}h_{Y}\right)"']
			&{}& \left(X{,}HTY\right) \arrow[rr, "\left(1{,}Hy\right)"'] && \left(X{,}HY\right)
		\end{tikzcd}$$
		\item In particular, this $2$-natural isomorphism is the identity if $h: T'H \Rightarrow HT$ is an invertible icon.
		\item Suppose $(H', h', h_{0}, h_{2}): (\B', T') \rightarrow (\B'', T'')$ is another morphism of pseudomonads and let $M_{1, h'\circ h}$ denote the $2$-natural isomorphism given as in part (1) for the composite morphism of pseudomonads from $(\B, T)$ to $(\B'', T'')$. Then there is an equation of $2$-natural isomorphisms as depicted below.
		
		$$\begin{tikzcd}[column sep = 18, font=\fontsize{9}{6}]
			\mathcal{A}^{S} \otimes \mathcal{B}^{T} \arrow[rr, "M_{S{,}T}"]
			\arrow[dd, "\mathcal{A}^{S}\otimes \overline{H}"'] &{}\arrow[dd, Rightarrow, shorten = 10, "M_{1{,}h}"]
			& \left(\mathcal{A} \otimes \mathcal{B}\right)^{S\otimes T}
			\arrow[dd, "\overline{\mathcal{A} \otimes H}"]
			\\
			\\
			\mathcal{A}^{S} \otimes \mathcal{B'}^{T'}\arrow[dd, "\mathcal{A}^{S}\otimes \overline{H'}"'] \arrow[rr, "M_{S{,}T'}"'] &{}\arrow[dd, Rightarrow, shorten = 10, "M_{1{,}h'}"]& \left(\mathcal{A} \otimes \mathcal{B}'\right)^{S\otimes T'}\arrow[dd, "\overline{\mathcal{A} \otimes H'}"]&=&
			\\
			\\
			\mathcal{A}^{S} \otimes \mathcal{B''}^{T''}	\arrow[rr, "M_{S{,}T''}"'] &{}& \left(\mathcal{A} \otimes \mathcal{B}''\right)^{S\otimes T''}
		\end{tikzcd}\begin{tikzcd}[column sep = 18, font=\fontsize{9}{6}]
		\mathcal{A}^{S} \otimes \mathcal{B}^{T} \arrow[rr, "M_{S{,}T}"]
		\arrow[dddd, "\mathcal{A}^{S}\otimes \overline{H'\circ H}"'] &{}\arrow[dddd, Rightarrow, shorten = 20, "M_{1{,}h'\circ h}"]
		& \left(\mathcal{A} \otimes \mathcal{B}\right)^{S\otimes T}
		\arrow[dddd, "\overline{\mathcal{A} \otimes H'\circ H}"]
		\\
		\\
		\\
		\\
		\mathcal{A}^{S} \otimes \mathcal{B''}^{T''} \arrow[rr, "M_{S{,}T''}"'] &{}& \left(\mathcal{A} \otimes \mathcal{B}''\right)^{S\otimes T''}
	\end{tikzcd}$$
		
	\end{enumerate} 
\end{proposition}

\begin{proof}
	The proof of part (1) is given in Appendix \ref{Appendix M 1 h well-defined}. For part (2), it is clear from the description given that if this $2$-natural transformation is well-defined, it will be the identity if the $1$-cell components $h_{Y}: T'HY \to HTY$ are identities, and hence if $h$ is an invertible icon. For part (3), the component of both $2$-natural transformations are tight pseudomorphisms for the pseudomonad $(\A\otimes \B'', S\otimes T'')$, whose $2$-cell components are given as displayed below.
	
	$$\begin{tikzcd}[column sep = 18, font=\fontsize{9}{6}]
		\left(SX{,}T''H'HY\right) \arrow[rr, "\left(1{,}{h'}_{HY}\right)"]
		\arrow[dd, "\left(x{,}1\right)"']
		&{}\arrow[dd, Rightarrow, shorten = 15, "\left(x{,}{h'}_{HY}\right)"']&\left(SX{,}H'T'HY\right) \arrow[rr, "\left(1{,}H'h_{Y}\right)"]
		\arrow[dd, "\left(x{,}1\right)"'description]
		&{}\arrow[dd, Rightarrow, shorten = 15, "\left(x{,}H'h_{Y}\right)"]
		& \left(SX{,} H'HTY\right) \arrow[dd, "\left(x{,}1\right)"]
		\\
		\\
		\left(X{,}T''H'HY\right) \arrow[rr, "\left(1{,}{h'}_{HY}\right)"']&{}&\left(X{,}H'T'HY\right) \arrow[rr, "\left(1{,}H'h_{Y}\right)"']
		&{}& \left(X{,}H'HTY\right) \arrow[rr, "\left(1{,}H'Hy\right)"']&&\left(X{,}H'HY\right)
	\end{tikzcd}$$
\end{proof}

\begin{proposition}\label{Proposition naturality of compositor in 2-cells as second component varies}
	Let $\left(H, h, h_{0}, h_{2}\right):\left(\mathcal{B}, T\right) \rightarrow \left(\mathcal{B}', T'\right)$ be as above and let $\left(1_{H}, \tilde{\psi}\right): \left(H, h, h_{0},h_{2}\right) \rightarrow \left(H, h', {h'}_{0}, {h'}_{2}\right)$ be an ambistrict pseudomonad transformation. Then there is an equation between $2$-natural transformations as displayed below.
	
	$$\begin{tikzcd}[font =\fontsize{9}{6}]
		\mathcal{A}^{S}\otimes \mathcal{B}^{T} \arrow[dd, Rightarrow, shorten = 10, shift left = 5, "M_{1{,}h'}"]\arrow[rr, "M_{S{,}T}"]\arrow[dd, bend right = 60, "\mathcal{A}^{S}\otimes \overline{H'}"'] && {\left(\mathcal{A}\otimes\mathcal{B}\right)}^{S\otimes T}\arrow[dd, bend right = 60, "\overline{\left(\mathcal{A}\otimes H\right)}"'name=A]
		\arrow[dd, bend left = 60, "\overline{\left(\mathcal{A}\otimes H\right)}"name=B]
		\arrow[from=B, to=A, Rightarrow, shift left = 5, shorten = 5, "\overline{\mathcal{A}\otimes\psi}"]
		\\
		&&&=
		\\
		\mathcal{A}^{S}\otimes \mathcal{B'}^{T'} \arrow[rr, "M_{S{,}T'}"'] && {\left(\mathcal{A}\otimes\mathcal{B'}\right)}^{S\otimes T'}
	\end{tikzcd}\begin{tikzcd}[font=\fontsize{9}{6}]
		\mathcal{A}^{S}\otimes \mathcal{B}^{T}\arrow[rr, "M_{S{,}T}"]\arrow[dd, bend right = 60, "\mathcal{A}^{S}\otimes \overline{H'}"'name=A]\arrow[dd, bend left = 60, "\mathcal{A}^{S}\otimes \overline{H}"name=B] && {\left(\mathcal{A}\otimes\mathcal{B}\right)}^{S\otimes T}\arrow[dd, bend left = 60, "\overline{\left(\mathcal{A}\otimes H\right)}"]\arrow[dd, Rightarrow, shorten = 10, shift right = 8, "M_{1{,}h'}"]
		\arrow[from=B, to=A, Rightarrow, shift left = 5, shorten = 5, "\mathcal{A}^{S}\otimes\overline{\psi}"]
		\\
		\\
		\mathcal{A}^{S}\otimes \mathcal{B'}^{T'} \arrow[rr, "M_{S{,}T'}"'] && {\left(\mathcal{A}\otimes\mathcal{B'}\right)}^{S\otimes T'}
	\end{tikzcd}$$
	
\end{proposition}

\begin{proof}
	We first note that both sides are indeed $2$-natural transformations, since the $1$-cell components of the pseudonatural transformation $\A^S \otimes \overline{\psi}$ are tight pseudomorphisms, and $M_{S, T'}$ sends these to identities. This is similar to the comments in the first paragraph of the proof of Proposition \ref{Proposition 2-naturality of compositor for EM in the first variable}. The $2$-natural transformation on the left side of the equation has components on $(\mathbf{X}, \mathbf{Y})$ given by the pasting depicted below left, while the $2$-natural transformation on the right side of the equation has components on $(\mathbf{X}, \mathbf{Y})$ given by the pasting depicted below right. By the interchange law in $\A\otimes \B$ for the morphism $x: SX \to X$ in $\A$ on the $2$-cell $\overline{\psi}_{Y}$ in $\B$, these are equal.  

	$$\begin{tikzcd}[column sep = 18, font=\fontsize{9}{6}]
		\left(SX{,}T'HY\right) \arrow[rr,bend left= 30, "\left(1{,}h_{Y}\right)"name=A]
		 \arrow[rr,bend right = 30, "\left(1{,}h_{Y}\right)"'name=B]
		\arrow[dd, "\left(x{,}1\right)"']
		&{}\arrow[from=A, to=B, Rightarrow, shorten = 12,shift right = 3, "\left(x{,}\psi_{Y}\right)"]
		& \left(SX{,} HTY\right) \arrow[dd, "\left(x{,}1\right)"]
		\\
		\\
		\left(X{,}T'HY\right) \arrow[rr,bend right = 30, "\left(1{,}{h'}_{Y}\right)"'name=C]
		\arrow[from=B, to=C, Rightarrow, shorten = 8,shift right = 2, "\left(x{,}h'_{Y}\right)"]
		&{}& \left(X{,}HTY\right) \arrow[d, "\left(1{,}Hy\right)"]
		\\
		&& \left(X{,}HY\right)&{}
	\end{tikzcd}\begin{tikzcd}[column sep = 18, font=\fontsize{9}{6}]
	\left(SX{,}T'HY\right) \arrow[rr,bend left= 30, "\left(1{,}h_{Y}\right)"name=A]
	\arrow[dd, "\left(x{,}1\right)"']
	&{}\arrow[from=A, to=B, Rightarrow, shorten = 12, "\left(x{,}h_{Y}\right)"]
	& \left(SX{,} HTY\right) \arrow[dd, "\left(x{,}1\right)"]
	\\
	\\
	\left(X{,}T'HY\right) \arrow[rr,bend left = 30, "\left(1{,}h_{Y}\right)"name=B] \arrow[rr,bend right = 30, "\left(1{,}{h'}_{Y}\right)"'name=C]
	\arrow[from=B, to=C, Rightarrow, shorten = 8,shift right = 2, "\left(1{,}\psi_{Y}\right)"]
	&{}& \left(X{,}HTY\right) \arrow[d, "\left(1{,}Hy\right)"]
	\\
	&& \left(X{,}HY\right)
\end{tikzcd}$$
\end{proof}

\begin{proposition}\label{Proposition naturality of the compositor for EM on interchangers in Psmnd Gray otimes Psmnd Gray}
	\noindent Let $\left(G, g, g_{0}, g_{2}\right): \left(\mathcal{A}, S\right) \rightarrow \left(\mathcal{A}', S'\right)$ and $\left(H, h, h_{0}, h_{2}\right): \left(\mathcal{B}, T\right) \rightarrow \left(\mathcal{B}', T'\right)$ be morphisms of pseudomonads and consider the ambistrict pseudomonad transformation given by their interchange $2$-cell $(1_{HG}, h_{g})$ in the $\mathbf{Gray}$-monoid $\mathbf{Psmnd}(\mathbf{Gray})_\text{ast}$, as described in Lemma \ref{interchanger between morphisms of pseudomonads in Gray monoid structure on Psmnd Gray}. Let $\overline{h_{g}}$ denote $2$-natural transformation given as the image of this ambistrict pseudomonad transformation under $\mathbf{EM}$. Then there is an equation of $2$-natural isomorphisms as displayed below. 
	
	$$\begin{tikzcd}[font=\fontsize{9}{6}, column sep = 1, row sep = 10]
		&\A^S \otimes \B^T\arrow[ldd, "\overline{G}\otimes \B^T"']\arrow[rr, "M_{S{,}T}"] 
		&{}&{(\A \otimes \B)}^{(S\otimes T)}\arrow[ldd, "\overline{G\otimes \B}"]\arrow[rdd, "\overline{\A\otimes H}"]\arrow[dddd, Rightarrow, shorten = 10, "\overline{(1_{HG}{,}h_{g})}"]
		\\
		&&=
		\\
		{\A'}^{S'} \otimes \B^T
		\arrow[rdd, "{\A'}^{S'}\otimes \overline{H}"']
		\arrow[rr, "M_{S'{,}T}"] 
		&&{(\A' \otimes \B)}^{(S'\otimes T)}\arrow[rdd, "\overline{\A'\otimes H}"]\arrow[dd, Rightarrow, shorten = 10, "M_{S'{,}h}"]
		&&{(\A \otimes \B')}^{(S\otimes T')}\arrow[ldd, "\overline{G\otimes \B'}"]
		\\
		\\
		&{\A'}^{S'} \otimes {\B'}^{T'}	\arrow[rr, "M_{S'{,}T}"'] 
		&{}&{(\A' \otimes \B')}^{(S'\otimes T')}
	\end{tikzcd}$$

$$=\begin{tikzcd}[font=\fontsize{9}{6}, column sep = 1, row sep = 10]
	&\A^S \otimes \B^T
	\arrow[ldd, "\overline{G}\otimes \B^T"']
	\arrow[rdd, "\A^{S}\otimes \overline{H}"]
	\arrow[rr, "M_{S{,}T}"] 
	&{}&{(\A \otimes \B)}^{(S\otimes T)}\arrow[rdd, "\overline{\A\otimes H}"]\arrow[dd, Rightarrow, shorten = 10, "M_{S{,}h}"]
	\\
	\\
	{\A'}^{S'} \otimes \B^T
	\arrow[rdd, "{\A'}^{S'}\otimes \overline{H}"']
	&=&\A^{S}\otimes {\B'}^{T'}
	\arrow[rr,"M_{S{,}T'}"]
	\arrow[ldd, "\overline{G}\otimes {\B'}^{T'}"]
	&{}&{(\A \otimes \B')}^{(S\otimes T')}\arrow[ldd, "\overline{G\otimes \B'}"]
	\\
	&&&=
	\\
	&{\A'}^{S'} \otimes {\B'}^{T'}	\arrow[rr, "M_{S'{,}T}"'] 
	&{}&{(\A' \otimes \B')}^{(S'\otimes T')}
\end{tikzcd}$$
	
\end{proposition}

\begin{proof}
	Both $2$-natural isomorphisms have tight component pseudomorphisms on $(\mathbf{X}, \mathbf{Y})$, given as displayed below.

	$$\begin{tikzcd}[column sep = 18, font=\fontsize{9}{6}]
		(S'GX{,}T'HY)\arrow[dd, "(g_{X}{,}1)"']\arrow[rr, "(1{,}h_{Y})"]
		&{}\arrow[dd, Rightarrow, shorten = 15, "(g_{X}{,}h_{Y})"]
		&(S'GX{,}HTY)\arrow[dd, "(g_{X}{,}1)"]
		\\
		\\
		\left(GSX{,}T'HY\right) \arrow[rr, "\left(1{,}h_{Y}\right)"]
		\arrow[dd, "\left(Gx{,}1\right)"']
		&{}\arrow[dd, Rightarrow, shorten = 15, "\left(Gx{,}h_{Y}\right)"]
		& \left(GSX{,} HTY\right) \arrow[dd, "\left(Gx{,}1\right)"]
		\\
		\\
		\left(GX{,}T'HY\right) \arrow[rr, "\left(1{,}h_{Y}\right)"']
		&{}& \left(GX{,}HTY\right) \arrow[rr, "\left(1{,}Hy\right)"'] && \left(GX{,}HY\right)
	\end{tikzcd}$$
\end{proof}

\noindent We are finally ready to establish that the compositor for $\mathbf{EM}: {\mathbf{Psmnd}\left(\mathbf{Gray}\right)}_\text{ast} \rightarrow \mathbf{Gray}_\text{2-nat}$ is well-defined as a pseudonatural transformation.

\begin{corollary}\label{Corollary compositor for EM is a well-defined pseudonatural transformation}
	Consider the $2$-category ${\mathbf{Psmnd}\left(\mathbf{Gray}\right)}_\text{ast}$ as a $\mathbf{Gray}$ -category in which all $3$-cells are identities, and similar for the $\mathbf{Gray}$ tensor product ${\mathbf{Psmnd}\left(\mathbf{Gray}\right)}_\text{ast}\bigotimes \mathbf{Psmnd}\left(\mathbf{Gray}\right)_\text{ast}$. There is a trinatural transformation as depicted below.
	
	$$\begin{tikzcd}[font=\fontsize{9}{6}]
		{\mathbf{Psmnd}\left(\mathbf{Gray}\right)}_\text{ast}\bigotimes \mathbf{Psmnd}\left(\mathbf{Gray}\right)_\text{ast}
		\arrow[dd, "\mathbf{EM}\bigotimes \mathbf{EM}"'] \arrow[rr, "\otimes"] &{}\arrow[dd, Rightarrow, shorten = 10, "M_{-{,}?}"]& {\mathbf{Psmnd}\left(\mathbf{Gray}\right)}_\text{ast} \arrow[dd, "\mathbf{EM}"]
		\\
		\\
		\mathbf{Gray}^\mathbf{2} \bigotimes \mathbf{Gray}^\mathbf{2} \arrow[rr, "\otimes"']
		&{}& \mathbf{Gray}^\mathbf{2}
	\end{tikzcd}$$
	\noindent The image of all components of this trinatural transformation under the target projection $\mathbf{GRAY}$-functor $d_{0}:\Gr^\mathbf{2} \to \Gr$ are all identities. Moreover, the trinatural transformation $M_{-, ?}$ has
	
	\begin{itemize}
		\item $1$-cell components on a pair of pseudomonads $\left(\mathcal{A}, S\right)$, $\left(\mathcal{B}, T\right)$ given by the $2$-functors described in Remark \ref{Remark explicit description of compositor for EM},
		\item identity $2$-cell components on generators of the form $\left(\left(G, g, g_{0}, g_{2}\right), 1_{\left(\mathcal{B}, T\right)}\right): \left(\left(\mathcal{A}, S\right), \left(\mathcal{B}, T\right)\right) \rightarrow \left(\left(\mathcal{A}', S'\right), \left(\mathcal{B}, T\right)\right)$,
		\item $2$-cell components on generators of the form $\left(1_{\left(\mathcal{A}, S\right)}, \left(H, h, h_{0}, h_{2}\right)\right): \left(\left(\mathcal{A}, S\right), \left(\mathcal{B}, T\right)\right) \rightarrow \left(\left(\mathcal{A}, S\right), \left(\mathcal{B}', T'\right)\right)$ given by the $2$-natural isomorphisms described in Proposition \ref{Proposition 2-cell component of compositor as the second pseudomonad varies} part (1).
		\item identity $3$-cells mediating the usual axioms for a pseudonatural transformation.
	\end{itemize} 
\end{corollary}

\begin{proof}
	Strict naturality in pseudomonad morphisms from the first component has been checked in Proposition \ref{Proposition M is natural in the first component}, while $2$-naturality in ambistrict pseudomonad transformations from the first component is Proposition \ref{Proposition 2-naturality of compositor for EM in the first variable}. Naturality in ambistrict pseudomonad transformations from the second component is shown in Proposition \ref{Proposition naturality of compositor in 2-cells as second component varies}. Naturality in the interchanger $2$-cells in the $\Gr$-tensor product ${\mathbf{Psmnd}\left(\mathbf{Gray}\right)}_\text{ast} \bigotimes {\mathbf{Psmnd}\left(\mathbf{Gray}\right)}_\text{ast}$ is Proposition \ref{Proposition naturality of the compositor for EM on interchangers in Psmnd Gray otimes Psmnd Gray}. Proposition \ref{Proposition 2-cell component of compositor as the second pseudomonad varies} part (2) shows that the usual unit law for pseudonatural transformations holds on the nose, while Proposition \ref{Proposition 2-cell component of compositor as the second pseudomonad varies} part (3) shows that the usual composition law holds on composable pairs in which only the second pseudomonad is varied. Since the $2$-cell components for $M_{-, ?}$ on morphisms from the first argument are identities, as per Proposition \ref{Proposition M is natural in the first component}, it is easy to see that the usual composition law for pseudonatural transformations also holds strictly on all other composable pairs of morphisms. The trinaturality axioms also follow easily since all $3$-cell components for $M_{-, ?}$ are identities.
\end{proof}

\noindent Were $\Gr^\mathbf{2}$ to restrict to a $\Gr$-monoid whose $2$-cells were $2$-natural transformations, then the trinatural transformation of Corollary \ref{Corollary compositor for EM is a well-defined pseudonatural transformation} would be a $3$-pseudonatural transformation in the sense of Definition 3.2.3 of \cite{Miranda strictifying operational coherences}. We spend the rest of this subsection checking the remaining conditions involved normal lax monoidal $2$-functoriality of the Eilenberg-Moore construction for pseudomonads. A precise statement of this property for $\mathbf{EM}$ will be given later in Theorem \ref{Theorem EM is normal lax Gray functor}, after the subtleties involved in taking $\Gr$ tensor products of $2$-natural transformations are addressed in Subsection \ref{Subsection gray tensors of 2-naturals}.

\begin{proposition}\label{associativity and unitality of EM}
	Consider the family of $2$-functors of the form $M_{S, T}: \mathcal{A}^{S} \otimes \B^T \to {(\A\otimes \B)}^{(S \otimes T)}$ described in Remark \ref{Remark explicit description of compositor for EM}, and the family of $2$-natural isomorphisms of the form $M_{1, h}: \overline{\A\otimes H}\circ M_{S, T} \Rightarrow M_{S, T'}\circ \A^S \otimes \overline{H}$ described in Proposition \ref{Proposition 2-cell component of compositor as the second pseudomonad varies}.
	\begin{enumerate}
		\item If either $\mathcal{A}$ or $\mathcal{B}$ are terminal, then both $M_{S, T}$ and $M_{1, H}$ are identities.
		\item Let $(\mathcal{C}, R)$ be another pseudomonad. Then the following diagram of $2$-functors commutes.
		
		$$\begin{tikzcd}[font=\fontsize{9}{6}]
			\mathcal{A}^{S}\otimes \mathcal{B}^{T}\otimes \mathcal{C}^{R} \arrow[rr, "\mathcal{A}^{S}\otimes M_{T{,} R}"]\arrow[dd, "M_{S{,}T}\otimes \mathcal{C}^{R}"']
			&& 	\mathcal{A}^{S}\otimes {\left(\mathcal{B}\otimes \mathcal{C}\right)}^{\left(T\otimes R\right)}\arrow[dd, "M_{S{,}T\otimes R}"]
			\\
			\\
			{\left(\mathcal{A}\otimes \mathcal{B}\right)}^{\left(S\otimes T\right)}\otimes \mathcal{C}^{R}\arrow[rr, "M_{S\otimes T{,}R}"']
			&&{\left(\mathcal{A}\otimes \mathcal{B}\otimes \mathcal{C}\right)}^{\left(S\otimes T\otimes R\right)}
		\end{tikzcd}$$
	\end{enumerate}
	
\end{proposition}

\begin{proof}
	This follows from the unital and associative properties of the $\mathbf{Gray}$-tensor product.
\end{proof}

\section{Lifting monoidal structures to pseudoalgebras}\label{Section lifting monoidal structures}

\noindent In Section \ref{Section monoidal pseudomonads as pseudomonoids in pseudomonads} we discussed how opmonoidal pseudomonads could be viewed as pseudomonoids in an appropriate $\mathbf{Gray}$-monoid of pseudomonads, whose $2$-cells are ambistrict. One the one hand, this suggests that when we try to lift the $\mathbf{Gray}$-monoid structure $\left(\mathcal{A}, \bigoplus, I\right)$ to the $2$-category of pseudoalgebras of an opmonoidal pseudomonad $\left(\mathcal{A}, S\right)$, we should expect some weakness corresponding to the modification components $\omega$, $\gamma$ and $\delta$ for the monoidal $2$-functor $S: \mathcal{A} \rightarrow \mathcal{A}$ underlying the opmonoidal pseudomonad being considered. On the other hand, Corollary \ref{EM of ambistrict is 2-natural} suggests that the weakness of the lifted monoidal structure may be significantly simpler than in an arbitrary monoidal bicategory.
\\
\\
\noindent The first goal of this Section is to collect the results of Subsection \ref{Subsection normal lax monoidality of EM} into a precise statement about $\mathbf{EM}$ being a lax morphism of $\mathbf{Gray}$-monoids. From this it will follow that a pseudomonoid in $\left({\mathbf{Psmnd}\left(\mathbf{Gray}\right)}_\text{ast}, \otimes\right)$ will be sent to some monoidal structure whose associator and left and right unitors are $2$-natural isomorphisms with the following properties.

\begin{itemize}
	\item their component pseudomorphisms are tight, in the sense of Definition \ref{definition tight pseudomorphism}.
	\item they satisfy the usual axioms for a pseudomonoid on the nose rather than up to any invertible modifications.
\end{itemize}

\noindent The goal of Subsection \ref{Subsection gray tensors of 2-naturals} is to describe an appropriate two-dimensional target for the Eilenberg-Moore construction for pseudomonads. This is in lieu of the perhaps more natural three-dimensional target given by the enriched functor category $[\mathbf{2}, \Gr]$; using this as a target would require more complex calculations and the development of the theory of three-dimensional monoidal structures. A step towards describing a simpler two-dimensional target is taken in Definition \ref{definition 2-category G} and Proposition \ref{Proposition EM restricted to ambistrict 2-cells lands in G}. However, the key obstruction is that $2$-natural transformations fail to be closed under composition. This is discussed further in Remark \ref{Remark subtleties in Gray tensor of 2-natural transformations}. The target we will use will finally be described in Definition \ref{Definition Gray 2-nat}.

\subsection{Gray tensor products of $2$-natural transformations}\label{Subsection gray tensors of 2-naturals}

\begin{definition}\label{definition 2-category G}
	Let $\mathbf{G}$ denote the $2$-category whose
	
	\begin{itemize}
		\item objects are $2$-functors $F: \A' \to \A$,
		\item morphisms $F \to G$ are pairs of $2$-functors $P: \A \to \B$ and $Q: \A' \to \B'$ satisfying $PF = GQ$,
		\item $2$-cells $\phi: (P, Q) \Rightarrow (P', Q')$ may only exist if $P=P'$, in which case they consist of a $2$-natural transformation $\phi: Q \Rightarrow Q'$ such that $G.\phi$ is the identity. 
	\end{itemize} 
\end{definition}

\begin{proposition}\label{Proposition EM restricted to ambistrict 2-cells lands in G}
	The $\mathbf{Gray}$-functor $\mathbf{EM}: \mathbf{Psmnd}(\mathbf{Gray}) \to \mathbf{Gray}^\mathbf{2}$ satisfies the following conditions.
	
	\begin{enumerate}
		\item It is tri-fully-faithful.
		\item It has triessential image given by pseudomonadic $2$-functors.
		\item It restricts to a $2$-functor $\mathbf{EM}: \mathbf{Psmnd}(\mathbf{Gray})_\text{ast} \to \mathbf{G}$. 
	\end{enumerate} 
\end{proposition}

\begin{proof}
	Parts (1) and (2) are a consequence of Theorem 3.4 of \cite{Formal Theory of Pseudomonads}; since $\mathbf{Lift}(\mathbf{Gray})$ has objects consisting of pseudomonads, and we may identify these with their pseudomonadic $2$-functors. For part (3), first observe that when $A$ and $B$ are pseudomonadic, then the condition that $B.\phi$ should be the identity is exactly to say that the components of $\phi$ are tight. But by Corollary \ref{EM of ambistrict is 2-natural} $\mathbf{EM}(\phi, \tilde{\phi})$ is $2$-natural and has tight components if and only if $\phi$ is an identity, which is to say that $(\phi, \tilde{\phi})$ is ambistrict.
\end{proof}

\noindent Note that the condition on $2$-cells in Definition \ref{definition 2-category G} that $B.\phi$ should be an identity exactly corresponds to tightness of the components of the $2$-natural transformation $\phi$. In light of Proposition \ref{Proposition EM restricted to ambistrict 2-cells lands in G}, we may try to show that the restriction of $\mathbf{EM}$ is lax monoidal with respect to the $\mathbf{Gray}$-tensor product. However, there are some subtleties which need to be addressed before an appropriate target $\mathbf{Gray}$-monoid containing $\mathbf{G}$ can be given. These subtleties are described in Remark \ref{Remark subtleties in Gray tensor of 2-natural transformations}, to follow.

\begin{remark}\label{Remark subtleties in Gray tensor of 2-natural transformations}
	\noindent Recall that any symmetric monoidal closed, complete and cocomplete category $\mathcal{V}$ is itself a monoidal $\mathcal{V}$-category \cite{Kelly Basic Concepts of Enriched Category Theory}. If particular when $\mathcal{V}= \mathbf{Gray}$, there is a notion of the $\mathbf{Gray}$-tensor product of pseudonatural transformations. If $p: F \Rightarrow G: \mathcal{A} \rightarrow \mathcal{B}$ and $q: H \Rightarrow K: \mathcal{C} \rightarrow \mathcal{D}$ are pseudonatural transformations between $2$-functors then their $\mathbf{Gray}$-tensor product is the pseudonatural transformation $p\otimes q: F\otimes H \Rightarrow G \otimes K$ whose
	
	\begin{itemize}
		\item component on $\left(X, Y\right) \in \mathcal{A} \otimes \mathcal{C}$ is the morphism in $\mathcal{C} \otimes \mathcal{D}$ displayed below.
		
		$$\begin{tikzcd}[font=\fontsize{9}{6}]
			\left(FX{,}HY\right) \arrow[rr, "\left(p_{X}{,}HY\right)"] && \left(GX{,}HY\right) \arrow[rr, "\left(GX{,}q_{Y}\right)"] && \left(GX{,}KY\right) 
			\end{tikzcd}$$
		\item component on $\left(f, Y\right): \left(X, Y\right) \rightarrow \left(X', Y\right)$ is the $2$-cell depicted below.
	\end{itemize}
	
	\begin{equation}\label{Equation component of p otimes q as X varies}
		\begin{tikzcd}[font=\fontsize{9}{6}]
		\left(FX{,}HY\right) \arrow[rr, "\left(p_{X}{,}HY\right)"]\arrow[dd, "\left(Ff{,}HY\right)"']
		&{}\arrow[dd, Rightarrow, shorten = 5, "\left(p_{f}{,}HY\right)"]
		& \left(GX{,}HY\right) \arrow[rr, "\left(GX{,}q_{Y}\right)"] \arrow[dd, "\left(Gf{,}HY\right)"description]
		&{}\arrow[dd, Rightarrow, shorten = 5, "\left(Gf{,}q_{Y}\right)"]
		& \left(GX{,}KY\right)
		\arrow[dd, "\left(Gf{,}KY\right)"]
		\\
		\\
		\left(FX'{,}HY\right) \arrow[rr, "\left(p_{X'}{,}HY\right)"'] &{}& \left(GX'{,}HY\right) \arrow[rr, "\left(GX'{,}q_{Y}\right)"'] &{}& \left(GX'{,}KY\right)
	\end{tikzcd}
\end{equation}
	
	\begin{itemize}
		\item component on $\left(X, g\right): \left(X, Y\right) \rightarrow \left(X, Y'\right)$ is the $2$-cell depicted below.
	\end{itemize} 
	
	\begin{equation}\label{Equation component of p otimes q as Y varies}
	\begin{tikzcd}[font=\fontsize{9}{6}]
		\left(FX{,}HY\right) \arrow[rr, "\left(p_{X}{,}HY\right)"]\arrow[dd, "\left(FX{,}Hg\right)"']
		&{}\arrow[dd, Rightarrow, shorten = 5, "\left(p_{X}{,}Hg\right)"]
		& \left(GX{,}HY\right) \arrow[rr, "\left(GX{,}q_{Y}\right)"] \arrow[dd, "\left(GX{,}Hg\right)"description]
		&{}\arrow[dd, Rightarrow, shorten = 5, "\left(GX{,}q_{g}\right)"]
		& \left(GX{,}KY\right)
		\arrow[dd, "\left(GX{,}Kg\right)"]
		\\
		\\
		\left(FX{,}HY'\right) \arrow[rr, "\left(p_{X}{,}HY'\right)"'] &{}& \left(GX{,}HY'\right) \arrow[rr, "\left(GX{,}q_{Y'}\right)"'] &{}& \left(GX{,}KY'\right)
	\end{tikzcd}
\end{equation}
	
	\noindent In particular, the class of $2$-natural transformations is not closed under $\mathbf{Gray}$-tensor product. Any attempt at restricting this to a $\mathbf{Gray}$-monoid structure on $2$-$\mathbf{Cat}$ would require $2$-natural transformations of the form $\A \otimes q: \A \otimes H \Rightarrow \A \otimes K$. But taking $p$ to be the identity in Equation \ref{Equation component of p otimes q as Y varies}, it is clear that this still fails to be strictly natural in generators in $\mathcal{A}\otimes \mathcal{C}$ from $\mathcal{A}$, unless $q$ is in fact an identity. Rather, its component is given by an interchanger in the $\mathbf{Gray}$-tensor product $\A \otimes \mathcal{D}$. Note that similar comments apply for tensors of the form $p\otimes \mathcal{C} : F \otimes \mathcal{C} \Rightarrow G \otimes \mathcal{C}$, as can be seen in Equation \ref{Equation component of p otimes q as X varies}. Even when $p:F \Rightarrow G$ is $2$-natural, this will be a pseudonatural transformation whose components are interchangers of its $1$-cell components in the $\mathbf{Gray}$-tensor product $\mathcal{C}\otimes \mathcal{D}$. 
\end{remark}

\noindent Definition \ref{Definition Gray 2-nat}, to follow, addresses the issue just raised in Remark \ref{Remark subtleties in Gray tensor of 2-natural transformations}, to form a monoidal $2$-category generated by $2$-natural transformations. This will be used to collect the results of Subsection \ref{Subsection normal lax monoidality of EM} into a precise statement about monoidality of $\mathbf{EM}$, and to clarify what further calculations are needed to establish a lifted monoidal structure on pseudoalgebras.

\begin{definition}\label{Definition Gray 2-nat}
	Let $\mathbf{Gray}_{2\text{-nat}}$ be the strict monoidal $2$-category whose 
	\begin{itemize}
		\item underlying monoidal category is $\left(2\text{-}\mathbf{Cat}, \otimes, \mathbf{1}\right)$,
		\item generating $2$-cells are $2$-natural transformations,
		\item relations between $2$-cells are equalities of pseudonatural transformations, and the middle-four interchange relation.
	\end{itemize} 

\noindent Hence the $2$-cells in $\mathbf{Gray}_{2\text{-nat}}$ are congruence classes of pseudonatural transformations. Each representative in such a congruence class admits some decomposition into (composites of whiskerings of) $\mathbf{Gray}$-tensor products of $2$-natural transformations. The congruence is generated by the requirement of middle four interchange in the underlying $2$-category of $\mathbf{Gray}_{2\text{-nat}}$.
\end{definition} 

\begin{remark}\label{Remark lit review on gray tensors of 2-nats}
	Notions suggestive of enrichment over a would-be monoidal $2$-category structure \cite{Garner Shulman Enriched Categories as a Free Cocompletion} involving $\Gr$ tensor products and $2$-natural transformations have found applications in the literature, despite the fact that $2$-natural transformations are not closed under $\Gr$ tensor products. As briefly discussed in Remark 3.3.3 of \cite{Miranda strictifying operational coherences}, the $3$-pseudofunctors and $3$-pseudonatural transformations introduced there are reminiscent of enriched pseudofunctors and enriched pseudonatural transformations. The weak maps between $\Gr$-categories described in \cite{Bourke Lobbia Skew Approach} are special cases of $3$-pseudofunctors in which the unit constraints are identities.
	\\
	\\
	\noindent All of these notions can indeed be seen as special cases of enriched notions over $\mathcal{V} = \mathbf{Gray}_\text{2-nat}$. In particular, their $2$-cell constraints are generating $2$-cells (i.e. $2$-natural transformations). Moreover, the examples genuinely satisfy the necessary equations on $2$-cells, rather than just satisfying these conditions up to the congruence generated by middle four interchange. In Definition \ref{Definition Gray pseudomonoid} we introduce a similarly special pseudomonoid over $\Gr_\text{2-nat}$, as the lifted monoidal structure on $\A^S$ will be of this form. An explicit description of the $2$-cells in $\mathbf{Gray}_{2\text{-nat}}$ will not be needed to lift symmetric monoidal structures to Eilenberg-Moore objects of pseudomonads, since most equations of relevance will turn out to hold even before quotienting by the congruence generated by middle-four interchange.
\end{remark}

\subsection{Monoidal structure on pseudoalgebras}\label{subsection monoidal structure on pseudoalgebras}

\begin{theorem}\label{Theorem EM is normal lax Gray functor}
	The $2$-functor $\mathbf{EM}: {\mathbf{Psmnd}\left(\mathbf{Gray}\right)}_\text{ast} \rightarrow \mathbf{Gray}_\text{2-nat}$ admits the structure of a normal lax monoidal $2$-functor when equipped with a compositor given by the pseudonatural transformation $M_{-, ?}$ of Corollary \ref{Corollary compositor for EM is a well-defined pseudonatural transformation}. All further data of this normal lax monoidal $2$-functor are given by identities.
\end{theorem}

\begin{proof}
	Following Proposition \ref{associativity and unitality of EM}, it suffices to show that there are equalities of $2$-natural isomorphisms as the pseudomonads are varied along morphisms of pseudomonads. We give details for associativity, as the arguments for left and right unit law are similar and indeed simpler.
	\\
	\\
	\noindent The condition for associativity is trivial when varying $(\A, S)$, since the components $M_{g, 1}$ are identities by Proposition \ref{Proposition M is natural in the first component}. To see that associativity holds as $(\B, T)$ is varied, observe that both $M_{S\otimes T', R}.(M_{1, h}\otimes \mathcal{C}^{R})$ and $M_{1, h\otimes R}.(\mathcal{A}^S\otimes M_{T, R})$ have component on $(\mathbf{X}, \mathbf{Y}, \mathbf{Z})$ given by the tight pseudomorphism whose corresponding $2$-cell in $\A\otimes \B\otimes \mathcal{C}$ is as depicted below.
	
	$$\begin{tikzcd}[column sep = 20, font=\fontsize{9}{6}]
		\left(SX{,}T'HY{,}RZ\right) \arrow[rr, "\left(1{,}h_{Y}{,}1\right)"]
		\arrow[dd, "\left(x{,}1{,}1\right)"']
		&{}\arrow[dd, Rightarrow, shorten = 15, "\left(x{,}h_{Y}{,}1\right)"]
		& \left(SX{,} HTY{,}RZ\right) \arrow[dd, "\left(x{,}1{,}1\right)"]
		\\
		\\
		\left(X{,}T'HY{,}RZ\right) \arrow[rr, "\left(1{,}h_{Y}{,}1\right)"']
		&{}& \left(X{,}HTY{,}RZ\right) \arrow[rr, "\left(1{,}Hy{,}1\right)"'] && \left(X{,}HY{,}RZ\right)\arrow[rr, "\left(1{,}1{,}z\right)"'] &&\left(X{,}HY{,}Z\right)
	\end{tikzcd}$$
	
	\noindent Let $(K ,k ,k_{0}, k_{2}): (\mathcal{C}, R) \to (\mathcal{C}', R')$ be a morphism of pseudomonads. To see that associativity holds on $(K ,k ,k_{0}, k_{2})$, observe that both of the $2$-natural transformations $M_{S\otimes T, k}.(M_{S, T}\otimes \mathcal{C}^{R})$ and $M_{S, T\otimes k}.(\A^S \otimes M_{T, k})$ have component on $(\mathbf{X}, \mathbf{Y}, \mathbf{Z}) \in \A^S\otimes \B^T \otimes \mathcal{C}^{R}$ given by the tight pseudomorphism whose corresponding $2$-cell is as displayed below.
	
	$$\begin{tikzcd}[column sep = 18, font=\fontsize{9}{6}]
		\left(SX{,}HTY{,}R'KZ\right) \arrow[rr, "\left(1{,}1{,}k_{Z}\right)"]
		\arrow[dd, "\left(x{,}1{,}1\right)"']
		&{}\arrow[dd, Rightarrow, shorten = 15, "\left(x{,}1{,}k_Z\right)"]
		& \left(SX{,} HTY{,}RZ\right) \arrow[dd, "\left(x{,}1{,}1\right)"]
		\\
		\\
		\left(X{,}HTY{,}RZ\right)
		\arrow[dd, "(1{,}HTy{,}1)"'] \arrow[rr, "\left(X{,}1{,}k_{Z}\right)"']
		&{}\arrow[dd, Rightarrow, shorten = 15, "(X{,}HTy{,}k_{Z})"]
		& \left(X{,}HTY{,}RZ\right) \arrow[dd, "\left(1{,}HTy{,}1\right)"]
		\\
		\\
		(X{,}HY{,}RZ)\arrow[rr, "(1{,}1{,}k_{Z})"']
		&{}& \left(X{,}HY{,}RZ\right)\arrow[rr, "\left(1{,}1{,}z\right)"'] &&\left(X{,}HY{,}Z\right)
	\end{tikzcd}$$
\end{proof}

\noindent Proposition \ref{Proposition pseudomonoid structure on pseudoalgebras in Gray 2nat}, to follow, gives a pseudomonoid structure on the $2$-category of pseudoalgebras. It is an easy consequence of Theorem \ref{Theorem EM is normal lax Gray functor}. Following this, we will carefully analyse the associator and left and right unitor of the pseudomonoid, and see that in fact the lifted structure is a particularly simple kind of monoidal bicategory.

\begin{proposition}\label{Proposition pseudomonoid structure on pseudoalgebras in Gray 2nat}
	Let $(\A, \bigoplus, I, S, \eta, \mu)$ be a semi-strict opmonoidal pseudomonad and consider the $2$-category of pseudoalgebras $\A^S$. There is a pseudomonoid structure on $\A^S$ in the $\mathbf{Gray}$-monoid $\mathbf{Gray}_\text{2-nat}$.
\end{proposition}

\begin{proof}
	Recall (see for example the comments following Definition 1 in \cite{Low dimensional structures formed by tricategories}) that a pseudomonoid in a $\mathbf{Gray}$-monoid $\mathcal{V}$ is a lax trihomomorphism $\mathbf{1} \to \Sigma\mathcal{V}$, where $\Sigma$ denotes the suspension. As a lax monoidal $2$-functor, $\mathbf{EM}: {\mathbf{Psmnd}(\mathbf{Gray})}_\text{ast} \to \mathbf{Gray}_{2-nat}$ is in particular a lax trihomomorphism between the suspensions. Recall also that lax trihomomorphisms are closed under composition. The desired pseudomonoid in $\mathbf{Gray}_\text{2-nat}$ is given as the composite displayed below.
	
	$$\begin{tikzcd}
		\mathbf{1} \arrow[rr, "\mathbb{S}"]&& {\mathbf{Psmnd}(\mathbf{Gray})}_\text{ast}\arrow[rr, "\mathbf{EM}"]&& \mathbf{Gray}_\text{2-nat}
	\end{tikzcd}$$
	
\end{proof}

\noindent A careful analysis of the pseudomonoid structure on $\mathcal{A}^{S}$ is needed in order to show that it is in fact a monoidal bicategory, rather than one whose associators and left and right unitors are given by equivalence classes of certain pseudonatural transformations. The structure will turn out to be a fairly simple kind of monoidal bicategory, to be described in Definition \ref{Definition Gray pseudomonoid}, to follow.

\begin{definition}\label{Definition Gray pseudomonoid}
	A \emph{$\mathbf{Gray}$-pseudomonoid} is a cubical tricategory with one object (Definition 8.1 of \cite{Gurski Coherence in Three Dimensional Category Theory}), in which the unitors and associators $\alpha$, $\lambda$ and $\rho$ are $2$-natural isomorphisms, and the invertible modification components are all identities.
\end{definition}

\begin{remark}\label{Remark ways of thinking about Gray pseudomonoids}
	We outline several ways of thinking about $\Gr$ pseudomonoids. They are like $\Gr$ monoids except that they only satisfy the usual monoid axioms up to $2$-natural isomorphisms which satisfy the usual monoidal category axioms on the nose. Alternatively, $\Gr$-pseudomonoids are like $\mathcal{V} = \Cat$ enriched monoidal categories, except that their tensor product is defined out of the $\Gr$-tensor product rather than the cartesian product of $2$-categories. Finally, they are pseudomonoids in $\Gr_\text{2-nat}$ whose associators and left and right unitors are the generating $2$-cells described in Definition \ref{Definition Gray 2-nat}, and for which the usual pseudomonoid axioms hold as equations between generating $2$-cells rather than up to the congruence generated by middle four interchange.
\end{remark} 

\begin{theorem}\label{monoidal bicategory structure on pseudoalgebras}
	There is a $\mathbf{Gray}$-pseudomonoid, with data as follows.
	
		\begin{itemize}
		\item The unit is given by $\overline{I}: \mathbf{1} \rightarrow \mathcal{A}^{S}$, which corresponds to the pseudoalgebra described in Example \ref{Example pseudoalgebra structure on the unit}.
		\item The tensor product of pseudoalgebras is given by the $2$-functor displayed below.
		
		$$\begin{tikzcd}[font=\fontsize{9}{6}]
			\mathcal{A}^{S}\otimes \mathcal{A}^{S} \arrow[rr, "M_{S{,}S}"] && {\left(\mathcal{A}\otimes \mathcal{A}\right)}^{\left(S \otimes S\right)} \arrow[rr, "\overline{\left(\bigoplus{,}\chi\right)}"] && \mathcal{A}^{S}
		\end{tikzcd}$$

		\item The associator $\alpha$ is given by the following pasting in $2$-$\mathbf{Cat}$, where the $2$-natural isomorphism in the top right square is as described in Proposition \ref{Proposition 2-cell component of compositor as the second pseudomonad varies}, while the $2$-natural isomorphism in the bottom right square corresponds to $\mathbf{EM}\left(1, \omega\right)$ as described in Remark \ref{Remark explicit EM}.
		
		\begin{equation}\label{Equation associator for lifted monoidal structure}
			\begin{tikzcd}[font=\fontsize{9}{6}]
				\mathcal{A}^{S}\otimes \mathcal{A}^{S}\otimes \mathcal{A}^{S} \arrow[rr, "\mathcal{A}^{S}\otimes M_{S{,}S}"]\arrow[dd, "M_{S{,}S} \otimes \mathcal{A}^{S}"']
				&& \mathcal{A}^{S}\otimes {\left(\mathcal{A}\otimes \mathcal{A}\right)}^{\left(S \otimes S\right)} \arrow[rr, "\mathcal{A}^{S}\otimes \overline{\bigoplus}"] \arrow[dd, "M_{S{,}S\otimes S}"description]
				&{}\arrow[dd, Rightarrow, shorten = 15, "M_{1{,}\chi}"]& \mathcal{A}^{S}\otimes \mathcal{A}^{S}\arrow[dd, "M_{S{,}S}"]
				\\
				&=
				\\
				{\left(\mathcal{A}\otimes \mathcal{A}\right)}^{\left(S \otimes S\right)} \otimes \mathcal{A}^{S} \arrow[dd, "\overline{\bigoplus}\otimes \mathcal{A}^{S}"']\arrow[rr, "M_{S\otimes S{,}S}"]
				&&{\left(\mathcal{A}\otimes \mathcal{A}\otimes \mathcal{A}\right)}^{\left(S \otimes S\otimes S\right)}\arrow[dd, "\overline{\bigoplus\otimes \mathcal{A}}"description]\arrow[rr, "\overline{ \mathcal{A}\otimes\bigoplus}"]
				&{}\arrow[dd, Rightarrow, shorten = 15, "\overline{\omega}"]
				&{\left(\mathcal{A}\otimes \mathcal{A}\right)}^{\left(S \otimes S\right)}\arrow[dd, "\overline{\bigoplus}"]
				\\
				&=
				\\
				\mathcal{A}^{S} \otimes \mathcal{A}^{S}\arrow[rr, "M_{S{,}S}"'] && {\left(\mathcal{A}\otimes \mathcal{A}\right)}^{\left(S\otimes S\right)}\arrow[rr, "\overline{\bigoplus}"'] &{}& \mathcal{A}^{S}
			\end{tikzcd}
		\end{equation}
	
		\item left unitor $\lambda$ is given by the pasting in $2$-$\mathbf{Cat}$ depicted below left,
		
		\item right unitor $\rho$ is given by the pasting in $2$-$\mathbf{Cat}$ depicted below right,
		
		\begin{equation}\label{Equation left and right unitors for lifted monoidal structure}
			\begin{tikzcd}[font=\fontsize{9}{6}]
				\mathcal{A}^{S} \arrow[rr, "\overline{I} \otimes 1_{\mathcal{A}^{S}}"]
				\arrow[rrdd, bend right = 20, "\overline{I\otimes 1}"description]
				\arrow[rrdddd, bend right = 40, "1_{\mathcal{A}^{S}}"']
				&{}
				&\mathcal{A}^{S}\otimes \mathcal{A}^{S} \arrow[dd, "M_{S{,}S}"]
				\\
				&=
				\\
				&&{\left(\mathcal{A}\otimes \mathcal{A}\right)}^{\left(S \otimes S\right)} \arrow[dd, "\overline{\bigoplus}"]\arrow[dd, shift right = 10, Rightarrow, shorten =10, "\overline{\gamma}"]
				\\
				\\
				&&\mathcal{A}^{S}&{}
			\end{tikzcd}\begin{tikzcd}[font=\fontsize{9}{6}]
				\mathcal{A}^{S} \arrow[rr, "1_{\mathcal{A}^{S}}\otimes \overline{I}"]
				\arrow[rrdd, bend right = 20, "\overline{1\otimes I}"description]
				\arrow[rrdddd, bend right = 40, "1_{\mathcal{A}^{S}}"']
				&{}&\mathcal{A}^{S}\otimes \mathcal{A}^{S} \arrow[dd, "M_{S{,}S}"]
				\\
				&{}\arrow[u, Rightarrow, shorten = 5, "M_{1{,}\iota}"']
				\\
				&&{\left(\mathcal{A}\otimes \mathcal{A}\right)}^{\left(S \otimes S\right)} \arrow[dd, "\overline{\bigoplus}"]\arrow[dd, shift right = 10, Leftarrow, shorten =10, "\overline{\delta}"]
				\\
				\\
				&&\mathcal{A}^{S}
			\end{tikzcd}
		\end{equation}
		
	\end{itemize}
\end{theorem}

\begin{proof}
	By Proposition \ref{Proposition pseudomonoid structure on pseudoalgebras in Gray 2nat}, there is a pseudomonoid structure on $\A^S$ in the $\mathbf{GRAY}$-monoid $\Gr_\text{2-nat}$. We give more explicit descriptions of the structure, to verify that this is in fact a $\Gr$-pseudomonoid. The tensor product of a pair of pseudoalgebras $(\mathbf{X}, \mathbf{Y})$ is the pseudoalgebra with underlying object $X \bigoplus Y$, structure map displayed below
	
	$$\begin{tikzcd}[font=\fontsize{9}{6}]
		S(X\bigoplus Y) \arrow[rr, "\chi_{X{,}Y}"]&& SX\bigoplus SY \arrow[rr, "x\bigoplus SY"] && X \bigoplus SY \arrow[rr, "X\bigoplus y"] && X \bigoplus Y
	\end{tikzcd}$$
	
	\noindent and unit and multiplication constraints given as displayed below.
	
	$$\begin{tikzcd}[font=\fontsize{9}{6}]
		S^{2}(X\bigoplus Y)\arrow[dddddd, "\mu_{X\bigoplus Y}"'] \arrow[rr, "S(\chi_{X{,}Y})"]
		&{}\arrow[dddddd, Rightarrow, shorten = 60, "\mu_{2{,}X{,}Y}"']
		& S(SX\bigoplus SY)\arrow[dd, "\chi_{SX{,}SY}"'] \arrow[rr, "S(x\bigoplus SY)"]
		&{}\arrow[dd, Rightarrow, shorten = 15, "\chi_{x{,}SY}"]
		& S(X \bigoplus SY) \arrow[rr, "S(X\bigoplus y)"]\arrow[dd, "\chi_{X{,}SY}"]
		&{}\arrow[dd, Rightarrow, shorten = 15, "\chi_{X{,}y}"]
		& S(X \bigoplus Y)\arrow[dd, "\chi_{X{,}Y}"]
		\\
		\\
		&&S^2X\bigoplus S^2Y\arrow[rr, "Sx\bigoplus S^2Y"]
		\arrow[dd, "\mu_X\bigoplus S^2{Y}"']
		&{}\arrow[dd, Rightarrow, shorten = 15, "x_{2}\bigoplus S^2Y"]
		&SX\bigoplus S^2Y\arrow[rr, "SX\bigoplus Sy"]
		\arrow[dd, "x\bigoplus S^2Y"]
		&{}\arrow[dd, Rightarrow, shorten = 15, "x_{Sy}"]
		&SX\bigoplus SY\arrow[dd, "x\bigoplus SY"]
		\\
		\\
		&&SX\bigoplus S^2Y\arrow[rr, "x\bigoplus S^2Y"]
		\arrow[dd, "SX\bigoplus \mu_{Y}"']
		&{}\arrow[dd, Rightarrow, shorten = 15, "x_{\mu_{Y}}"]&X\bigoplus S^{2}Y
		\arrow[rr, "X\bigoplus Sy"]
		\arrow[dd, "X\bigoplus \mu_{Y}"]
		&{}\arrow[dd, Rightarrow, shorten = 15, "X\bigoplus y_{2}"]
		&X\bigoplus SY
		\arrow[dd, "X\bigoplus y"]
		\\
		\\
		S(X\bigoplus Y)\arrow[rr, "\chi_{X{,}Y}"']
		&{}& SX\bigoplus SY \arrow[rr, "x\bigoplus SY"']
		&{}& X \bigoplus SY \arrow[rr, "X\bigoplus y"'] 
		&{}& X \bigoplus Y
	\end{tikzcd}$$

	$$\begin{tikzcd}[font=\fontsize{9}{6}]
		X\bigoplus Y\arrow[rrdd, "\eta_{X}\bigoplus Y"]
		\arrow[rrrr, "\eta_{X\bigoplus Y}"]
		\arrow[rrdddd, bend right = 30, "1_{X\bigoplus Y}"']
		&&{}\arrow[dd, Leftarrow, shorten = 15, "\eta_{2{,}X{,}Y}"]
		&&S(X\bigoplus Y) \arrow[dd, "\chi_{X{,}Y}"]
		\\
		&{}\arrow[dd, Leftarrow, shorten = 15, "x_{0}\bigoplus Y"]
		\\
		&&SX\bigoplus Y\arrow[dd, "x\bigoplus Y"]
		\arrow[rr, "SX\bigoplus\eta_{Y}"]
		&{}\arrow[dd, Leftarrow, shorten = 15, "x\bigoplus \eta_{Y}"]
		&SX\bigoplus SY \arrow[dd, "x\bigoplus SY"] 
		\\
		&{}
		\\
		&&X\bigoplus Y\arrow[rr, "X\bigoplus \eta_{Y}"]\arrow[rrdd, bend right = 30, "1_{X\bigoplus Y}"']
		&{}\arrow[dd, Leftarrow, shorten = 20, "X\bigoplus y_{0}"]
		&X \bigoplus SY \arrow[dd, "X\bigoplus y"] 
		\\
		\\
		&&&{}&X \bigoplus Y
	\end{tikzcd}$$

\noindent Regarding the associator, observe that the top left square in Equation \ref{Equation associator for lifted monoidal structure} commutes by Proposition \ref{associativity and unitality of EM} part (2), and the bottom left square in Equation \ref{Equation associator for lifted monoidal structure} commutes by Proposition \ref{Proposition M is natural in the first component}. The component of the associator on a triple $(\mathbf{X}, \mathbf{Y}, \mathbf{Z})$ is the tight pseudomorphism whose $2$-cell component given by the pasting in $\mathcal{A}$ depicted below.

$$\begin{tikzcd}[font=\fontsize{9}{6}, column sep = 18]
	S(X \bigoplus Y \bigoplus Z)
	\arrow[rr, "\chi_{X{,}Y\bigoplus Z}"]
	\arrow[dd, "\chi_{X\bigoplus Y{,}Z}"']
	&{}\arrow[dd, Leftarrow, shorten = 15, "\omega_{X{,}Y{,}Z}"]
	&SX\bigoplus S(Y\bigoplus Z)
	\arrow[rr, "x\bigoplus 1"]
	\arrow[dd, "1\bigoplus \chi_{Y{,}Z}"]
	&{}\arrow[dd, Leftarrow, shorten = 15, "x_{\chi_{Y{,}Z}}"]
	&X\bigoplus S(Y\bigoplus Z)
	\arrow[dd, "1\bigoplus \chi_{Y{,}Z}"]
	\\
	\\
	S(X\bigoplus Y)\bigoplus SZ 
	\arrow[rr, "\chi_{X{,}Y}\bigoplus 1"']
	&{}&SX\bigoplus SY\bigoplus SZ 
	\arrow[rr, "x\bigoplus 1\bigoplus 1"'] 
	&{}& X \bigoplus SY \bigoplus SZ
	\arrow[rr, "1\bigoplus y\bigoplus 1"'] 
	&& X \bigoplus Y\bigoplus SZ
	\arrow[rr, "1 \bigoplus z"']
	&&X\bigoplus Y\bigoplus Z
\end{tikzcd}$$

\noindent Regarding the description of left unitor, observe that the upper triangle depicted to the left of Equation \ref{Equation left and right unitors for lifted monoidal structure} indeed commutes by Proposition \ref{Proposition M is natural in the first component}. The components of the left and right unitor are also both tight pseudomorphisms. The tight pseudomorphism giving the left unitor on $\mathbf{X}$ has $2$-cell component as displayed below.

$$\begin{tikzcd}[font=\fontsize{9}{6}]
	S(I\bigoplus X) \arrow[rr, "\chi_{I{,}X}"]\arrow[rrrr, bend right = 30, "1_{SX}"']
	&& SI\bigoplus SX\arrow[d, Rightarrow, shorten = 7, "\gamma_{X}"]
	\arrow[rr, "\iota\bigoplus 1"] && I \bigoplus SX \arrow[rr, "1\bigoplus x"] && I \bigoplus X
	\\
	&&{}
\end{tikzcd}$$

\noindent Meanwhile, the tight pseudomorphism giving the right unitor on $\mathbf{X}$ has $2$-cell component as displayed below.

$$\begin{tikzcd}[font=\fontsize{9}{6}]
	S(X\bigoplus I) \arrow[rr, "\chi_{X{,}I}"]\arrow[rrrdd, bend right = 20, "1_{SX}"']
	&& SX\bigoplus SI\arrow[ldd, Rightarrow, shorten = 20, "\delta_{X}"]
	\arrow[rdd, "1\bigoplus \iota"]
	\arrow[rr, "x\bigoplus 1"] 
	&& X \bigoplus SI
	\arrow[ldd, Rightarrow, shorten = 20, "x_{\iota}"]
	\arrow[rr, "1\bigoplus \iota"] && X \bigoplus I
	\\
	\\
	&{}&&SX\bigoplus I
	\arrow[rrruu,bend right = 20,  "x\bigoplus 1"']
\end{tikzcd}$$

\noindent Note that $2$-naturality of $\overline{\omega}$, $\overline{\gamma}$ and $\overline{\delta}$ follows from Corollary \ref{EM of ambistrict is 2-natural}, since the components $\omega$, $\gamma$ and $\delta$ of the opmonoidal $2$-functor $S: \mathcal{A} \to \A$ determine ambistrict pseudomonad transformations, as discussed in Subsection \ref{Section monoidal pseudomonads as pseudomonoids in pseudomonads}. Observe that the associator and left and right unitors for the pseudomonoid structure on $\A^S$ described above are indeed well-defined pasting diagrams in $2$-$\mathbf{Cat}$, not just $\mathbf{Gray}_\text{2-nat}$.
\\
\\
\noindent As such, all $2$-cells in sight are genuine $2$-natural transformations, rather than equivalence classes of certain pseudonatural transformations formed as $\mathbf{Gray}$-tensor products of $2$-natural transformations. It remains to verify that the pseudomonoid axioms hold as equations between $2$-natural transformations, rather than just as equations between $2$-cells in $\Gr_\text{2-nat}$. The proof is via direct calculations using the relations in the presentation of the $\Gr$ tensor product, and the axioms for the opmonoidal $2$-functor $\left(S, \iota, \chi, \gamma, \omega, \delta\right): (\A, \bigoplus, I) \to (\A, \bigoplus, I)$. Details are deferred to Appendix \ref{Appendix monoidal bicategory structure on pseudoalgebras}.
\end{proof}

\noindent We emphasise that the only obstruction to restricting the $\mathbf{Gray}$-pseudomonoid structure of Theorem \ref{monoidal bicategory structure on pseudoalgebras} to the underlying category of $\A^S$ is that the tensor product emanates from $\A^S \otimes \A^S$, rather than the cartesian product $\A^S \times \A^S$. As such, there is only a cubical pseudofunctor $\A^S \times \A^S \to \A^S$, rather than a genuine $2$-functor, and so this does not restrict to a genuine functor between underlying categories. Corollary \ref{Corollary tensor product of pseudoalgebras special properties}, to follow, records some further coherence properties of the lifted monoidal structure on pseudoalgebras.

\begin{corollary}\label{Corollary tensor product of pseudoalgebras special properties}
	\noindent The $\mathbf{Gray}$-pseudomonoid of Theorem \ref{monoidal bicategory structure on pseudoalgebras} has the following properties. 
	\begin{enumerate}
		\item the forgetful functor $\A^S \to \A$ is strict monoidal.
		\item the left unitor is strict if $\gamma$ is an identity.
		\item the right unitor is strict if $\delta$ and $\iota$ are identities.
		\item the associator is strict if $\omega$ and the $1$-cell components of $\chi$ are identities.
		\item the unit is a strict algebra if the monoidal pseudonatural transformations $\left(\eta, \eta_{0}, \eta_{2}\right)$ and $\left(\mu, \mu_{0}, \mu_{2}\right)$ are unital.
		\item even if the opmonoidal pseudomonad is a strictly monoidal $2$-monad, the tensor product of strict algebras typically does not have identity unit or multiplication constraints.
	\end{enumerate}
\end{corollary}

\begin{proof}
	Part (1) amounts to the fact that the lifted tensor product of pseudoalgebras, and between pseudomorphisms and pseudoalgebras, is given at the level of underlying objects and morphisms in terms of their corresponding tensors in $(\A, \bigoplus, I)$. This is easy to inspect from the explicit descriptions given in the proof of Theorem \ref{monoidal bicategory structure on pseudoalgebras}. Parts (2), (3) and (4) are easy observations from the explicit descriptions given in the proof of Theorem \ref{monoidal bicategory structure on pseudoalgebras}. Part (5) is also an easy observation from the explicit description of the unit, given in Example \ref{Example pseudoalgebra structure on the unit}. For part (6), observe that the tensor product of strict algebras $(X, x)$ and $(Y, y)$ has unit constraint given in terms of the interchanger $x_{\eta_{Y}}$, and has multiplication constraint given in terms of interchangers such as $x_{Sy}$. These will typically not be identities.
\end{proof}

\section{Symmetries for opmonoidal pseudomonads}\label{Section symmetries for opmonoidal pseudomonads}

\noindent We now extend our results about lifting monoidal structures to $2$-categories of pseudoalgebras to analogous results about lifting braidings, syllapses and symmetries. In this setting, braidings typically give rise to pseudomonad transformations that need not be ambistrict. As such, the simplifications in Section \ref{Section lifting monoidal structures} via which we were able to work with a $\mathbf{Gray}$-monoid of pseudomonads are too `low dimensional' to also apply to lifting braidings, syllapses and symmetries.
\\
\\
\noindent One may hope for a relatively simple three-dimensional symmetric monoidal structure on $\mathbf{Psmnd}(\mathbf{Gray})$. However, this would need to be more complicated than a $\mathcal{V} = \mathbf{Gray}$ enriched monoidal structure. This is because $\mathcal{V}$-enriched monoidal categories have an underlying monoidal category but as we have already observed, the interchange law in any monoidal structure on $\PsG$ needs to be weak as per Lemma \ref{interchanger between morphisms of pseudomonads in Gray monoid structure on Psmnd Gray}.
\\
\\
\noindent We leave the pursuit of a three-dimensional monoidal structure on $\PsG$ to future research. For our goals of lifting symmetric monoidal structures to Eilenberg-Moore objects of pseudomonads, it is enough to do the following.

\begin{itemize}
	\item Recognise compatibility structures between pseudomonads and braidings, syllapses or symmetries as pseudomonad transformations and pseudomonad modifications. This is done in Subsection \ref{Braided and symmetric opmonoidal monads}, which also reviews compatibility between these structures and, transfors between them such as $2$-functors and pseudonatural transformations.
	\item Establish a relatively simple fragment of the interaction between the Eilenberg-Moore construction for pseudomonads and the $\mathbf{Gray}$-tensor product. This amounts to an invertible modification whose components are tight $2$-natural transformations, witnessing symmetric monoidality of $\mathbf{EM}: \PsG_\text{ast} \to \mathbf{Gray}_\text{2-nat}$. This structure is established in Subsection \ref{Subsection symmetric opmonoidality of EM}.
	\item Use the components of $\Xi$ to build the braid and symmetry structure on $\A^S$. This is done in Subsection \ref{Subsection on lifting braidings, syllapses and symmetries}. The axioms for the symmetric monoidal bicategory structure on $\A^S$ the follow immediately from their analogues for the structure on $\A$, by faithfulness on $2$-cells of $U^S: \A^S \to \A$.
\end{itemize}
\subsection{Braided, Sylleptic and Symmetric Opmonoidal Pseudomonads}\label{Braided and symmetric opmonoidal monads}

\noindent Recall the notion of a braided monoidal bicategory from 2.4 of \cite{Gurski Loop Spaces}, and the notion of a braided $\Gr$-monoid from Definition 2.1 of \cite{Crans Generalised Centers}, where it is called a semi-strict monoidal $2$-category. We recall appropriate structure or properties for opmonoidal $2$-functors or pseudonatural transformations to be compatible with braid, syllepsis or symmetry structures. As we describe in Remark \ref{Remark braiding and syllapsis comparison to literature}, these are routine adaptations of existing definitions in the literature. The tricategory structures that these form have been established in \cite{Schommer Piers Classification of 2D ETFTs}, and guide our definition of braided, sylleptic and symmetric opmonoidal pseudomonads in Definition \ref{Definition braided opmonoidal pseudomonad}.

\begin{definition}\label{Braided opmonoidal 2-functors}
	Let $\left(\mathcal{A}, \bigoplus, I\right)$ and $\left(\mathcal{B}, \bigoplus, I\right)$ be braided $\mathbf{Gray}$-monoids. A \emph{braided opmonoidal $2$-functor} consists of an opmonoidal $2$-functor $\left(F, \chi, \iota, \gamma, \delta, \omega\right): \left(\mathcal{A}, \bigoplus, I\right) \rightarrow \left(\mathcal{B}, \bigoplus, I\right)$ between underlying $\mathbf{Gray}$-monoids, as well as an invertible modification $\Omega$ as depicted below, subject to the axioms listed below.
		
		$$\begin{tikzcd}[column sep = 15, font=\fontsize{9}{6}]
			\mathcal{A}\otimes \mathcal{A} \arrow[ddd, "F \otimes F"']
			\arrow[rrrr, "\bigoplus"]\arrow[rrd, "\tau"] &&{}\arrow[d, Rightarrow, shorten = 5, "\beta"]
			&&\mathcal{A}\arrow[ddd, "F"]
			\\
			&&\mathcal{A}\otimes\mathcal{A}\arrow[rru, "\bigoplus"']
			\arrow[ddd, "F\otimes F"description]
			&{}\arrow[dd, Rightarrow, shorten = 10, "\chi"]
			\\
			&=&&&&&\cong_{\Omega}&{}
			\\
			\mathcal{B}\otimes \mathcal{B}
			\arrow[rrd, "\tau"'] &&{}&{}&\mathcal{B}
			\\
			&&\mathcal{B}\otimes \mathcal{B}\arrow[rru, "\bigoplus"']
		\end{tikzcd}\begin{tikzcd}[row sep = 25, column sep = 15,  font=\fontsize{9}{6}]
			\mathcal{A}\otimes \mathcal{A} \arrow[ddd, "F \otimes F"']
			\arrow[rrrr, "\bigoplus"] &&{}\arrow[ddd, Rightarrow, shorten = 25, "\chi"]
			&&\mathcal{A}\arrow[ddd, "F"]
			\\
			&&&
			\\
			\\
			\mathcal{B}\otimes \mathcal{B}
			\arrow[rrd, "\tau"']\arrow[rrrr, "\bigoplus"]
			&&{}\arrow[d, Rightarrow, shorten = 5, "\beta"]
			&{}&\mathcal{B}
			\\
			&&\mathcal{B}\otimes \mathcal{B}\arrow[rru, "\bigoplus"']
		\end{tikzcd}$$
		
		\noindent AXIOMS 
		\\
		\\
\noindent \emph{(right braiding law)}
		
	$$\begin{tikzcd}[font = \fontsize{9}{6}, column sep = 10]
			&&F(X\bigoplus Y\bigoplus Z)\arrow[dddd, Rightarrow, shorten = 20, "\chi_{\beta_{X{,}Y}{,}1}"]
			\arrow[lldd, "\chi_{X\bigoplus Y{,}Z}"']
			\arrow[rrdd, "F(\beta_{X{,}Y}\bigoplus 1)"description]
			\arrow[rrrr, "F\beta_{X{,}Y\bigoplus Z}"]
			&&{}\arrow[dd, Rightarrow, shift right = 2,shorten = 15, "Fr_{Y{,}Z}^{X}"]
			&&F(Y\bigoplus Z\bigoplus X)
			\arrow[dddd, Rightarrow, shorten =20, "\chi_{1{,}\beta_{X{,}Z}}"]
			\arrow[ddrr, "\chi_{Y{,}Z\bigoplus X}"]
			\\
			\\
			F(X\bigoplus Y)\bigoplus FZ
			\arrow[dddd, "\chi_{X{,}Y}\bigoplus 1"']
			\arrow[rrdd, "F\beta_{X{,}Y}\bigoplus 1"description]
			&&&&F(Y\bigoplus X\bigoplus Z)
			\arrow[dddddd, Rightarrow, shorten = 50, "\omega_{Y{,}X{,}Z}"]
			\arrow[lldd, "\chi_{Y\bigoplus X{,}Z}"description]
			\arrow[rrdd, "\chi_{Y{,}X\bigoplus Z}"description]
			\arrow[rruu, "F(1\bigoplus \beta_{X{,}Z})"description]
			&&&&FY\bigoplus F(Z\bigoplus X)
			\arrow[dddd, "1\bigoplus \chi_{Z{,}X}"]
			\\
			\\
			&&F(Y\bigoplus X)\bigoplus FZ\arrow[rrdddd, "\chi_{Y{,}X}\bigoplus 1"description]
			\arrow[ddd, Rightarrow, shorten = 10, "\Omega_{X{,}Y}\bigoplus 1"']
			&&&&FY\bigoplus F(X \bigoplus Z)
			\arrow[ddd,Rightarrow, shorten = 10, "1\bigoplus \Omega_{X{,}Z}"]
			\arrow[rruu, "1\bigoplus F\beta_{X{,}Z}"description]\arrow[lldddd, "1\bigoplus \chi_{X{,}Z}"description]
			\\
			\\
			FX\bigoplus FY\bigoplus FZ\arrow[rrrrdd, "\beta_{FX{,}FY}\bigoplus 1"']
			&&&&&&&&FY\bigoplus FZ \bigoplus FX
			\\
			&&{}&&&&{}
			\\
			&&&&FY\bigoplus FX \bigoplus FZ\arrow[rrrruu, "1\bigoplus\beta_{FX{,}FZ}"']
	\end{tikzcd}$$

	$$\begin{tikzcd}[font = \fontsize{9}{6}, column sep = 10]
	&&F(X\bigoplus Y\bigoplus Z)
	\arrow[dddd, "\chi_{X{,}Y\bigoplus Z}"description]
	\arrow[lldd, "\chi_{X\bigoplus Y{,}Z}"']
	\arrow[rrrr, "F\beta_{X{,}Y\bigoplus Z}"]
	&&{}\arrow[dddd, Rightarrow, shorten = 40, "\Omega_{X{,}Y\bigoplus Z}"]
	&&F(Y\bigoplus Z\bigoplus X)
	\arrow[dddd, "\chi_{Y\bigoplus Z{,}X}"description]
	\arrow[ddrr, "\chi_{Y{,}Z\bigoplus X}"]
	\\
	\\
	F(X\bigoplus Y)\bigoplus FZ
	\arrow[dd, Rightarrow, shift left = 10, shorten = 10, "\omega_{X{,}Y{,}Z}"]
	\arrow[dddd, "\chi_{X{,}Y}\bigoplus 1"']
	&&&&&&&&FY\bigoplus F(Z\bigoplus X)
	\arrow[dd, Rightarrow, shorten = 15, shift right =10, "\omega_{Y{,}Z{,}X}"']
	\arrow[dddd, "1\bigoplus \chi_{Z{,}X}"]
	\\
	\\
	{}&&FX\bigoplus F(Y\bigoplus Z)
	\arrow[lldd, "1\bigoplus \chi_{Y{,}Z}"description]
	\arrow[rrrr, "\beta_{FX{,}F(Y\bigoplus Z)}"]
	&&{}\arrow[dd, Rightarrow, shorten = 15, "\beta_{1{,}\chi_{Y{,}Z}}"]
	&&F(Y\bigoplus Z)\bigoplus FX
	\arrow[rrdd, "\chi_{Y{,}Z}\bigoplus 1"description]
	&&{}
	\\
	\\
	FX\bigoplus FY\bigoplus FZ\arrow[rrrrdd, "\beta_{FX{,}FY}\bigoplus 1"']
	\arrow[rrrrrrrr, "\beta_{FX{,}FY\bigoplus FZ}"']
	&&&&{}\arrow[dd, Rightarrow, shorten = 15,"r_{FY{,}FZ}^{FX}"]
	&&&&FY\bigoplus FZ \bigoplus FX
	\\
	&&{}&&&&{}
	\\
	&&&&FY\bigoplus FX \bigoplus FZ\arrow[rrrruu, "1\bigoplus\beta_{FX{,}FZ}"']
\end{tikzcd}$$

\noindent \emph{(left braiding law)}

	$$\begin{tikzcd}[font = \fontsize{9}{6}, column sep = 10]
	&&F(X\bigoplus Y\bigoplus Z)\arrow[dddd, Rightarrow, shorten = 20, "\chi_{1{,}\beta_{Y{,}Z}}"]
	\arrow[lldd, "\chi_{X{,}Y\bigoplus Z}"']
	\arrow[rrdd, "F(1\bigoplus \beta_{Y{,}Z})"description]
	\arrow[rrrr, "F\beta_{X\bigoplus Y{,} Z}"]
	&&{}\arrow[dd, Rightarrow, shift right = 2,shorten = 15, "Fl_{X{,}Y}^{Z}"]
	&&F(Z\bigoplus X\bigoplus Y)
	\arrow[dddd, Rightarrow, shorten =20, "\chi_{\beta_{X{,}Z}{,}1}"]
	\arrow[ddrr, "\chi_{Z\bigoplus X{,}Y}"]
	\\
	\\
	FX\bigoplus F(Y\bigoplus Z)
	\arrow[dddd, "1\bigoplus \chi_{Y{,}Z}"']
	\arrow[rrdd, "1\bigoplus F\beta_{Y{,}Z}"description]
	&&&&F(X\bigoplus Z\bigoplus Y)
	\arrow[dddddd, Rightarrow, shorten = 50, "\omega_{X{,}Z{,}Y}"]
	\arrow[lldd, "\chi_{X{,}Z\bigoplus Y}"description]
	\arrow[rrdd, "\chi_{X\bigoplus Z{,}Y}"description]
	\arrow[rruu, "F(\beta_{X{,}Z}\bigoplus 1)"description]
	&&&&F(Z\bigoplus Y)\bigoplus FX
	\arrow[dddd, "\chi_{Z{,}Y}\bigoplus 1"]
	\\
	\\
	&&FX\bigoplus F(Z\bigoplus Y)\arrow[rrdddd, "1\bigoplus \chi_{Z{,}Y}"description]
	\arrow[ddd, Rightarrow, shorten = 10, "1\bigoplus \Omega_{Y{,}Z}"']
	&&&&F(X \bigoplus Z)\bigoplus FY
	\arrow[ddd,Rightarrow, shorten = 10, "\Omega_{X{,}Z}\bigoplus 1"]
	\arrow[rruu, " F\beta_{X{,}Z}\bigoplus 1"description]
	\arrow[lldddd, "\chi_{X{,}Z}\bigoplus 1"description]
	\\
	\\
	FX\bigoplus FY\bigoplus FZ\arrow[rrrrdd, "1\bigoplus\beta_{FY{,}FZ}"']
	&&&&&&&&FZ\bigoplus FX \bigoplus FY
	\\
	&&{}&&&&{}
	\\
	&&&&FX\bigoplus FZ \bigoplus FY\arrow[rrrruu, "\beta_{FX{,}FZ}\bigoplus 1"']
\end{tikzcd}$$

$$\begin{tikzcd}[font = \fontsize{9}{6}, column sep = 10]
	&&F(X\bigoplus Y\bigoplus Z)\arrow[dddd, "\chi_{X\bigoplus Y{,}Z}"description]
	\arrow[lldd, "\chi_{X{,}Y\bigoplus Z}"']
	\arrow[rrrr, "F\beta_{X\bigoplus Y{,} Z}"]
	&&{}\arrow[dddd, Rightarrow, shorten = 40, "\Omega_{X\bigoplus Y{,}Z}"]
	&&F(Z\bigoplus X\bigoplus Y)
	\arrow[dddd, "\chi_{Z{,}X\bigoplus Y}"description]
	\arrow[ddrr, "\chi_{Z\bigoplus X{,}Y}"]
	\\
	\\
	FX\bigoplus F(Y\bigoplus Z)
	\arrow[dd, Rightarrow, shorten = 15, shift left = 10, "\omega_{X{,}Y{,}Z}"]
	\arrow[dddd, "1\bigoplus \chi_{Y{,}Z}"']
	&&&&&&&&F(Z\bigoplus Y)\bigoplus FX
	\arrow[dd, Rightarrow, shorten = 15,shift right = 10, "\omega_{Z{,}X{,}Y}"']
	\arrow[dddd, "\chi_{Z{,}Y}\bigoplus 1"]
	\\
	\\
	{}&&F(X\bigoplus Y)\bigoplus FZ
	\arrow[lldd, "\chi_{X{,}Y}\bigoplus 1"description]
	\arrow[rrrr, "\beta_{F(X\bigoplus Y){,}FZ}"]
	&&{}\arrow[dd, Rightarrow, shorten = 15, "\beta_{\chi_{X{,}Y}{,}1}"]
	&&FZ\bigoplus F(X\bigoplus Y)\arrow[rrdd, "1\bigoplus \chi_{X{,}Y}"description]
	&&{}
	\\
	\\
	FX\bigoplus FY\bigoplus FZ
	\arrow[rrrrrrrr, "\beta_{FX\bigoplus FY{,}FZ}"]
	\arrow[rrrrdd, "1\bigoplus\beta_{FY{,}FZ}"']
	&&&&{}\arrow[dd, Rightarrow, shorten = 15, "l_{FX{,}FY}^{FZ}"]
	&&&&FZ\bigoplus FX \bigoplus FY
	\\
	&&{}&&&&{}
	\\
	&&&&FX\bigoplus FZ \bigoplus FY\arrow[rrrruu, "\beta_{FX{,}FZ}\bigoplus 1"']
\end{tikzcd}$$

\end{definition}

\begin{definition}\label{Definition braided pseudonatural transformation}
	A \emph{braided monoidal pseudonatural transformation} $\left(\phi, \phi_{0}, \phi_{2}\right): \left(F, \Omega\right) \rightarrow \left(G, \Theta\right)$ is a monoidal pseudonatural transformation  between underlying opmonoidal $2$-functors $\left(\phi, \phi_{0}, \phi_{2}\right): F \rightarrow G$ which satisfies the following equation of pasting diagrams in $\mathcal{B}$ for every $X, Y \in \mathcal{A}$.
	
	$$\begin{tikzcd}[font=\fontsize{9}{6}]
		F\left(X \bigoplus Y\right) \arrow[rrrr, "F\beta_{X{,}Y}"]
		\arrow[dddd, "\phi_{X\bigoplus Y}"']
		\arrow[rrd, "\chi_{X{,}Y}"] &&&{}\arrow[d, Rightarrow, "\Omega_{X{,} Y}"]& F\left(Y\bigotimes X\right)\arrow[rrd, "\chi_{Y{,}X}"]
		\\
		&{}\arrow[ddd, Rightarrow, shorten = 10, "\phi_{2{,}X{,}Y}"']
		&FX\bigoplus FY \arrow[rrrr, "\beta_{FX{,}FY}"]\arrow[dd, "\phi_{X}\bigoplus FY"description]
		&{}&{}\arrow[dd, Rightarrow, shift right = 22, shorten = 10, "\beta_{\phi_{X}{,}FY}"]&&
		FY\bigoplus FX
		\arrow[lldd, "FY \bigoplus \phi_{X}"']
		\arrow[dd, "\phi_{Y}\bigoplus FX"]
		\\
		&&&&&{}\arrow[dd, Rightarrow, shorten = 10, "{(\phi_{X})}_{(\phi_{Y})}"description]
		\\
		&&GX\bigoplus FY
		\arrow[rr, "\beta_{GX{,}FY}"]
		\arrow[dd, "GX \bigoplus \phi_{Y}"description]
		&&FY \bigoplus GX\arrow[dd, Rightarrow, shorten = 10, shift right = 22, "\beta_{GX{,}\phi_{Y}}"]
		\arrow[rrdd, "\phi_{Y}\bigoplus GX"']
		&&GY \bigoplus FX\arrow[dd, "GY \bigoplus \phi_{X}"]
		\\
		G\left(X \bigoplus Y\right)\arrow[rrd, "\chi_{X{,}Y}"']&{}&&&&{}
		\\
		&&GX\bigoplus GY\arrow[rrrr, "\beta_{GY{,}GY}"'] &&{}&& GY \bigoplus GX
	\end{tikzcd}$$
	
	$$= \begin{tikzcd}[font=\fontsize{9}{6}]
		F\left(X \bigoplus Y\right) \arrow[rrrr, "F\beta_{X{,}Y}"]
		\arrow[dddd, "\phi_{X\bigoplus Y}"'] &&{}\arrow[dddd, Rightarrow, shorten = 35, "\phi_{\beta_{X{,}Y}}"]
		&{}& F\left(Y\bigoplus X\right)
		\arrow[dddd, "\phi_{Y\bigoplus X}"']
		\arrow[rrd, "\chi_{Y{,}X}"]
		\\
		&{}
		&{}
		&{}&{}&{}\arrow[ddd, Rightarrow, shorten = 10, "\phi_{2{,}Y{,}X}"]
		&
		FY\bigoplus FX
		\arrow[dd, "\phi_{Y}\bigoplus FX"]
		\\
		&&&&&{}
		\\
		&&{}
		&&{}
		&&GY \bigoplus FX\arrow[dd, "GY \bigoplus \phi_{X}"]
		\\
		G\left(X \bigoplus Y\right)\arrow[rrd, "\chi_{X{,}Y}"']\arrow[rrrr, "G\beta_{X{,}Y}"]
		&{}&{}\arrow[d, Rightarrow, shift left = 20, "\Theta_{X{,}Y}"]
		&&G\left(Y\bigoplus X\right)
		\arrow[rrd, "\chi_{X{,}Y}"']
		&{}
		\\
		&&GX\bigoplus GY\arrow[rrrr, "\beta_{GY{,}GY}"'] &&{}&& GY \bigoplus GX
	\end{tikzcd}$$
	
\end{definition}

\begin{definition}\label{definition sylleptic opmonoidal functor}
	If $(\A, \bigoplus , I)$ and $(\B, \bigoplus I)$ from Definition \ref{Braided opmonoidal 2-functors} have syllepsis $\Sigma$, then $F$ is called \emph{sylleptic} if for every pair of objects $X, Y \in\A$, the following diagram commutes.
	
	$$\begin{tikzcd}[font=\fontsize{9}{6}, column sep = 18]
		F(X\bigoplus Y)\arrow[rr, "F\beta_{X{,}Y}"]\arrow[dd, "\chi_{X{,}Y}"']
		&{}\arrow[dd, Rightarrow, shorten = 15, "\Omega_{X{,}Y}"]
		&F(Y\bigoplus X)\arrow[rd, "F\beta_{Y{,}X}"]
		\arrow[dd, "\chi_{Y{,}X}"description]
		\\
		&&&F(X\bigoplus Y)\arrow[dd, "\chi_{X{,}Y}"]\arrow[d, Rightarrow, shorten = 5, shift right = 8, "\Omega_{Y{,}X}"']
		\\
		FX\bigoplus FY
		\arrow[rrrd, bend right = 20, "1_{FX\bigoplus FY}"']
		\arrow[rr, "\beta_{FX{,}FY}"]
		&{}&FY\bigoplus FX
		\arrow[d, Rightarrow, shorten = 5, shift right = 3, "\Sigma_{FX{,}FY}"']
		\arrow[rd, "\beta_{FY{,}FX}"description]&{}
		\\
		&&{}&FX\bigoplus FY&{}
	\end{tikzcd}=\begin{tikzcd}[font=\fontsize{9}{6}, column sep = 18]
	F(X\bigoplus Y)
	\arrow[rrrd, bend right = 20, "1_{F(X\bigoplus Y)}"']
	\arrow[rr, "F\beta_{X{,}Y}"]
	\arrow[dd, "\chi_{X{,}Y}"']
	&{}
	&F(Y\bigoplus X)
	\arrow[d, Rightarrow, shorten = 5, shift right = 3, "F\Sigma_{X{,}Y}"']
	\arrow[rd, "F\beta_{Y{,}X}"]
	\\
	&&{}&F(X\bigoplus Y)\arrow[dd, "\chi_{X{,}Y}"]
	\\
	FX\bigoplus FY
	\arrow[rrrd, bend right = 20, "1_{FX\bigoplus FY}"']
	&{}&=
	\\
	&&{}&FX\bigoplus FY
\end{tikzcd}$$
	
\end{definition}

\begin{remark}\label{Remark braiding and syllapsis comparison to literature}
	The axioms we gave in Definition \ref{Braided opmonoidal 2-functors} for braided opmonoidal $2$-functors are dual to those listed in Definition 14 of \cite{monoidal bicategories and hopf algebroids}, which are axioms for the monoidal variants in which compositors are of the form $\chi_{X, Y}: FX \oplus FY \to F(X \oplus Y)$. The axiom given at the top of page 126 in \cite{monoidal bicategories and hopf algebroids} for braided opmonoidal pseudonatural transformations is also a dual version, but also suppresses notation for interchanger ${\left(\phi_{X}\right)}_{\left(\phi_{Y}\right)}$. It is otherwise analogous to Definition \ref{Definition braided pseudonatural transformation}. Finally, the axiom in Definition \ref{definition sylleptic opmonoidal functor} specialises Definition 1.5 of \cite{Gurski Osorno Infinite Loop Spaces and coherence for symmetric monoidal bicategories} to the setting where the monoidal bicategories are symmetric $\Gr$-monoids, and the opmonoidal pseudofunctor is a $2$-functor. We have restated the axioms in detail rather than simply describing adaptations to existing literature to assist the reader in verifying details in the proof of Proposition \ref{Proposition Braidings Syllapses and Symmetries as data in Psmnd Gray}.
\end{remark}

\begin{proposition}\label{Tricategories of E n bicategories}
	There are tricategories
	\begin{enumerate}
		\item $\mathbf{opBrMonBicat}$, of braided monoidal bicategories, braided opmonoidal pseudofunctors, braided monoidal transformations and monoidal modifications.
		\item $\mathbf{opSylMonBicat}$, of sylleptic monoidal bicategories, sylleptic opmonoidal pseudofunctors, braided monoidal transformations and monoidal modifications.
		\item $\mathbf{opSymMonBicat}$, of symmetric monoidal bicategories, sylleptic opmonoidal pseudofunctors, braided monoidal transformations and monoidal modifications..
	\end{enumerate} 
\end{proposition}

\begin{proof}
	See Sections 2.3 and 2.4 of \cite{Schommer Piers Classification of 2D ETFTs} for sub-tricategory structures involving strong (braided resp. sylleptic) monoidal pseudofunctors, and Lemma 3.5 of \cite{Gurski Johnson Osorno K Theory for 2-categories} for the symmetric case. It is easy to observe that the same constructions as described in Section 2.4 of \cite{Schommer Piers Classification of 2D ETFTs} are also well-defined for (braided resp. sylleptic) opmonoidal pseudofunctors, and braided opmonoidal transformations.
\end{proof}

\begin{example}\label{Example opmonoidal pseudomonad structure on Prof}
	Recall the bicategory $\mathbf{Prof}$ of categories, profunctors, and natural transformations. The usual cartesian product of categories extends to a symmetric monoidal bicategory structure on $\mathbf{Prof}$, and the free symmetric monoidal category $2$-monad on $\Cat$ also extends to a symmetric opmonoidal pseudomonad on $(\mathbf{Prof}, \times, \mathbf{1})$. This opmonoidal pseudomonad is moreover symmetric, in the sense that it is a pseudomonad in the bicategory $\mathbf{opSymMonBicat}$ of Proposition \ref{Tricategories of E n bicategories}. See Theorem 8.2 of \cite{Fiore Gambino Hyland Monoidal Bicategories differential linear logic and analytic functors} for details on the compatibility between the various structures on $\mathbf{Prof}$.
\end{example}

\noindent As justified by coherence arguments similar to what we discussed in Subsection \ref{Subsection opmonoidal pseudomonads as pseudomonads in the tricategory of monoidal bicategories}, we restrict our attention from pseudomonads internal to the tricategories of Proposition \ref{Tricategories of E n bicategories} to simpler, stricter, notions.

\begin{definition}\label{Definition braided opmonoidal pseudomonad}
	\hspace{1mm}
	\begin{enumerate}
		\item A \emph{braided opmonoidal pseudomonad} consists of the following data.
		\begin{itemize}
			\item A braided $\mathbf{Gray}$-monoid $\left(\mathcal{A}, \bigoplus, I, \beta, r, l\right)$,
			\item A braided opmonoidal $2$-functor $\left(S, \chi, \iota, \gamma, \omega, \delta, \Omega\right): \left(\mathcal{A}, \bigoplus, I, \beta, r, l\right) \rightarrow \left(\mathcal{A}, \bigoplus, I, \beta\right) $, which will be abbreviated as $S$ whenever the rest of the data is clear from context,
			\item A braided pseudonatural transformations $\left(\eta, \eta_{0}, \eta_{2}\right): 1 \rightarrow S$ and $\left(\mu, \mu_{0}, \mu_{2}\right): S^{2} \rightarrow S$
			\item invertible monoidal modifications $\mathbf{a}$, $\mathbf{l}$ and $\mathbf{r}$ which satisfy the axioms needed for $\left(\mathcal{A}, S, \eta, \mu, \mathbf{l}, \mathbf{a}, \mathbf{r}\right)$ to be a pseudomonad.
		\end{itemize}
		\item A \emph{sylleptic opmonoidal pseudomonad} consists of the following data.
		\begin{itemize}
			\item A sylleptic $\mathbf{Gray}$-monoid $\left(\mathcal{A}, \bigoplus, I, \beta, r, l, \Sigma\right)$,
			\item A braided opmonoidal pseudomonad on $\left(\mathcal{A}, \bigoplus, I, \beta,r, l\right)$, for which the underlying braided opmonoidal $2$-functor is moreover sylleptic.
		\end{itemize}
	\item A \emph{symmetric opmonoidal pseudomonad} is a sylleptic opmonoidal pseudomonad whose underlying sylleptic $\Gr$-monoid is moreover symmetric.
	\end{enumerate}
\end{definition}

\begin{remark}\label{Remark commutativity}
	A different notion of compatibility between monads and symmetric monoidal structures in the two-dimensional setting has been considered in \cite{Hyland Power Pseudo commutative monads and pseudo closed categories}. This is known as \emph{pseudo-commutativity}, and is expressed in terms of compatibility between left and right strengths, which are defined in terms of enrichment over $\Cat$. Theorem 7 of \cite{Hyland Power Pseudo commutative monads and pseudo closed categories} relates the notion of symmetric pseudo-commutativity to what they call a \emph{symmetric pseudo-monoidal $2$-monad}, however a definition of the latter notion is left `for another occasion'. Part (3) of Definition \ref{Definition braided opmonoidal pseudomonad} can be adapted to a precise definition by asking the underlying pseudomonad to be a $2$-monad, and the pseudonatural transformations $\chi$ and $\iota$ to be adjoint equivalences.
\end{remark}

\noindent In Proposition \ref{Proposition Braidings Syllapses and Symmetries as data in Psmnd Gray}, to follow, we consider semi-strict opmonoidal pseudomonads which are moreover either braided or sylleptic, and describe corresponding data extending the pseudomonoid in $\mathbf{Psmnd}(\Gr)_\text{ast}$ of Theorem \ref{theorem semi-strict opmonoidal pseudomonads are pseudomonoids in the Gray monoid of pseudomonads}. These are analogous to data up to which the symmetric pseudomonoid axioms described in Definitions 13 and 17 of \cite{monoidal bicategories and hopf algebroids} hold.

\begin{proposition}\label{Proposition Braidings Syllapses and Symmetries as data in Psmnd Gray}
	Consider a braided opmonoidal pseudomonad as in Definition \ref{Definition braided opmonoidal pseudomonad}. For a pseudomonad $(\A, S)$, write ${(\A, S)}^n$ for the $n$-fold tensor product as per Lemma \ref{Gray tensor of two pseudomonads}. Punctuate the $\mathbf{Gray}$ tensor product with $``."$ wherever necessary.
	
	\begin{enumerate}
		\item The pair $\left(\beta, \Omega\right)$ form a pseudomonad transformation as depicted below.
		
		$$\begin{tikzcd}[font = \fontsize{9}{6}]
			{(\A{,}S)}^{2}
			\arrow[rr, "(\tau_{\A{,}\A}{,}1)"]
			\arrow[rd, "(\bigoplus{,}\chi)"']
			&{}\arrow[d, Rightarrow, shorten = 5, "(\beta{,}\Omega)"]
			&{\left(\A{,}S\right)}^{2}\arrow[ld, "(\bigoplus{,}\chi)"]
			\\
			&\left(\A{,} S\right)
		\end{tikzcd}$$
		
		\item The modification $r$ forms a pseudomonad modification as displayed below.
		
		$$\begin{tikzcd}[font = \fontsize{9}{6}, column sep =18]
			&(\A{,}S)^3
			\arrow[rrrr, "(\tau_{\A{,}\A^2}{,}1)"]
			\arrow[ldd, "(\bigoplus{,}\chi). 1"']
			\arrow[rd, "1. (\bigoplus{,}\chi)"']
			&&&&(\A{,}S)^3
			\arrow[rdd, "1.(\bigoplus{,}\chi)"]
			\arrow[ld, "(\bigoplus{,}\chi).1 "]
			\\
			&&(\A{,}S)^2
			\arrow[rr, Rightarrow, shorten = 15, shift right = 10, "(\beta{,}\Omega)"]
			\arrow[rr, equal,shorten = 30, shift left = 7]
			\arrow[rr, "(\tau_{\A{,}\A}{,}1)"]
			\arrow[rdd, "(\bigoplus{,}\chi)"']
			&&(\A{,}S)^2
			\arrow[ldd, "(\bigoplus{,}\chi)"]
			\\
			(\A{,}S)^2
			\arrow[rr, shorten = 15, Rightarrow, "(1{,}\omega)"]
			\arrow[rrrd, "(\bigoplus{,}\chi)"']
			&&{}&&{}\arrow[rr, shorten = 15, Rightarrow, "(1{,}\omega)"]
			&&(\A{,}S)^2\arrow[llld, "(\bigoplus{,}\chi)"]
			\\
			&&&(\A{,}S)
		\end{tikzcd}$$
	$$\Rrightarrow_{r}
	\begin{tikzcd}[font = \fontsize{9}{6}, column sep =25]
		&(\A{,}S)^3
		\arrow[rrrr, equal, shorten = 70, shift right = 5]
		\arrow[rrd, "(\tau_{\A{,}\A}{,}1).1"']
		\arrow[rrrr, "(\tau_{\A{,}\A^2}{,}1)"]
		\arrow[ldd, "(\bigoplus{,}\chi). 1"']
		&&&&(\A{,}S)^3
		\arrow[rdd, "1.(\bigoplus{,}\chi)"]
		\\
		{}\arrow[rrr, Rightarrow, shorten = 60, shift right = 5, "(\beta{,}\Omega).1"]
		&&&(\A{,}S)^{3}
		\arrow[rru, "1.(\tau_{\A{,}\A}{,}1)"']
		\arrow[llld, near end,"1. (\bigoplus{,}\chi)"]
		\arrow[rrrd, near end,"1. (\bigoplus{,}\chi)"']
		\arrow[rrr, Rightarrow, shorten = 60, shift right = 5,"1.(\beta{,}\Omega)"]
		&&&{}
		\\
		(\A{,}S)^2
		\arrow[rrrrrr, shorten = 130, Rightarrow, "(1{,}\omega)"]
		\arrow[rrrd, "(\bigoplus{,}\chi)"']
		&&{}&&{}
		&&(\A{,}S)^2\arrow[llld, "(\bigoplus{,}\chi)"]
		\\
		&&&(\A{,}S)
	\end{tikzcd}$$

		\item The modification $l$ forms a pseudomonad modification as displayed below.

		$$\begin{tikzcd}[font = \fontsize{9}{6}, column sep =18]
			&(\A{,}S)^3
			\arrow[rrrr, "(\tau_{\A^2{,}\A}{,}1)"]
			\arrow[ldd, "1.(\bigoplus{,}\chi)"']
			\arrow[rd, "(\bigoplus{,}\chi).1"']
			&&&&(\A{,}S)^3
			\arrow[rdd, "(\bigoplus{,}\chi).1"]
			\arrow[ld, "1.(\bigoplus{,}\chi)"]
			\\
			&&(\A{,}S)^2
			\arrow[rr, Rightarrow, shorten = 15, shift right = 10, "(\beta{,}\Omega)"]
			\arrow[rr, equal,shorten = 30, shift left = 7]
			\arrow[rr, "(\tau_{\A{,}\A}{,}1)"]
			\arrow[rdd, "(\bigoplus{,}\chi)"']
			&&(\A{,}S)^2
			\arrow[ldd, "(\bigoplus{,}\chi)"]
			\\
			(\A{,}S)^2
			\arrow[rr, shorten = 15, Rightarrow, "{(1{,}\omega)}^{-1}"]
			\arrow[rrrd, "(\bigoplus{,}\chi)"']
			&&{}&&{}\arrow[rr, shorten = 15, Rightarrow, "{(1{,}\omega)}^{-1}"]
			&&(\A{,}S)^2\arrow[llld, "(\bigoplus{,}\chi)"]
			\\
			&&&(\A{,}S)
		\end{tikzcd}$$
		$$\Rrightarrow_{l}
		\begin{tikzcd}[font = \fontsize{9}{6}, column sep =25]
			&(\A{,}S)^3
			\arrow[rrrr, equal, shorten = 70, shift right = 5]
			\arrow[rrd, "1.(\tau_{\A{,}\A}{,}1)"']
			\arrow[rrrr, "(\tau_{\A^2{,}\A}{,}1)"]
			\arrow[ldd, "1.(\bigoplus{,}\chi)"']
			&&&&(\A{,}S)^3
			\arrow[rdd, "(\bigoplus{,}\chi).1"]
			\\
			{}\arrow[rrr, Rightarrow, shorten = 60, shift right = 5, "1.(\beta{,}\Omega)"]
			&&&(\A{,}S)^{3}
			\arrow[rru, "(\tau_{\A{,}\A}{,}1).1"']
			\arrow[llld, near end,"(\bigoplus{,}\chi).1"]
			\arrow[rrrd, near end,"(\bigoplus{,}\chi).1"']
			\arrow[rrr, Rightarrow, shorten = 60, shift right = 5,"(\beta{,}\Omega).1"]
			&&&{}
			\\
			(\A{,}S)^2
			\arrow[rrrrrr, shorten = 130, Rightarrow, "{(1{,}\omega)}^{-1}"]
			\arrow[rrrd, "(\bigoplus{,}\chi)"']
			&&{}&&{}
			&&(\A{,}S)^2\arrow[llld, "(\bigoplus{,}\chi)"]
			\\
			&&&(\A{,}S)
		\end{tikzcd}$$

		\item If the braided $\mathbf{Gray}$-monoid moreover has a syllepsis $\Sigma$ and $S: \A \to \A$ is sylleptic, then $\Sigma$ is a pseudomonad modification, as displayed below.
	
\end{enumerate}
		$$\begin{tikzcd}[font = \fontsize{9}{6}, column sep =18]
			(\A{,}S)^2\arrow[rrrr, "(1{,}1)^{2}"]
			\arrow[rrdd, "(\tau{,}1)"]
			\arrow[rrdddd, bend right = 30,"(\bigoplus {,}\chi)"'] 
			&&{}&&(\A{,}S)^2
			\arrow[lldddd,bend left = 30, "(\bigoplus{.}\chi)"]
			\\
			&{}&=&{}
			\\
			&&(\A{,}S)^2
			\arrow[dd, "(\bigoplus{.}\chi)"]
			\arrow[rruu, "(\tau{.}1)"]
			&&&\Rrightarrow_{\Sigma}&{}
			\\
			&{}\arrow[uu, Rightarrow, shorten = 15, "(\beta{,}\Omega)"]
			&&{}\arrow[uu, Rightarrow, shorten = 15, "(\beta{,}\Omega)"]
			\\
			&&(\A{,}S)
		\end{tikzcd}\begin{tikzcd}[font = \fontsize{9}{6}, column sep = 18]
		(\A{,}S)^2\arrow[rrrr, "(1{,}1)^{2}"]
		\arrow[rrdd, "(\tau{,}1)"]
		\arrow[rrdddd, bend right = 30,"(\bigoplus {,}\chi)"'] 
		&&{}&&(\A{,}S)^2
		\arrow[lldddd,bend left = 30, "(\bigoplus{.}\chi)"]
		\\
		&{}&=&{}
		\\
		&&(\A{,}S)^2
		\arrow[dd, "(\bigoplus{.}\chi)"]
		\arrow[rruu, "(\tau{.}1)"]
		\\
		&{}\arrow[uu, Rightarrow, shorten = 15, "(\beta{,}\Omega)"]
		&&{}\arrow[uu, Rightarrow, shorten = 15, "(\beta{,}\Omega)"]
		\\
		&&(\A{,}S)
	\end{tikzcd}$$
		
\end{proposition}

\begin{proof}
	The unit and multiplication laws for the pseudomonad transformation in part (1) say precisely that $\left(\eta, \eta_{0}, \eta_{2}\right): 1_\mathcal{A} \rightarrow S$ and $\left(\mu, \mu_{0}, \mu_{2}\right): S^2 \rightarrow S$ are braided pseudonatural transformations. The braided opmonoidal $2$-functoriality laws say precisely that $r$ and $l$ are pseudomonad modifications, as needed for parts (2) and (3). Finally, the sylleptic opmonoidal functoriality law says precisely that $\Sigma$ is a $3$-cell of pseudomonads, as needed for part (4).
\end{proof}

 \noindent We do not pursue a symmetric monoidal $\Gr$-category structure on $\PsG$ in which braided (resp. sylleptic, symmetric) opmonoidal pseudomonads can be seen as certain three-dimensional monoids. Nonetheless, each of the pieces of data described in Proposition \ref{Proposition Braidings Syllapses and Symmetries as data in Psmnd Gray} gives analogous data between $2$-categories of pseudoalgebras, via $\mathbf{EM}: \PsG \to \mathbf{Gray}^\mathbf{2}$. In the next two Subsections we directly show that these pseudonatural equivalences and invertible modifications assemble into a symmetric monoidal bicategory structure on $\A^S$.

\subsection{Symmetric opmonoidality of $\mathbf{EM}$}\label{Subsection symmetric opmonoidality of EM}

\noindent The proof of Proposition \ref{Symmetriser's component 2-natural transformation}, to follow, is given in Appendix \ref{Appendix Proposition symmetriser's component 2-nat trans}.

\begin{proposition}\label{Symmetriser's component 2-natural transformation}
	Let $\left(\mathcal{A}, S\right)$ and $\left(\mathcal{B}, T\right)$ be pseudomonads. Then there is a $2$-natural isomorphism as depicted below, whose component on an object $\left(\mathbf{X}, \mathbf{Y}\right) \in \mathcal{A}^{S} \otimes \mathcal{B}^{T}$ is the tight pseudomorphism in ${\left(\mathcal{A}\otimes \mathcal{B}\right)}^{S\otimes T}$ with $2$-cell component given by the interchanger $\left(x, y\right)$.
	
	$$\begin{tikzcd}[font=\fontsize{9}{6}]
		\mathcal{A}^{S} \otimes \mathcal{B}^{T}
		\arrow[rr, "M_{S{,}T}"]
		\arrow[dd, "\tau"']
		&{}\arrow[dd, Rightarrow, shorten = 10, "\Xi_{S{,}T}"]
		& {\left(\mathcal{A}\otimes \mathcal{B}\right)}^{S\otimes T}
		\arrow[dd, "\overline{\tau}"]
		\\
		\\
		\mathcal{B}^{T} \otimes \mathcal{A}^{S}
		\arrow[rr, "M_{T{,}S}"']
		&{}& {\left(\mathcal{B}\otimes \mathcal{A}\right)}^{T\otimes S}
	\end{tikzcd}$$
\end{proposition}

\begin{proposition}\label{Proposition assignment of pseudomonads to symmetriser for EM is a modification}
	The assignment of a pair of pseudomonads $\left((\A, S), (\B, T)\right)$ to the $2$-natural isomorphism $\Xi_{S, T}$ of Proposition \ref{Symmetriser's component 2-natural transformation} defines an invertible modification. 
\end{proposition}

\begin{proof}
	It suffices to verify the modification condition as either of the pseudomonads vary along a morphism of pseudomonads. The modification condition on such data corresponds to an equation of $2$-natural transformations with tight components, which in turn corresponds to an equation between $2$-cells in $\B \otimes \A$. As the second variable varies along $(H, h):(\B, T) \to (\B', T')$, the two tight pseudomorphisms which need to be equal both have $2$-cell components given as displayed below left. Similarly, as the first variable varies along $(G, g):(\A, S) \to (\A', S')$, the two tight pseudomorphisms which need to be equal both have $2$-cell components given as displayed below right. This completes the proof.
	
	$$\begin{tikzcd}[column sep = 16, font=\fontsize{9}{6}]
		\left(T'HY{,}SX\right) \arrow[rr, "\left(h_{Y}{,} 1\right)"]
		\arrow[dd, "\left(1{,}x\right)"']
		&{}\arrow[dd, Rightarrow, shorten = 15, "\left({h}_{Y}{,}x\right)"]
		& \left(HTY{,}SX\right) \arrow[dd, "\left(1{,}x\right)"] \arrow[rr, "\left(Hy{,}1\right)"]
		&{}\arrow[dd, Rightarrow, shorten = 15, "(Hy{,}x)"]
		& \left(HY{,} SX\right) \arrow[dd, "\left(1{,}x\right)"]
		\\
		\\
		\left(T'HY{,}X\right) \arrow[rr, "\left(h_{Y}{,}1\right)"']
		&{}& \left(HTY{,}X\right) \arrow[rr, "\left(Hy{,}1\right)"'] 
		&{}& \left(HY{,}X\right)&{}
	\end{tikzcd}\begin{tikzcd}[column sep = 16, font=\fontsize{9}{6}]
	\left(TY{,}S'GX\right) \arrow[rr, "\left(1{,}g_{X}\right)"]
	\arrow[dd, "\left(y{,}1\right)"']
	&{}\arrow[dd, Rightarrow, shorten = 15, "\left(y{,}g_{X}\right)"]
	& \left(TY{,} GSX\right) \arrow[dd, "\left(y{,}1\right)"]\arrow[rr, "(1{,}Gx)"]
	&{}\arrow[dd, Rightarrow, shorten = 15, "(y{,}Gx)"]
	&(TY{,}GX)\arrow[dd, "(y{,}1)"]
	\\
	\\
	\left(Y{,}S'GX\right) \arrow[rr, "\left(1{,}g_{X}\right)"']
	&{}& \left(Y{,}GSX\right) \arrow[rr, "\left(1{,}Gx\right)"'] &{}& \left(Y{,}GX\right)
\end{tikzcd}$$
\end{proof}

\begin{theorem}\label{theorem EM is symmetric monoidal}
	The lax monoidal $2$-functor $\mathbf{EM}: \PsG_\text{ast} \to \Gr_\text{2-nat}$ of Theorem \ref{Theorem EM is normal lax Gray functor} is symmetric, when equipped with the invertible modification of Proposition \ref{Proposition assignment of pseudomonads to symmetriser for EM is a modification}.
\end{theorem}

\begin{proof}
	We explain how the left and right braided lax monoidal $2$-functoriality conditions follow from the hexagon equations for the $\mathbf{Gray}$ tensor product, and the sylleptic monoidal $2$-functoriality equation follows from symmetry of the $\mathbf{Gray}$ tensor product.
	\\
	\\
	\noindent Consider first the axioms required for $\mathbf{EM}: \mathbf{Psmnd}(\Gr)_{\text{ast}} \to \Gr_{\text{2-nat}}^\mathbf{2}$ to be braided lax monoidal $2$-functor. These are dual to the axioms listed in Definition \ref{Braided opmonoidal 2-functors}. We explain how most of the $2$-cells appearing in the pasting diagrams in these axioms are actually identities.
	
	\begin{itemize}
		\item $\mathbf{Psmnd}(\Gr)_{\text{ast}}$ and $\Gr_{\text{2-nat}}^\mathbf{2}$ both satisfy the left and right braid equations on the nose, rather than up to invertible $2$-cells of the form $r_{Y, Z}^{X}$ or $l_{X, Y}^{Z}$.
		\item The underlying lax monoidal $2$-functor $\mathbf{EM}$ is strictly associative, by Proposition \ref{associativity and unitality of EM} part (3), and hence all instances of $\omega$ are identities.
		\item The braiding of $\Gr_{\text{2-nat}}^\mathbf{2}$ is strictly natural, and as such the instances of $\beta_{\chi, 1}$ and $\beta_{1, \chi}$ are identities.
		\item The compositor $\left((\A, S), (\B, T)\right) \mapsto M_{(\A, S), (\B, T)}$ is strictly natural as $(\A, S)$ varies, by Proposition \ref{Proposition M is natural in the first component}, and also strictly natural as $(\B, T)$ varies along morphisms of pseudomonads whose underlying pseudonatural transformations are invertible icons, by Proposition \ref{Proposition 2-cell component of compositor as the second pseudomonad varies} part (2). Since the braiding isomorphisms in $\mathbf{Psmnd}(\Gr)_{\text{ast}}$ are of this form, as per Lemma \ref{symmetry for Gray tensor product of pseudomonads}, it follows that all instances of $\chi_{1, \beta}$ and $\chi_{\beta, 1}$ are identities.
	\end{itemize}

\noindent The data that remain are of the form $\Omega_{-, ?}$. By Proposition \ref{Symmetriser's component 2-natural transformation}, it follows that the left and right braiding axioms for $\mathbf{EM}$ reduce to the equalities of pasting diagrams displayed below. But these both correspond to the left and right braid relations for the $\mathbf{Gray}$-tensor product, which indeed hold strictly as the below are also relations in the presentation of $\A \otimes \B\otimes \mathcal{C}$ discussed in Remark \ref{Remark Gray tensor product}.

$$\begin{tikzcd}[font=\fontsize{9}{6}]
	(SX{,}TY{,}RZ)
	\arrow[rr, "(x{,}1{,}1)"]
	\arrow[dd, "(1{,}y{,}1)"']
	&{}\arrow[dd, Rightarrow, shorten = 15,"(x{,}y{,}1)"]
	&(X{,}TY{,}RZ)\arrow[dd, "(1{,}y{,}1)"]
	\\
	\\
	(SX{,}Y{,}RZ)
	\arrow[rr, "(x{,}1{,}1)"]
	\arrow[dd, "(1{,}1{,}z)"']
	&{}\arrow[dd, Rightarrow, shorten = 15,"(x{,}1{,}z)"]
	&(X{,}Y{,}RZ)\arrow[dd, "(1{,}1{,}z)"]
	&=&{}
	\\
	\\
	(SX{,}Y{,}Z)
	\arrow[rr, "(x{,}1{,}1)"']
	&{}&(X{,}Y{,}Z)
\end{tikzcd}\begin{tikzcd}[font=\fontsize{9}{6}]
(SX{,}TY{,}RZ)
\arrow[rr, "(x{,}1{,}1)"]
\arrow[dddd, "(1{,}1{,}z)\circ(1{,}y{,}1)"description]
&{}\arrow[dddd, Rightarrow, shorten = 20,"\left(x{,}(1{,}z)\circ(y{,}1)\right)"description]
&(X{,}TY{,}RZ)\arrow[dddd, "(1{,}1{,}z)\circ(1{,}y{,}1)"description]
\\
\\
\\
\\
(SX{,}Y{,}Z)
\arrow[rr, "(x{,}1{,}1)"']
&{}&(X{,}Y{,}Z)
\end{tikzcd}$$

$$\begin{tikzcd}[font=\fontsize{9}{6}]
	(SX{,}TY{,}RZ)
	\arrow[rr, "(1{,}1{,}z)"]
	\arrow[dd, "(1{,}y{,}1)"']
	&{}\arrow[dd, Rightarrow, shorten = 15,"(1{,}y{,}z)"]
	&(X{,}TY{,}Z)\arrow[dd, "(1{,}y{,}1)"]
	\\
	\\
	(SX{,}Y{,}RZ)
	\arrow[rr, "(1{,}1{,}z)"]
	\arrow[dd, "(x{,}1{,}1)"']
	&{}\arrow[dd, Rightarrow, shorten = 15,"(x{,}1{,}z)"]
	&(SX{,}Y{,}Z)\arrow[dd, "(x{,}1{,}1)"]
	&=&{}
	\\
	\\
	(X{,}Y{,}RZ)
	\arrow[rr, "(1{,}1{,}z)"']
	&{}&(X{,}Y{,}Z)
\end{tikzcd}\begin{tikzcd}[font=\fontsize{9}{6}]
(SX{,}TY{,}RZ)
\arrow[rr, "(1{,}1{,}z)"]
\arrow[dddd, "(x{,}1{,}1)\circ(1{,}y{,}1)"description]
&{}\arrow[dddd, Rightarrow, shorten = 20,"\left((x{,}1)\circ(1{,}y){,}z\right)"description]
&(X{,}TY{,}Z)\arrow[dddd, "(x{,}1{,}1)\circ(1{,}y{,}1)"description]
\\
\\
\\
\\
(X{,}Y{,}RZ)
\arrow[rr, "(1{,}1{,}z)"']
&{}&(X{,}Y{,}Z)
\end{tikzcd}$$

\noindent Similarly, both $\mathbf{Psmnd}(\Gr)_{\text{ast}}$ and $\Gr_{\text{2-nat}}^\mathbf{2}$ satisfy the symmetry condition on the nose rather than up to invertible modifications $\Sigma$. As such, the condition in Definition \ref{definition sylleptic opmonoidal functor} follows by cancelling $(x, y)$ with its inverse. This corresponds to the fact that the braiding for the $\mathbf{Gray}$ tensor product strictly satisfies the symmetry condition.

\end{proof}

\subsection{On Lifting Braidings, Syllapses and Symmetries}\label{Subsection on lifting braidings, syllapses and symmetries}

\begin{remark}\label{remark data of lifted braid and syllepsis}
	Let $(\A, \bigoplus, I, \beta, r, l)$ be a braided $\mathbf{Gray}$-monoid equipped with a braided semi-strict opmonoidal pseudomonad $(S, \eta, \mu)$. The component pseudonatural equivalences and invertible modifications of the braided monoidal bicategory structure on $\A^S$ can be explicitly described. The braiding is given by the pseudonatural equivalence displayed below.
	
	$$\begin{tikzcd}[font=\fontsize{9}{6}]
		\mathcal{A}^{S}\otimes \mathcal{A}^{S}
		\arrow[dd, "M_{S{,}S}"'] \arrow[rr, "\tau"] &{}\arrow[dd, Leftarrow, shorten = 10, "\Xi_{S{,}S}"]
		&\mathcal{A}^{S}\otimes \mathcal{A}^{S}
		\arrow[dd, "M_{S{,}S}"]
		\\
		\\
		{\left(\mathcal{A}\otimes\mathcal{A}\right)}^{S\otimes S} \arrow[rd, "\overline{\bigoplus}"']
		\arrow[rr, "\overline{\tau}"]
		&{}\arrow[d,Leftarrow, "\overline{\left(\beta{,}\Omega\right)}"]
		& {\left(\mathcal{A}\otimes\mathcal{A}\right)}^{S\otimes S} \arrow[ld, "\overline{\bigoplus}"]
		\\
		&\mathcal{A}^{S}
	\end{tikzcd}$$

\noindent The remaining data can also be described via whiskering the images of the data in $\PsG$ described in Proposition \ref{Proposition Braidings Syllapses and Symmetries as data in Psmnd Gray}, with certain morphisms in $\mathbf{Gray}(\A^S\otimes \A^S\otimes \A^S, \A^S)$. These morphisms are $2$-natural transformations given as the pseudonaturality constraints of $M_{-, ?}$ or the components of $\Xi_{-, ?}$. The same is true for the syllepsis or symmetry if $(\A, \bigoplus, I, \beta)$ has a syllepsis or symmetry $\Sigma$, and $(S, \chi, \iota, \gamma, \delta, \omega, \Omega)$ is compatible with $\Sigma$. For example, the modification that will be the syllepsis is given by the $2$-cell in $\mathbf{Gray}(\A^S\otimes \A^S, \A^S)$ from the morphism displayed below left to the one displayed below right, given in terms of $\overline{\Sigma}$. Note that the pasting in $\mathbf{Gray}$ displayed below left indeed determines a unique pseudonatural transformation by Corollary \ref{EM of ambistrict is 2-natural} part (1), since the components of $\Xi_{S, S}$ are tight pseudomorphisms. Also, note that the $2$-natural transformation displayed below right is indeed the identity, as can be seen by inspecting the components of $\Xi_{S, S}$.

$$\begin{tikzcd}[font=\fontsize{8}{5}, column sep = 12]
	\mathcal{A}^{S}\otimes \mathcal{A}^{S}\arrow[rrrr, bend left = 30, "1"]\arrow[rrrr, equal, shorten = 60, shift left = 8]
	\arrow[dd, "M_{S{,}S}"'] \arrow[rr, "\tau"] &{}\arrow[dd, Leftarrow, shorten = 10, "\Xi_{S{,}S}"]
	&\mathcal{A}^{S}\otimes \mathcal{A}^{S}
	\arrow[dd, "M_{S{,}S}"] \arrow[rr, "\tau"]
	&{}\arrow[dd, Leftarrow, shorten = 10, "\Xi_{S{,}S}"]
	&\mathcal{A}^{S}\otimes \mathcal{A}^{S}
	\arrow[dd, "M_{S{,}S}"]
	\\
	\\
	{\left(\mathcal{A}\otimes\mathcal{A}\right)}^{S\otimes S} \arrow[rrd,bend right = 20, "\overline{\bigoplus}"']
	\arrow[rr, "\overline{\tau}"]
	&{}\arrow[d,Leftarrow,shorten = 10, "\overline{\left(\beta{,}\Omega\right)}"]
	& {\left(\mathcal{A}\otimes\mathcal{A}\right)}^{S\otimes S} \arrow[d, "\overline{\bigoplus}"]\arrow[rr, "\overline{\tau}"]
	&{}\arrow[d,Leftarrow,shorten = 10, "\overline{\left(\beta{,}\Omega\right)}"']
	& {\left(\mathcal{A}\otimes\mathcal{A}\right)}^{S\otimes S} \arrow[dll,bend left = 20, "\overline{\bigoplus}"]
	\\
	&{}&\mathcal{A}^{S}&{}&&\Rrightarrow^{\overline{\Sigma}}
\end{tikzcd}\begin{tikzcd}[font=\fontsize{9}{6}, column sep = 12]
\mathcal{A}^{S}\otimes \mathcal{A}^{S}\arrow[rrrr, bend left = 30, "1"]\arrow[rrrr, equal, shorten = 60, shift left = 8]
\arrow[dd, "M_{S{,}S}"'] \arrow[rr, "\tau"] &{}\arrow[dd, Leftarrow, shorten = 10, "\Xi_{S{,}S}"]
&\mathcal{A}^{S}\otimes \mathcal{A}^{S}
\arrow[dd, "M_{S{,}S}"] \arrow[rr, "\tau"]
&{}\arrow[dd, Leftarrow, shorten = 10, "\Xi_{S{,}S}"]
&\mathcal{A}^{S}\otimes \mathcal{A}^{S}
\arrow[dd, "M_{S{,}S}"]
\\
\\
{\left(\mathcal{A}\otimes\mathcal{A}\right)}^{S\otimes S}
\arrow[rrrr, equal, shorten = 60, red,shift right = 8] \arrow[rrd,bend right = 20, "\overline{\bigoplus}"']
\arrow[rr, "\overline{\tau}"]
&{}
& {\left(\mathcal{A}\otimes\mathcal{A}\right)}^{S\otimes S}
\arrow[rr, "\overline{\tau}"] 
&{}
& {\left(\mathcal{A}\otimes\mathcal{A}\right)}^{S\otimes S} \arrow[dll,bend left = 20, "\overline{\bigoplus}"]
\\
&{}&\mathcal{A}^{S}&{}
\end{tikzcd}$$

\end{remark}

\begin{theorem}\label{Theorem lifted braidings syllapses and symmetries}
	In the context of Remark \ref{remark data of lifted braid and syllepsis}, the following statements hold.
	\begin{enumerate}
		\item For a braided semi-strict opmonoidal pseudomonad, the data described in Remark \ref{remark data of lifted braid and syllepsis} comprise a braided monoidal bicategory structure, in the sense of Section 2.4 of \cite{Gurski Loop Spaces}, on $(\A^S, \overline{\bigoplus}, \overline{I}, \overline{\gamma}, \overline{\omega}, \overline{\delta})$.
		
		\item If $\Sigma$ is a syllepsis for $(\A, \bigoplus ,I)$ and $S: \A \to \A$ is sylleptic then the braided $\mathbf{Gray}$-pseudomonoid structure of part (1) is moreover sylleptic, with syllepsis given as described in Remark \ref{remark data of lifted braid and syllepsis}.
		
		\item If $\Sigma$ is moreover a symmetry then the syllepsis on $\A^S$ is also a symmetry.
	\end{enumerate} 
\end{theorem}

\begin{proof}
	The axioms for braiding, syllepsis and symmetry are all equations between $2$-cells in $\A^S$, given either by
	
	\begin{itemize}
		\item $2$-cell components of $\overline{(\beta, \Omega)}$,
		\item the components of $\overline{r}$, $\overline{l}$ or $\overline{\Sigma}$,
		\item interchangers in $\A^S$ between pseudomorphism components of the associator and the braiding.
		\item interchangers in $\A^S$ between two different pseudomorphism components of $\overline{(\beta, \Omega)}$.
	\end{itemize}

\noindent These $2$-cells are all the analogous pieces of data in $\A$. But the pseudomonadic $2$-functor $U^S: \A^S \to \A$ is faithful on $2$-cells, and as such the axioms hence follow immediately from the analogous axioms for these data in $\A$.
\end{proof}

\begin{corollary}\label{Corollary forgetful is braided/sylleptic}
	\hspace{1mm}
	\begin{enumerate}
		\item In the context of Theorem \ref{Theorem lifted braidings syllapses and symmetries} part (1), the forgetful $2$-functor $\A^S \to \A$ is strictly braided
		\item In the context of Theorem \ref{Theorem lifted braidings syllapses and symmetries} parts (2) and (3), the forgetful $2$-functor $\A^S \to \A$ it is strictly sylleptic.
	\end{enumerate}
\end{corollary}

\begin{proof}
	This is easy to inspect from the description of lifted braid and syllepsis structures given in Remark \ref{remark data of lifted braid and syllepsis}, since the components of the $2$-natural isomorphism $\Xi_{S, S}$ are tight pseudomorphisms.
\end{proof}

\section{Conclusions and future directions}

\noindent The proof strategies pursued in this paper for lifting structures are motivated by the following hypotheses.

\begin{enumerate}
	\item Compatibility structures between pseudomonads and symmetric $\mathbf{Gray}$-monoids should correspond to appropriately symmetric three-dimensional monoids in a three dimensional monoidal structure on the $\mathbf{Gray}$-category $\PsG$.
	\item The $\mathbf{Gray}$-functor $\mathbf{EM}: \PsG \to \Gr^\mathbf{2}$ should be symmetric lax monoidal.
	\item $\mathbf{EM}$ should hence carry symmetric three-dimensional monoids to symmetric monoidal bicategory structure on $2$-categories of pseudoalgebras.
\end{enumerate} 

\noindent Despite the fact that there is as yet no precise definition for symmetric monoidal $\mathbf{Gray}$-categories or lax functors between them, we have established enough aspects of the above hypotheses to lift symmetric $\mathbf{Gray}$-monoid structures on the base $\A$ to corresponding structures on $\A^S$.
\\
\\
\noindent We restricted the source and target to establish two-dimensional versions of the monoidal fragments of hypotheses (1) (Theorem \ref{theorem semi-strict opmonoidal pseudomonads are pseudomonoids in the Gray monoid of pseudomonads}) and (2) (Theorem \ref{Theorem EM is normal lax Gray functor}), however a direct verification of the axioms for the lifted monoidal structure was needed in the proof of Theorem \ref{monoidal bicategory structure on pseudoalgebras}. For braidings, syllapses and symmetries, compatibility with pseudomonads still corresponded to cells in the $\mathbf{Gray}$-category $\PsG$. Simple aspects of compatibility between $\mathbf{EM}$ and the symmetric fragments of the $\mathbf{Gray}$ tensor product of pseudomonads were established in Theorem \ref{theorem EM is symmetric monoidal} and used to leverage the results of \cite{Formal Theory of Pseudomonads} to lift structures on the base to pseudonatural equivalences and invertible modifications at the level of pseudoalgebras. The axioms needed for these lifted data to form a well-defined symmetric monoidal bicategory required separate verification, but followed easily using faithfulness on $2$-cells of the forgetful $2$-functor $\A^S \to \A$.
\\
\\
\noindent In the one dimensional setting, the $2$-categories $\mathbf{MonCat}_\text{oplax}$, $\mathbf{BrMonCat}_\text{oplax}$ and $\mathbf{SymMonCat}_\text{oplax}$ all have Eilenberg-Moore objects for monads, which are created by the forgetful $2$-functors to $\Cat$. This follows from Proposition 4.1.1 of \cite{Lack Limits and Lax Morphisms}, since all of these are $2$-categories of strict algebras and colax morphisms for appropriate $2$-monads on $\Cat$. See also Propositions 2.9 and 2.16 of \cite{McCrudden Opmonoidal Monads}. Moreover, these results also restrict to the setting of strict monoidal categories, and braided or symmetric strict monoidal categories. Some weakness will be needed in the analogues of such results in the two-dimensional setting, since as we have observed in Theorem \ref{monoidal bicategory structure on pseudoalgebras}, the lifted monoidal structure on $2$-categories of pseudoalgebras is weaker than a $\Gr$-monoid. See Chapter 6 of \cite{Miranda PhD} and the forthcoming \cite{Miranda Tricategorical Limits and Colimits}. The strictification of the monoidal (resp. braided, sylleptic, symmetric) bicategory $\A^S$ should have a universal property expressible in these terms.
\\
\\
\noindent Finally, we mention that many of the constructions and proofs in this paper are formal, in the sense that they hold for $\mathcal{V}$-enriched monoidal categories $\mathfrak{K}$ where $\mathcal{V}= \mathbf{Gray}$, provided that $\mathfrak{K}$ has Eilenberg-Moore objects. We have chosen to present our arguments for the case $\mathfrak{K} = \mathbf{Gray}$. This is because many of the constructions are simple enough to understand `pointwise', but become much harder to parse when expressed at a higher level of abstraction, despite still remaining valid there. Formal aspects will be described in detail in the forthcoming \cite{Miranda Extending Symmetric Monoidal Structures to Kleisli bicategories}, where an investigation of extensions of symmetric monoidal structures to Kleisli bicategories will also be undertaken.

\section{Appendix}\label{Section Appendix}

\subsection{Theorem \ref{theorem strictifying strong monoidal pseudomonads}}\label{Appendix Coherence}

\begin{proof}
	Recall from \cite{Low dimensional structures formed by tricategories} that every tricategory is biequivalent internally to the tricategory $\mathbf{Tricat}_\text{oplax}$ \cite{Gurski biequivalence in tricategories} to a $\mathbf{Gray}$-category. Indeed the $1$-cells in this biequivalence are genuine trihomomorphisms rather than just oplax ones. In particular, the monoidal bicategory is biequivalent in $\mathbf{MonBicat}_\text{oplax}$ to a $\mathbf{Gray}$-monoid via a strong monoidal biequivalence $E_{*}: \mathcal{V} \to \overline{\mathcal{V}}$, whose biequivalence inverse we will write as $E^{*}$. Recall also that $E^{*}E_{*}=1_{\mathcal{V}}$, and the pseudo-icon equivalence $\phi: 1_{\overline{\mathcal{V}}}\rightarrow E_{*}E^{*}$ satisfies $\phi.E_{*}= 1_{E_{*}}$. We write $\Phi$ for the costrict trimodification witnessing the other triangle identity for the biadjoint biequivalence $E_{*}\dashv E^{*}$.
	\\
	\\
	\noindent The opmonoidal pseudomonad structure transports along this biequivalence in $\mathbf{MonBicat}_\text{oplax}$ resulting in a pseudomonad on the $\mathbf{Gray}$-monoid $\left(\overline{\mathcal{V}}, \overline{\bigoplus}, \overline{I}\right)$. We give an explicit description of the transported pseudomonad. It has
	
	\begin{itemize}
		\item underlying opmonoidal endofunctor $E^{*}TE_{*}$,
		\item unit \begin{tikzcd}[font=\fontsize{9}{6}]
			1_{\overline{\mathcal{V}}}\arrow[rr, "\phi"] 
			&& E_{*}E^{*} \arrow[rr, "1.\eta .1"] 
			&&  E_{*}TE^{*}
		\end{tikzcd}
		\item multiplication \begin{tikzcd}[font=\fontsize{9}{6}]
			E_{*}TE^{*}E_{*}TE^{*}\arrow[rr, equal] && E_{*}T^{2}E^{*} \arrow[rr, "1.\mu. 1"] &&  E_{*}TE^{*}
		\end{tikzcd}
		\item Right unitor as given below.
		
		$$\begin{tikzcd}[column sep = 12, row sep = 10, font=\fontsize{9}{6}]
			E_{*}TE^{*}\arrow[rrrrdddd, bend right = 30, "1"']
			\arrow[rr, "1.\phi"] 
			&& E_{*}TE^{*}E_{*}E^{*} \arrow[rr, "1.\eta.1"] \arrow[rrdddd, "1"]\arrow[dddd, Rightarrow, shorten = 20, "\Phi"]
			&{}\arrow[dd, Rightarrow,shorten = 8, "\mathbf{r}"]&  E_{*}TE^{*}E_{*}TE^{*}\arrow[dd, equal]
			\\
			\\
			&&&{}&E_{*}T^{2}E^{*}\arrow[dd, "1.\mu.1"]
			\\
			\\
			&&{}&&E_{*}TE^{*}
		\end{tikzcd}$$
		\item Left unitor and associator given by $\mathbf{l}$ and $\mathbf{a}$.
	\end{itemize} 
	
	\noindent The associativity law for the new pseudomonad follows immediately from the analogous law for the original pseudomonad, while the unit law uses the analogous law for the original pseudomonad as well as the axioms for the biadjoint biequivalence $E_{*} \dashv E^{*}$, listed in Figures 1 and 2 of Definition 2.1 of \cite{Gurski biequivalence in tricategories}. Moreover, $E_{*} \dashv E^{*}$ lifts to a biequivalence internal to the tricategory $\mathbf{Psmnd}\left(\mathbf{MonBicat}_\text{oplax}\right)$. The morphisms of pseudomonads are strict in this biequivalence, while the unit is given in terms of $\phi$. The remaining details in the proof of part (1) are easy to verify and we leave this verification to the interested reader.
	\\
	\\
	\noindent For part (2), suppose the original pseudonatural transformations $\iota$ and $\chi$ are part of adjoint equivalences. The original pseudomonad $\mathbb{T}: \mathbf{Psmnd} \rightarrow \mathbf{MonBicat}_\text{oplax}$ factors through the wide sub-tricategory consisting of the strong monoidal pseudofunctors. Recall from part (4) of Proposition 4.3.4 of \cite{Miranda strictifying operational coherences} that the inclusion $I_{3}: \mathbf{Gray}$-$\mathbf{Cat}_{3} \rightarrow \mathbf{Tricat}_{3}$ has a left triadjoint $\mathbf{st}_{3}$. It is easy to see that this triadjunction restricts to one between $\mathbf{MonBicat}_\text{strong}$ and $\mathbf{Gray}$-$\mathbf{monoids}$. But since the $\mathbf{Gray}$-category $\mathbf{Psmnd}$ is cofibrant and $\mathbf{Gray}$-$\mathbf{monoids}$ is a $\mathbf{Gray}$-category, the trihomomorphism depicted below is equivalent to a $\mathbf{Gray}$-functor.
	
	$$
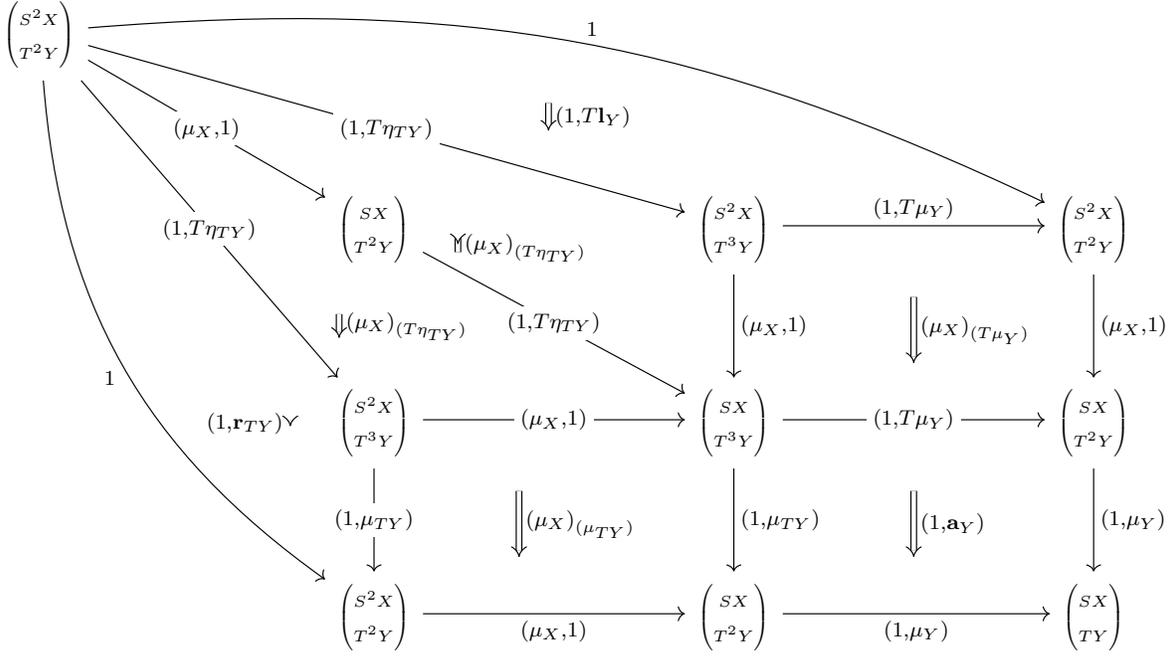
\begin{tikzcd}
		\mathbf{Psmnd} \arrow[rr, "\mathbb{T}"] && \mathbf{MonBicat}_\text{strong} \arrow[rr, "\mathbf{st}_{3}"] && \mathbf{Gray}\text{-}\mathbf{monoids}
	\end{tikzcd}$$
	
	\noindent This gives the required semi-strict monoidal pseudomonad $(\mathbf{st}_{3}(\mathcal{V}),\mathbf{st}_{3}(T))$ whose underlying endomorphism is strict and whose unit and multiplication are semi-strict. The biequivalence between $(\mathcal{V}, T)$ and $(\mathbf{st}_{3}(\mathcal{V}),\mathbf{st}_{3}(T))$ is given in terms of the data of the triadjunction $\mathbf{st}_{3}\dashv I_{3}$.
\end{proof} 

\subsection{Lemma \ref{Gray tensor of two pseudomonads}}\label{Appendix proof of Gray tensor product of two pseudomonads}

\begin{proof}
	To conserve space, we denote an element of the set $\A_{0}\times \B_{0}$ as a column rather. For the unit law, we need to show that the pasting diagram displayed below is the identity.
	
	$$\begin{tikzcd}[font=\fontsize{6}{3}, column sep = 35]
		\begin{pmatrix}S^2X
			\\
			T^2Y\end{pmatrix}
		\arrow[rrrrdd, bend left = 20, "1"]
		\arrow[rrdddd, bend right = 20, "1"']
		\arrow[rrdd, "(S\eta_{SX}{,}1)"description]
		&&{}\arrow[dddd, Rightarrow, shorten = 55, shift right = 15, "(\mathbf{r}_{SX}{,}1)"']
		\arrow[dd, Rightarrow, shorten = 20, "(S\mathbf{l}_{X}{,}1)"]
		\\
		\\
		&&\begin{pmatrix}S^3X
			\\
			T^2Y\end{pmatrix}
		\arrow[rr, "(S\mu_{X}{,}1)"]
		\arrow[dd, "(\mu_{SX}{,}1)"description]
		\arrow[rrdd, "(1{,}T\eta_{TY})"description]
		&&\begin{pmatrix}S^2X
			\\
			T^2Y\end{pmatrix}
		\arrow[dddd,shift right = 10, Rightarrow, shorten = 50, "{(\mu_{SX})}_{(T\eta_{TY})}"']
		\arrow[dd, Rightarrow, shorten = 15, shift right = 5, "{(S\mu_{X})}_{(T\eta_{TY})}"]
		\arrow[rrdd, "(1{,}T\eta_{TY})"]
		\arrow[rrrrdd, bend left = 20, "1"]
		&&{}\arrow[dd, Rightarrow, shorten = 20, "(1{,}T\mathbf{l}_{Y})"]		
		\\
		\\
		&&\begin{pmatrix}S^2X
			\\
			T^2Y\end{pmatrix}
		\arrow[rrdd, "(1{,}T\eta_{TY})"description]
		\arrow[rrdddd, bend right = 20, "1"']
		&&\begin{pmatrix}S^3X
			\\
			T^3Y\end{pmatrix}
		\arrow[dddd, Rightarrow,shift right = 12, shorten = 55, "(1{,}\mathbf{r}_{TY})"']
		\arrow[rr, "(S\mu_{X}
		{,}
		1)"]
		\arrow[dd, "(\mu_{SX}
		{,}
		1)"description]
		&{}\arrow[dd, Rightarrow,shift right = 5,  shorten = 20, "(\mathbf{a}_X
		{,}
		1)"]
		&\begin{pmatrix}S^2X
			\\
			T^3Y\end{pmatrix}
		\arrow[rr, "(1
		{,}
		T\mu_Y)"]
		\arrow[dd, "(\mu_X
		{,}
		1)"]
		&{}\arrow[dd, Rightarrow, shorten = 20, "{(\mu_{X})}_{(T\mu_{Y})}"]
		&\begin{pmatrix}S^2X
			\\
			T^2Y\end{pmatrix}
		\arrow[dd, "(\mu_X
		{,}
		1)"]
		\\
		\\
		&&&&\begin{pmatrix}S^2X
			\\
			T^3Y\end{pmatrix}
		\arrow[dd, "(1
		{,}
		\mu_{TY})"description]
		\arrow[rr, "(\mu_{X}{,}1)"description]
		&{}\arrow[dd, Rightarrow,shift right = 5,  shorten = 20, "{(\mu_{X})}_{(\mu_{TY})}"]
		&\begin{pmatrix}SX
			\\
			T^3Y\end{pmatrix}
		\arrow[dd, "(1
		{,}
		\mu_{TY})"]
		\arrow[rr, "(1{,}
		T\mu_{Y})"description]
		&{}\arrow[dd, Rightarrow, shorten = 20, "(1
		{,}
		\mathbf{a}_{Y})"]
		&\begin{pmatrix}SX
			\\
			T^2Y\end{pmatrix}
		\arrow[dd, "(1
		{,}
		\mu_{Y})"]
		\\
		\\
		&&&&\begin{pmatrix}S^2X
			\\
			T^2Y\end{pmatrix}
		\arrow[rr, "(\mu_{X}
		{,}
		1)"']
		&{}&\begin{pmatrix}SX
			\\
			T^2Y\end{pmatrix}
		\arrow[rr, "(1{,}
		\mu_{Y})"']
		&{}&
		\begin{pmatrix}SX
			\\
			TY\end{pmatrix}
	\end{tikzcd}$$
	
	\noindent Apply the interchange law for the morphism $T\eta_{TY}: T^2Y \to T^3 Y$ on the $2$-cell $\mathbf{a}_{X}$, followed by the unit law for the pseudomonad $(\A, S)$ to simplify this to the pasting diagram displayed below.

	$$\begin{tikzcd}[font=\fontsize{6}{3}, column sep = 45]
		\begin{pmatrix}S^2X
			\\
			T^2Y\end{pmatrix}
		\arrow[dddd, Rightarrow, shift left =60, shorten = 65, "{(\mu_{X})}_{(T\eta_{TY})}"]
		\arrow[rrdd,"(\mu_{X}{,}1)"description]
		\arrow[rrrrdd, "(1{,}T\eta_{TY})"description]
		\arrow[rrrrrrdd, bend left = 15, "1"]
		\arrow[rrdddd, "(1{,}T\eta_{TY})"description]
		\arrow[rrdddddd, bend right = 25, "1"']
		&&{}\arrow[dd, Rightarrow, shorten = 20,shift left = 25, "(1{,}T\mathbf{l}_{Y})"]		
		\\
		\\
		&&\begin{pmatrix}SX
			\\
			T^2Y\end{pmatrix}
		\arrow[dd, shorten = 15, shift right = 5, Rightarrow,"{(\mu_{X})}_{(T\eta_{TY})}"]
		\arrow[dddd, Rightarrow,shift right = 12, shorten = 55, "(1{,}\mathbf{r}_{TY})"']
		\arrow[rrdd,"(1{,}T\eta_{TY})"description]
		&{}
		&\begin{pmatrix}S^2X
			\\
			T^3Y\end{pmatrix}
		\arrow[rr, "(1
		{,}
		T\mu_Y)"]
		\arrow[dd, "(\mu_X
		{,}
		1)"]
		&{}\arrow[dd, Rightarrow, shorten = 20, "{(\mu_{X})}_{(T\mu_{Y})}"]
		&\begin{pmatrix}S^2X
			\\
			T^2Y\end{pmatrix}
		\arrow[dd, "(\mu_X
		{,}
		1)"]
		\\
		\\
		{}&&\begin{pmatrix}S^2X
			\\
			T^3Y\end{pmatrix}
		\arrow[dd, "(1
		{,}
		\mu_{TY})"description]
		\arrow[rr, "(\mu_{X}{,}1)"description]
		&{}\arrow[dd, Rightarrow,shift right = 5,  shorten = 20, "{(\mu_{X})}_{(\mu_{TY})}"]
		&\begin{pmatrix}SX
			\\
			T^3Y\end{pmatrix}
		\arrow[dd, "(1
		{,}
		\mu_{TY})"]
		\arrow[rr, "(1{,}
		T\mu_{Y})"description]
		&{}\arrow[dd, Rightarrow, shorten = 20, "(1
		{,}
		\mathbf{a}_{Y})"]
		&\begin{pmatrix}SX
			\\
			T^2Y\end{pmatrix}
		\arrow[dd, "(1
		{,}
		\mu_{Y})"]
		\\
		\\
		&&\begin{pmatrix}S^2X
			\\
			T^2Y\end{pmatrix}
		\arrow[rr, "(\mu_{X}
		{,}
		1)"']
		&{}&\begin{pmatrix}SX
			\\
			T^2Y\end{pmatrix}
		\arrow[rr, "(1{,}
		\mu_{Y})"']
		&{}&
		\begin{pmatrix}SX
			\\
			TY\end{pmatrix}
	\end{tikzcd}$$

	\noindent Next, apply the interchange law for $\mu_{X}: S^2X \to SX$ on the $2$-cells $T\mathbf{l}_{Y}$ and $\mathbf{r}_{TY}$. Finally, observe that the unit law for the pseudomonad $(\B, T)$ says that the resulting pasting is the identity on ${(\mu\otimes \mu)}_{X, Y}$, completing the proof.

	
\end{proof}

\subsection{Proposition \ref{Proposition 2-cell component of compositor as the second pseudomonad varies} part 1}\label{Appendix M 1 h well-defined}

\begin{proof}
	The $2$-natural transformation $M_{1, H}$ is induced by the universal property of ${\left(\mathcal{A}\otimes \mathcal{B}\right)}^{\mathcal{A}\otimes H}$ given by a pseudomorphism for $S\otimes T'\circ -$. This pseudomorphism has $1$-cell component given by the identity pseudonatural transformation, and has $2$-cell component given by the modification formed by whiskering the interchanger ${\left(\varepsilon^{S}\right)}_{h}$ with the pseudonatural transformation $\varepsilon^{T}$. This modification indeed has component $2$-cell given as specified in part (1). We explain how to verify that this pseudomorphism is well-defined, via calculations involving data in the $2$-category $\A \otimes \B'$. 
	\\
	\\
	\noindent Consider the pasting diagram depicted below and observe that the unit law follows by applying the interchange law for the morphism $x: SX \to X$ in $\A$ on the $2$-cell $h_{0, Y}: h_{Y}.\eta_{HY}' \Rightarrow H\eta_{Y}$ in $\mathcal{B}$.

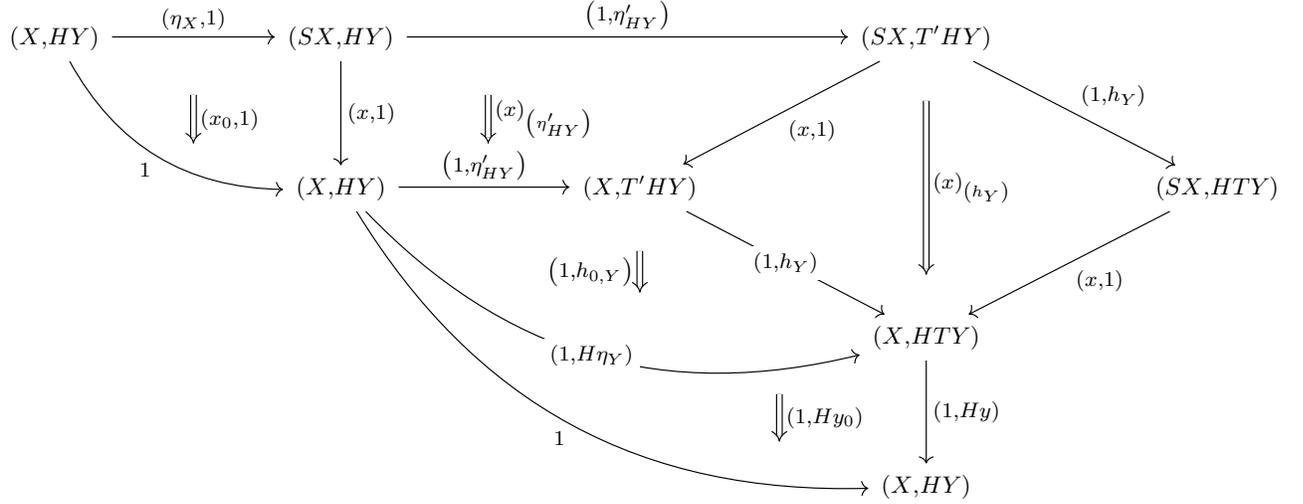
	$$\begin{tikzcd}[font=\fontsize{9}{6}]
		\left(X{,}HY\right)\arrow[rr, "\left(\eta_{X}{,}1\right)"]
		\arrow[rrdd,bend right=30, "1"']
		&{}\arrow[dd, Rightarrow, shorten = 15, "\left(x_{0}{,}1\right)"]
		&\left(SX{,}HY\right)
		\arrow[dd, "\left(x{,}1\right)"]
		\arrow[rrrr, "\left(1{,}\eta_{HY}'\right)"]
		&{}\arrow[dd, Rightarrow, shorten = 15, "(x)_{\left(\eta_{HY}'\right)}"]
		&&&\left(SX{,}T'HY\right)
		\arrow[lldd, "\left(x{,}1\right)"]
		\arrow[dddd, Rightarrow, shorten = 15, "{(x)}_{\left(h_{Y}\right)}"]
		\arrow[rrdd, "\left(1{,}h_{Y}\right)"]
		\\
		\\
		&{}&\left(X{,}HY\right)
		\arrow[rr, "\left(1{,}\eta_{HY}'\right)"]
		\arrow[rrrrdd, bend right = 30, "\left(1{,}H\eta_{Y}\right)"description]
		\arrow[rrrrdddd, bend right = 30, "1"']
		&{}&\left(X{,}T'HY\right)
		\arrow[dd, Rightarrow, shorten = 15, "\left(1{,}h_{0{,}Y}\right)"']
		\arrow[rrdd, "\left(1{,}h_{Y}\right)"description]
		&&&& \left(SX{,}HTY\right)
		\arrow[lldd, "\left(x{,}1\right)"]
		\\
		\\
		&&&&{}\arrow[dd, Rightarrow, shorten = 15, shift left = 20, "\left(1{,}Hy_{0}\right)"]
		&&\left(X{,}HTY\right)\arrow[dd, "\left(1{,}Hy\right)"]
		\\
		\\
		&&&&{}&&\left(X{,}HY\right)
	\end{tikzcd}$$
	
	\noindent The proof of the multiplication law is similar. Start with the pasting diagram displayed below, in which elements of $\A_{0} \times \B_{0}$ are denoted as columns rather than rows to conserve space. 
	
	$$
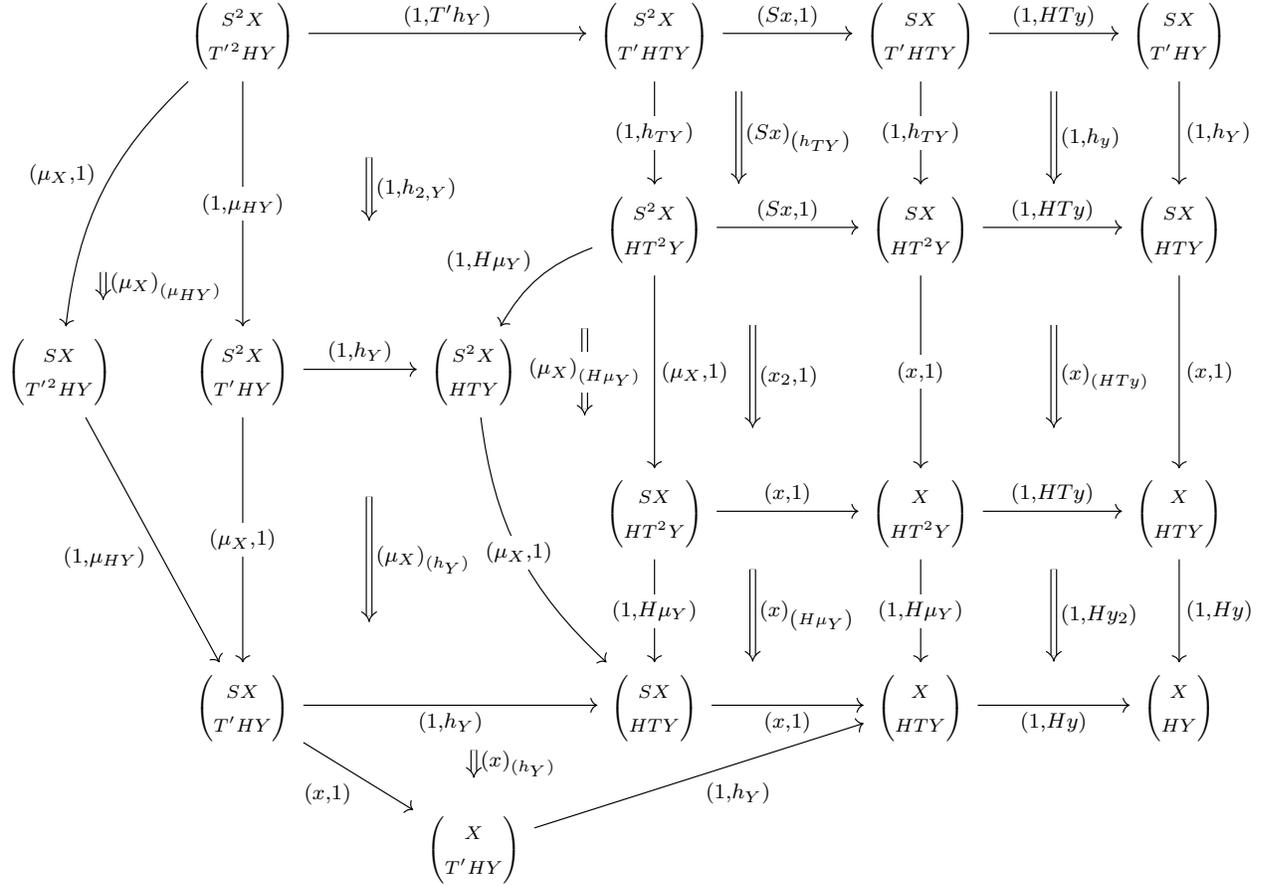
\begin{tikzcd}[font=\fontsize{7}{4}, column sep = 20]
		&\begin{pmatrix}
			S^2X\\
			{T'}^2HY
		\end{pmatrix}
		\arrow[ddd, "(1{,}\mu_{HY})"description]
		\arrow[rrr, "\left(1{,}T'h_{Y}\right)"]
		\arrow[dddl, bend right = 20,"\left(\mu_{X}{,}1\right)"']
		&&{}\arrow[ddd, Rightarrow, shorten = 40, shift right = 15, "(1{,}h_{2{,}Y})"]
		&\begin{pmatrix}
			S^2X\\
			{T'}HTY
		\end{pmatrix}
		\arrow[rr, "\left(Sx{,}1\right)"]
		\arrow[dd, "\left(1{,}h_{TY}\right)"description]
		&{}\arrow[dd, Rightarrow,shift right = 7,  shorten = 15, "{(Sx)}_{\left(h_{TY}\right)}"]
		&\begin{pmatrix}
			SX\\
			{T'}HTY
		\end{pmatrix}
		\arrow[rr, "\left(1{,}HTy\right)"]
		\arrow[dd, "\left(1{,}h_{TY}\right)"description]
		&{}\arrow[dd, Rightarrow, shorten = 15, "\left(1{,}h_{y}\right)"]
		&\begin{pmatrix}
			SX\\
			{T'}HY
		\end{pmatrix}
		\arrow[dd, "\left(1{,}h_{Y}\right)"]
		\\
		\\
		&{}
		\arrow[d, Rightarrow,shift right = 20, shorten = 10, "{(\mu_{X})}_{(\mu_{HY})}"]
		&&&\begin{pmatrix}
			S^2X\\
			H{T}^2Y
		\end{pmatrix}
		\arrow[ld, bend right = 20, "(1{,}H\mu_{Y})"']
		\arrow[dd, Rightarrow, shorten = 20, shift right = 10,"{(\mu_{X})}_{(H\mu_{Y})}"description]
		\arrow[rr, "\left(Sx{,}1\right)"]
		\arrow[dd, "\left(\mu_{X}{,}1\right)"]
		&{}\arrow[dd, Rightarrow,shift right = 5,  shorten = 30, "\left(x_{2}{,}1\right)"]
		&\begin{pmatrix}
			SX\\
			H{T}^2Y
		\end{pmatrix}
		\arrow[rr, "\left(1{,}HTy\right)"]
		\arrow[dd, "\left(x{,}1\right)"description]
		&{}\arrow[dd, Rightarrow, shorten = 30, "{(x)}_{\left(HTy\right)}"]
		&\begin{pmatrix}
			SX\\
			H{T}Y
		\end{pmatrix}
		\arrow[dd, "\left(x{,}1\right)"]
		\\
		\begin{pmatrix}
			SX\\
			{T'}^2HY
		\end{pmatrix}
		\arrow[dddr, "\left(1{,}\mu_{HY}\right)"']
		&\begin{pmatrix}
			S^2X\\
			{T}'HY
		\end{pmatrix}
		\arrow[ddd, "(\mu_{X}{,}1)"description]
		\arrow[rr, "(1{,}h_{Y})"]
		&&\begin{pmatrix}
			S^2X\\
			HTY
		\end{pmatrix}
		\arrow[ddd, Rightarrow, shorten = 30, shift right = 15,"{(\mu_{X})}_{(h_{Y})}"]
		\arrow[rddd, bend right = 20,"(\mu_{X}{,}1)"description]
		\\
		&{}&&&\begin{pmatrix}
			SX\\
			HT^2Y
		\end{pmatrix}
		\arrow[dd, "\left(1{,}H\mu_{Y}\right)"description]\arrow[rr, "\left(x{,}1\right)"]
		&{}\arrow[dd, Rightarrow,shift right =5, shorten = 15, "(x)_{\left(H\mu_{Y}\right)}"]
		&\begin{pmatrix}
			X\\
			HT^2Y
		\end{pmatrix}
		\arrow[dd, "\left(1{,}H\mu_{Y}\right)" description]
		\arrow[rr, "\left(1{,}HTy\right)"]
		&{}\arrow[dd, Rightarrow, shorten = 15, "\left(1{,}Hy_{2}\right)"]&
		\begin{pmatrix}
			X\\
			HTY
		\end{pmatrix}
		\arrow[dd, "\left(1{,}Hy\right)"]
		\\
		\\
		&\begin{pmatrix}
			SX\\
			T'HY
		\end{pmatrix}
		\arrow[rrd, "(x{,}1)"']
		\arrow[rrr, "\left(1{,}h_{Y}\right)"']
		&&{}\arrow[d, Rightarrow, shorten =10, "(x)_{(h_{Y})}"]
		&\begin{pmatrix}
			SX\\
			HTY
		\end{pmatrix}
		\arrow[rr, "\left(x{,}1\right)"']
		&{}&\begin{pmatrix}
			X\\
			HTY
		\end{pmatrix}\arrow[rr, "\left(1{,}Hy\right)"']
		&{}&\begin{pmatrix}
			X\\
			HY
		\end{pmatrix}
		\\
		&&&\begin{pmatrix}
			X\\T'HY
		\end{pmatrix}\arrow[rrru, "(1{,}h_{Y})"']
	\end{tikzcd}$$
		
	\noindent The proof proceeds via various applications of the interchange axioms involving morphisms in the boundary of $h_{2, Y}$ on the $2$-cell $x_{2}$, and morphisms in the boundary of $x_{2}$ on the $2$-cells $h_{2, Y}$ and $Hy_{2}$. We leave details to the interested reader.
		\end{proof}

\subsection{Theorem \ref{monoidal bicategory structure on pseudoalgebras}}\label{Appendix monoidal bicategory structure on pseudoalgebras}

\begin{proof}
	By Proposition \ref{Proposition pseudomonoid structure on pseudoalgebras in Gray 2nat}, the usual pentagon and middle unit laws for a monoidal category hold up to the congruence generated by middle-four interchange. We need to check that these axioms are actually satisfied as equations between $2$-natural transformations. For the unit law, we need to show that the pasting diagram displayed below is the identity. But this follows from the interchange law for the morphism $x: SX \to X$ on the $2$-cell $\gamma_{Z}$, and the unit law for the opmonoidal $2$-functor $S: \A \to \A$.

	$$\begin{tikzcd}[font=\fontsize{9}{6}, column sep = 18]
		S(X \oplus I \oplus Z)
		\arrow[rr, "\chi_{X{,}I\oplus Z}"]
		\arrow[dd, "\chi_{X\oplus I{,}Z}"']
		&{}\arrow[dd, Leftarrow, shorten = 15, "\omega_{X{,}I{,}Z}"]
		&SX\oplus S(I\oplus Z)
		\arrow[rr, "x\oplus 1"]
		\arrow[dd, "1\oplus \chi_{I{,}Z}"description]
		&{}\arrow[dd, Leftarrow, shorten = 15, "x_{\chi_{I{,}Z}}"]
		&X\oplus S(I\oplus Z)\arrow[rrdddd, bend left = 30, "1\oplus 1_{SZ}"]
		\arrow[dd, "1\oplus \chi_{I{,}Z}"description]
		&{}\arrow[dddd, Rightarrow, shorten = 40, "1\oplus \gamma_{Z}"]
		\\
		\\
		S(X\oplus I)\oplus SZ
		\arrow[rrrrdd,bend right = 20, "1_{SX}\oplus 1"'] 
		\arrow[rr, "\chi_{X{,}I}\oplus 1"']
		&{}&SX\oplus SI\oplus SZ
		\arrow[dd, Leftarrow, shorten = 15, "\delta_{X}\oplus 1"] 
		\arrow[rrdd, "1\oplus \iota\oplus 1"']
		\arrow[rr, "x\oplus 1\oplus 1"'] 
		&{}& X \oplus SI \oplus SZ
		\arrow[rrdd, "1\oplus \iota\oplus 1"description] 
		\arrow[dd, Leftarrow, shorten = 15, "x_{\iota}\oplus 1"]
		\\
		\\
		&&{}&&SX\oplus I \oplus SZ\arrow[rr, "x\oplus 1 \oplus 1"']
		&{}& X \oplus I\oplus SZ
		\arrow[rr, "1 \oplus z"']
		&& X\oplus I\oplus Z
	\end{tikzcd}$$

	\noindent For the pentagon law, observe that the composite of tight pseudomorphisms $({\alpha}_{\mathbf{W}{,}\mathbf{X}{,}\mathbf{Y}}\oplus \mathbf{Z})\circ {\alpha}_{\mathbf{W}, \mathbf{X\oplus Y}, \mathbf{Z}}\circ (\mathbf{W}\oplus {\alpha}_{\mathbf{X}, \mathbf{Y}, \mathbf{Z}})$ has $2$-cell component given by the following pasting diagram in $\mathcal{A}\otimes \mathcal{A}$.

	$$\begin{tikzcd}[font=\fontsize{9}{6}, column sep = 15]
		S(W \oplus X\oplus Y \oplus Z)
		\arrow[rr, "\chi_{W{,}X\oplus Y\oplus Z}"]
		\arrow[dd, "\chi_{W\oplus X\oplus Y{,}Z}"']
		&{}\arrow[dd, Leftarrow, shorten = 15,shift right = 8, "\omega_{W{,}X\oplus Y{,}Z}"]
		&SW\oplus S(X\oplus Y\oplus Z)
		\arrow[rr, "w\oplus 1"]
		\arrow[dd, "1\oplus \chi_{X\oplus Y{,}Z}"description]
		&{}\arrow[dd, Leftarrow, shorten = 15,shift right = 5, "w_{\chi_{X\oplus Y{,}Z}}"]
		&W\oplus S(X\oplus Y\oplus Z)
		\arrow[dd, "1\oplus \chi_{X\oplus Y{,}Z}"description]
		\arrow[dd, Leftarrow, shorten = 15, shift left = 15, "1\oplus \omega_{X{,}Y{,}Z}"]
		\arrow[rrdd, bend left = 30, "1\oplus \chi_{X{,}Y\oplus Z}"]
		\\
		\\
		S(W\oplus X\oplus Y)\oplus SZ 
		\arrow[dd, "\chi_{W\oplus X{,}Y}\oplus 1"']
		\arrow[rr, "\chi_{W{,}X\oplus Y}\oplus 1"']
		&{}\arrow[dd, Leftarrow, shift right = 8,shorten = 15, "\omega_{W{,}X{,}Y}\oplus 1"]
		&SW\oplus S(X\oplus Y)\oplus SZ 
		\arrow[rr, "w\oplus 1\oplus 1"'] 
		\arrow[dd, "1\oplus \chi_{X{,}Y}\oplus 1"description]
		&{}\arrow[dd, Leftarrow, shorten = 15, shift right = 5, "w_{\chi_{X{,}Y}}\oplus 1"]
		& W \oplus S(X\oplus Y) \oplus SZ
		\arrow[dd, "1\oplus \chi_{X{,}Y}\oplus 1"description]
		&&W\oplus SX\oplus S(Y\oplus Z)
		\arrow[lldd, "1\oplus 1\oplus \chi_{Y{,}Z}"description]
		\arrow[dd, "1\oplus x\oplus 1"]
		\\
		\\
		S(W\oplus X)\oplus SY \oplus SZ
		\arrow[rr, "\chi_{W{,}X}\oplus 1\oplus 1"']
		&{}&SW\oplus SX\oplus SY \oplus SZ\arrow[rr, "w\oplus 1\oplus 1"']
		&{}& W \oplus SX\oplus SY \oplus SZ
		\arrow[dd, "1\oplus x\oplus  1\oplus 1"']
		\arrow[dd, Leftarrow, shorten = 15, shift left = 10, "1\oplus x_{\chi_{Y{,}Z}}"]
		&&W\oplus X \oplus S(Y\oplus Z)
		\arrow[lldd, bend left = 30,"1\oplus 1\oplus \chi_{Y{,}Z}"] 
		\\
		\\
		W\oplus X\oplus  Y\oplus Z
		&& W \oplus X\oplus Y\oplus SZ
		\arrow[ll, "1 \oplus 1\oplus 1\oplus z"]
		&& W \oplus X\oplus SY\oplus SZ
		\arrow[ll, "1\oplus 1\oplus  y\oplus 1"] 
	\end{tikzcd}$$
	
	\noindent Apply the interchange law for the morphism $w: SW \to W$ on the $2$-cell $\omega_{X, Y, Z}$ in the $\mathbf{Gray}$-monoid $\A$, followed by the pentagon law for the opmonoidal $2$-functor $S: \A \to \A$ to arrive at the pasting diagram below. Observe that this is the required composite ${\alpha}_{\mathbf{W}, \mathbf{X}, \mathbf{Y}\oplus \mathbf{Z}} \circ {\alpha}_{\mathbf{W}\oplus \mathbf{X}, \mathbf{Y}, \mathbf{Z}}$, completing the proof.

	$$\begin{tikzcd}[font=\fontsize{9}{6}, column sep = 15]
		S(W \oplus X\oplus Y \oplus Z)
		\arrow[dddd, Leftarrow, shorten = 20,shift left = 25, red, "\omega_{W\oplus X{,}Y{,}Z}"description]
		\arrow[rrdd, blue, bend left = 20,"\chi_{W\oplus X{,}Y\oplus Z}"description]
		\arrow[rr, "\chi_{W{,}X\oplus Y\oplus Z}"]
		\arrow[dd, "\chi_{W\oplus X\oplus Y{,}Z}"']
		&{}
		&SW\oplus S(X\oplus Y\oplus Z)
		\arrow[dd, Leftarrow, shorten =15,red, "\omega_{W{,}X{,}Y\oplus Z}"]
		\arrow[rrdd, bend left = 20,blue,"1\oplus \chi_{X{,}Y\oplus Z}"description]
		\arrow[rr, "w\oplus 1"]
		&{}
		&W\oplus S(X\oplus Y\oplus Z)
		\arrow[dd, Leftarrow, red,shorten = 15, "w_{\chi_{X{,}Y\oplus Z}}"]
		\arrow[rrdd, bend left = 30, "1\oplus \chi_{X{,}Y\oplus Z}"]
		\\
		\\
		S(W\oplus X\oplus Y)\oplus SZ 
		\arrow[dd, "\chi_{W\oplus X{,}Y}\oplus 1"']
		&{}
		&\color{blue}S(W\oplus X)\oplus S(Y\oplus Z)
		\arrow[rr,blue, "\chi_{W{,}X}\oplus 1"]
		\arrow[lldd,blue, bend left = 20,"1\oplus \chi_{Y{,}Z}"description]
		\arrow[dd, red, Leftarrow, shorten = 15, "{(\chi_{W{,}X})}_{(\chi_{Y{,}Z})}"]
		&{}
		&\color{blue}SW\oplus SX \oplus S(Y\oplus Z)
		\arrow[lldd,bend left = 20, blue,"1\oplus 1\oplus \chi_{Y{,}Z}"description]
		\arrow[dd, Leftarrow,shorten = 5, red, "w_{1\oplus \chi_{Y{,}Z}}"]
		\arrow[rr, blue,"w\oplus 1\oplus 1"]
		&&W\oplus SX\oplus S(Y\oplus Z)
		\arrow[dd, "1\oplus x\oplus 1"]
		\arrow[lldd,"1\oplus 1\oplus \chi_{Y{,}Z}"]
		\\
		\\
		S(W\oplus X)\oplus SY \oplus SZ
		\arrow[rr, "\chi_{W{,}X}\oplus 1\oplus 1"']
		&{}&SW\oplus SX\oplus SY \oplus SZ\arrow[rr, "w\oplus 1\oplus 1"']
		&{}& W \oplus SX\oplus SY \oplus SZ
		\arrow[dd, "1\oplus x\oplus  1\oplus 1"']
		\arrow[dd, Leftarrow, shorten = 15, shift left = 10, "1\oplus x_{\chi_{Y{,}Z}}"]
		&&W\oplus X \oplus S(Y\oplus Z)
		\arrow[lldd, bend left = 30,"1\oplus 1\oplus \chi_{Y{,}Z}"] 
		\\
		\\
		W\oplus X\oplus  Y\oplus Z
		&& W \oplus X\oplus Y\oplus SZ
		\arrow[ll, "1 \oplus 1\oplus 1\oplus z"]
		&& W \oplus X\oplus SY\oplus SZ
		\arrow[ll, "1\oplus 1\oplus  y\oplus 1"] 
	\end{tikzcd}$$

\end{proof}

\subsection{Proposition \ref{Symmetriser's component 2-natural transformation}}\label{Appendix Proposition symmetriser's component 2-nat trans}

\begin{proof}
	The required $2$-natural transformation is determined by the universal property of the Eilenberg-Moore object ${\left(\mathcal{B}\otimes \mathcal{A}\right)}^{(T\otimes S)}$ by the tight pseudomorphism for $\left(T\otimes S\right) \circ -$ whose $2$-cell component is given by the interchanger of $\varepsilon^{S}$ and $\varepsilon^{T}$. This has component on the pseudoalgebra $\left(\mathbf{X}, \mathbf{Y}\right) \in \mathcal{A}^{S} \otimes \mathcal{B}^{T}$ is the tight pseudomorphism given by the interchanger $\left(x, y\right)$ in the $\mathbf{Gray}$ tensor product $\mathcal{B}\otimes \mathcal{A}$. It suffices to check that this is well-defined as a pseudomorphism. For the unit law, begin with the pasting diagram displayed below.

	$$\begin{tikzcd}[font = \fontsize{9}{6}, column sep = 16]
		&&(TY{,}X)
		\arrow[dd, Rightarrow, shorten = 15, "{(\eta_{Y})}_{(\eta_{X})}"]
		\arrow[rrdd, "(1{,}\eta_{X})"]
		\\
		\\
		\left(Y{,}X\right)
		\arrow[rruu, "(\eta_{Y}{,}1)"]
		\arrow[rr, "\left(1{,}\eta_{X}\right)"]
		\arrow[rrdd, bend right = 30, "1"']
		&{}\arrow[dd, shorten = 15, Rightarrow, "\left(1{,}x_{0}\right)"]
		& \left(Y{,}SX\right)
		\arrow[rr, "\left(\eta_{Y}{,}1\right)"]
		\arrow[dd, "\left(1{,}x\right)"]
		&{}\arrow[dd, Rightarrow, shorten = 15, "{(\eta_{Y})}_{(x)}"]
		& \left(TY{,}SX\right)
		\arrow[dd, "\left(1{,}x\right)"]
		\arrow[rrdd, "(y{,}1)"]
		\arrow[dddd, Rightarrow, shorten = 30, shift left = 10, "{(y)}_{(x)}"]
		\\
		\\
		&{}& \left(Y{,}X\right)\arrow[rr, "\left(\eta_{Y}{,}1\right)"]
		\arrow[rrdd, bend right =30, "1"']
		&{}\arrow[dd, Rightarrow, shorten = 15, "\left(y_{0}{,}1\right)"]
		&\left(TY{,}X\right)
		\arrow[dd, "\left(y{,}1\right)"]
		&&(Y{,}SX)
		\arrow[lldd, "(1{,}x)"]
		\\
		\\
		&&&{}&\left(Y{,}X\right) &{}
	\end{tikzcd}$$
	
	\noindent Apply the interchange law for the morphism $x: SX \to X$ and the $2$-cell $y_{0}$ to arrive at the pasting displayed below left. Finally, apply the interchange law for the morphism $\eta_{X}: X \to SX$ on the $2$-cell $y_{0}$ to arrive at the pasting depicted below right and observe that this completes the proof of the unit law.

	$$\begin{tikzcd}[font = \fontsize{9}{6}, column sep = 15]
		&&(TY{,}X)
		\arrow[dd, Rightarrow, shorten = 15, "{(\eta_{Y})}_{(\eta_{X})}"]
		\arrow[rrdd, "(1{,}\eta_{X})"]
		\\
		\\
		\left(Y{,}X\right)
		\arrow[rruu, "(\eta_{Y}{,}1)"]
		\arrow[rr, "\left(1{,}\eta_{X}\right)"]
		\arrow[rrdd, bend right = 30, "1"']
		&{}\arrow[dd, shorten = 15, Rightarrow, "\left(1{,}x_{0}\right)"]
		& \left(Y{,}SX\right)
		\arrow[rrrr, bend right = 45,blue, "1"']
		\arrow[rr, "\left(\eta_{Y}{,}1\right)"]
		\arrow[dd, "\left(1{,}x\right)"]
		&{}
		& \left(TY{,}SX\right)
		\arrow[d, Rightarrow, shorten = 5, red, "(y_{0}{,}1)"']
		\arrow[rr, "(y{,}1)"]
		&&(Y{,}SX)
		\arrow[dd, "(1{,}x)"]
		\\
		&&&&{}
		\\
		&{}& \left(Y{,}X\right)
		\arrow[rrrr, equal, red, shorten = 60, shift left = 5]
		\arrow[rrrr, "1"']
		&{}
		&
		&&\left(Y{,}X\right)&{}
	\end{tikzcd}\begin{tikzcd}[font = \fontsize{9}{6}, column sep = 15]
		\left(Y{,}X\right)
		\arrow[rr, "\left(\eta_{Y}{,}1\right)"]
		\arrow[rrdd, bend right = 30, "1"']
		&{}\arrow[dd, shorten = 15, red,Rightarrow, "\left(y_{0}{,}1\right)"]
		& \left(TY{,}X\right)\arrow[rr, "\left(1{,}\eta_{X}\right)"]
		\arrow[dd, blue,"\left(y{,}1\right)"]
		&{}\arrow[dd, Rightarrow, red,shorten = 15, "{(y)}_{\left(\eta_{X}\right)}"]
		& \left(TY{,}SX\right)
		\arrow[dd, "\left(y{,}1\right)"]
		\\
		\\
		&{}& \left(Y{,}X\right)
		\arrow[rr, "\left(1{,}\eta_{X}\right)"]
		\arrow[rrdd, bend right =30, "1"']
		&{}\arrow[dd, Rightarrow, shorten = 15, "\left(1{,}x_{0}\right)"]
		&\left(Y{,}SX\right)
		\arrow[dd, "\left(1{,}x\right)"]
		\\
		\\
		&&&{}&\left(Y{,}X\right)
	\end{tikzcd}$$	
	
	\noindent For the multiplication law, begin with the pasting diagram displayed below.

	$$\begin{tikzcd}[font=\fontsize{9}{6}]
		&&\left(T^{2}Y{,}S^{2}X\right)
		\arrow[dddd, Rightarrow, shorten = 25,shift right = 15, "{(\mu_{Y})}_{(\mu_{X})}^{-1}"']
		\arrow[lldd, "(\mu_{Y}{,}1)"']
		\arrow[rr, "\left(1{,}Sx\right)"]
		\arrow[dd, "\left(1{,}\mu_{X}\right)"description]
		&{}\arrow[dd, Rightarrow,shift right = 7,  shorten = 15, "\left(1{,}x_{2}\right)"]
		&\left(T^2Y{,}SX\right)
		\arrow[rr, "\left(Ty{,}1\right)"]
		\arrow[dd, "\left(1{,}x\right)"description]
		&{}\arrow[dd, Rightarrow, shorten = 15, "{(Ty)}_{(x)}"]
		&\left(TY{,}SX\right)
		\arrow[dd, "\left(1{,}x\right)"description]
		\\
		\\
		(TY{,}S^2X)
		\arrow[rrdd, "(1{,}\mu_{X})"']
		&&\left(T^2Y{,}SX\right)
		\arrow[dd, "\left(\mu_{Y}{,}1\right)"description]
		\arrow[rr, "\left(1{,}x\right)"]
		&{}\arrow[dd, Rightarrow, shorten = 15, "{(\mu_{Y})}_{\left(x\right)}^{-1}"]
		&\left(T^2Y{,}SX\right)
		\arrow[dd, "\left(\mu_{Y}{,}1\right)" description]
		\arrow[rr, "\left(Ty{,}1\right)"]
		&{}\arrow[dd, Rightarrow, shorten = 15, "\left(y_{2}{,}1\right)"]
		&\left(TY{,}X\right)
		\arrow[dd, "\left(y{,}1\right)"description]
		\\
		\\
		&&\left(TY{,}SX\right)
		\arrow[rrdd, "(y{,}1)"']
		\arrow[rr, "\left(1{,}x\right)"']
		&{}&\left(TY{,}X\right)
		\arrow[dd, Rightarrow,shorten = 15, "{{(y)}_{(x)}}^{-1}"]
		\arrow[rr, "\left(y{,}1\right)"']
		&{}&\left(Y{,}X\right)
		\\
		\\
		&&&&(Y{,}SX)\arrow[rruu, "(1{,}x)"']
	\end{tikzcd}$$
	
	\noindent Apply the interchange law for the morphism $\mu_{Y}: T^2Y \to TY$ on the $2$-cell $x_{2}$, followed by the interchange law for the morphism $x: SX \to X$ on the $2$-cell $y_{2}$  to arrive at the pasting diagram displayed below.


	$$\begin{tikzcd}[font=\fontsize{9}{6}]
		&&\left(T^{2}Y{,}S^{2}X\right)
		\arrow[dd, Rightarrow, shorten = 15,red,"{(\mu_{Y})}_{(Sx)}"]
		\arrow[lldd, "(\mu_{Y}{,}1)"']
		\arrow[rr, "\left(1{,}Sx\right)"]
		&{}
		&\left(T^2Y{,}SX\right)
		\arrow[lldd,blue,"(\mu_{Y}{,}1)"description]
		\arrow[rr, "\left(Ty{,}1\right)"]
		\arrow[dd, Rightarrow, shorten = 15, red,"\left(y_{2}{,}1\right)"]
		&{}
		&\left(TY{,}SX\right)
		\arrow[lldd,blue, "(y{,}1)"description]
		\arrow[dddd, Rightarrow, shorten = 35, shift right = 15,red, "{(y)}_{(x)}"]
		\arrow[dd, "\left(1{,}x\right)"description]
		\\
		\\
		(TY{,}S^2X)
		\arrow[rr,blue,"(1{,}Sx)"']
		\arrow[rrdd, "(1{,}\mu_{X})"']
		&&\color{blue}\left(TY{,}SX\right)
		\arrow[rr, blue, "(y{,}1)"]
		\arrow[dd, Rightarrow, red, shorten = 15, "(1{,}x_{2})"]
		\arrow[rrdd,blue, "(1{,}x)"description]
		&{}
		&\color{blue}\left(Y{,}SX\right)
		\arrow[dd, Rightarrow, shorten = 15, red, "{(y)}_{(x)}"]
		\arrow[rrdd, blue,"(1{,}x)"description]
		&{}
		&\left(TY{,}X\right)
		\arrow[dd, "\left(y{,}1\right)"description]
		\\
		\\
		&&\left(TY{,}SX\right)
		\arrow[rrdd, "(y{,}1)"']
		\arrow[rr, "\left(1{,}x\right)"']
		&{}&\left(TY{,}X\right)
		\arrow[dd, Rightarrow,shorten = 15, "{{(y)}_{(x)}}^{-1}"]
		\arrow[rr, "\left(y{,}1\right)"']
		&{}&\left(Y{,}X\right)
		\\
		\\
		&&&&(Y{,}SX)\arrow[rruu, "(1{,}x)"']
	\end{tikzcd}$$
	
	\noindent Apply the interchange law for the morphism $y: TY \to Y$ on the $2$-cell $x_{2}$ followed by the interchange law for the morphism $Sx: S^2X \to SX$ on the $2$-cell $y_{2}$ to arrive at the pasting displayed below. Finally, observe that this completes the proof of the multiplication law.


	$$\begin{tikzcd}[font=\fontsize{9}{6}]
		&&(T^2Y{,}SX)
		\arrow[dd, Rightarrow, shorten = 15, red, "(Ty{,}Sx)"]
		\arrow[rrdd, "(Ty{,}1)"]
		\\
		\\
		\left(T^{2}Y{,}S^{2}X\right)
		\arrow[rruu, "(1{,}Sx)"]
		\arrow[rr,blue, "\left(Ty{,}1\right)"]
		\arrow[dd, "\left(\mu_{Y}{,}1\right)"']
		&{}\arrow[dd, Rightarrow,shift right = 7, red, shorten = 15, "\left(y_{2}{,}1\right)"]
		&\color{blue}\left(TY{,}S^{2}X\right)
		\arrow[rr,blue, "\left(1{,}Sx\right)"]
		\arrow[dd, blue,"\left(y{,}1\right)"description]
		&{}\arrow[dd, Rightarrow, shorten = 15, red,"{(y)}_{(Sx)}"]
		&\left(TY{,}SX\right)
		\arrow[dd, "\left(y{,}1\right)"description]
		\arrow[rrdd, "(1{,}x)"]
		\arrow[dddd, Rightarrow, shorten = 35, shift left = 10, "(y{,}x)"]
		\\
		\\
		\left(TY{,}S^{2}X\right)
		\arrow[dd, "\left(1{,}\mu_{X}\right)"']
		\arrow[rr, blue, "\left(y{,}1\right)"]
		&{}\arrow[dd, red, Rightarrow, shorten = 15, "{(y)}_{\left(\mu_{X}\right)}"]
		&\color{blue}\left(Y{,}S^{2}X\right)
		\arrow[dd, blue, "\left(1{,}\mu_{X}\right)" description]
		\arrow[rr, blue, "\left(1{,}Sx\right)"]
		&{}\arrow[dd, red, Rightarrow, shorten = 15, "\left(1{,}x_{2}\right)"]
		&\left(Y{,}SX\right)
		\arrow[dd, "\left(1{,}x\right)"description]
		&&(TY{,}X)\arrow[lldd, "(y{,}1)"]
		\\
		\\
		\left(TY{,}SX\right)
		\arrow[rr, "\left(y{,}1\right)"']
		&{}&\left(Y{,}SX\right)
		\arrow[rr, "\left(1{,}x\right)"']
		&{}&\left(Y{,}X\right)
	\end{tikzcd}$$

\end{proof}


\begin{thebibliography}{9}
	\bibitem{BDSV Extended 3D bordism theory of modular objects} Bartlett, B., Douglas, C., Schommer-Pries, C., Vicary, J., Extended 3-dimensional bordism as the theory of modular objects, available at \url{https://arxiv.org/abs/1411.0945}
	\bibitem{BDSV Modular categories as representations of the 3-dimensional bordism 2-category} Bartlett, B., Douglas, C., Schommer-Pries, C., Vicary, J., Modular categories as representations of the 3-dimensional bordism 2-category, available at \url{https://arxiv.org/abs/1509.06811}
	\bibitem{Bourke Lobbia Skew Approach} Bourke, J., Lobbia, G., A skew approach for enrichment over $\mathbf{Gray}$-categories, \textit{Advances in Mathematics}, Vol. 434, (2023)
	\bibitem{Bruguieres Lack Virelizier} Bruguières, A., Lack, S., Virelizier, A., Hopf monads on monoidal categories, \textit{Advances in Mathematics}, Vol. 227, Issue 2, pp. 745-800, (2011)
	\bibitem{From Thin Concurrent Games to Generalized Species of Structures} Clairambault, P. Paquet, H., Olimpieri, F. From Thin Concurrent Games to Generalized Species of Structures, 38th Annual ACM/IEEE Symposium on Logic in Computer Science, (2023) 
	\bibitem{Crans Generalised Centers} Crans, S., Generalized centers of braided and
	sylleptic monoidal 2-categories, \textit{Advances in Mathematics}, 136, pp 183-223, (1998)
	\bibitem{crans tensor of gray categories} Crans, S., Tensor products of Gray-categories, \textit{Theory and Applications of Categories}, Vol. 5, No. 2, , pp. 12–69. (1999)
	\bibitem{monoidal bicategories and hopf algebroids} Day, B., Street, R., Monoidal Bicategories and Hopf Algebroids, \textit{Advances in Mathematics}, Vol. 129, Issue 1, pp. 99-157, (1997)
	\bibitem{Gray categories model algebraic tricategories} Ferrer, G., Gray-categories model algebraic tricategories, \textit{Theory and Applications of Categories}, Vol. 38, No. 29, pp. 1136–1155, (2022)
	\bibitem{Fiore Gambina Hyland Winskel generalised species} Fiore, M., Gambino, N., Hyland, M., Winskel, G., The cartesian closed bicategory of generalised species of structures, \textit{Journal of the London Mathematical Society}, Vol. 77, Issue 1, pp 203–220, (2008), available at \url{https://doi.org/10.1112/jlms/jdm096}
	\bibitem{Fiore Gambino Hyland Monoidal Bicategories differential linear logic and analytic functors} Fiore, M., Gambino, N., Hyland, M., Monoidal Bicategories, Differential Linear Logic, and Analytic Functors, \textit{In preparation}
	\bibitem{Galal Bicategorical Finite Nondeterminism} Galal, Z., A Bicategorical Model for Finite Nondeterminism, \textit{Proceedings of the 6th International Conference on Formal Structures for Computation and Deduction}, (2021)
	\bibitem{Formal Theory of Pseudomonads} Gambino, N., Lobbia, G. On the formal theory of pseudomonads and pseudodistributive laws, \textit{Theory and Applications of Categories}, Vol. 37, No. 2, pp 14-56, (2021)
	\bibitem{Low dimensional structures formed by tricategories} Garner, R., Gurski, N., The low dimensional structures formed by tricategories, \textit{Mathematical Proceedings of the Cambridge Mathematical Society}, (2009).
	\bibitem{Garner Shulman Enriched Categories as a Free Cocompletion} Garner, R., Shulman, M., Enriched categories as a free cocompletion, \textit{Advances in Mathematics}, Vol. 289, pp 1-94, (2016)
	\bibitem{GPS tricategory} Gordon, R., Power, A.J., Street, R., \textit{Coherence For Tricategories}, American Mathematical Society, Vol. 81, (1995)
	\bibitem{Gurski PhD} Gurski, N. An algebraic theory of tricategories, Thesis (Ph. D.)--University of Chicago, Dept. of Mathematics, (2006), available at \url{https://ncatlab.org/nlab/files/Gurski-AlgebraicTricategories.pdf}
	\bibitem{Gurski biequivalence in tricategories} Gurski, N., Biequivalences in Tricategories, \textit{Theory and Applications of Categories}, Vol. 26, No. 14, pp. 349-384 (2012)
	\bibitem{Gurski Coherence in Three Dimensional Category Theory} Gurski, N., \textit{Coherence in Three-Dimensional Category Theory}, Cambridge Tracts in Mathematics, Vol. 201, (2013) 
	\bibitem{Gurski Monoidal Structure of Strictification} Gurski, N., The monoidal structure of strictification, \textit{Theory and Applications of Categories}, Vol. 28, No. 1, pp. 1–23, (2013)
	\bibitem{Gurski Osorno Infinite Loop Spaces and coherence for symmetric monoidal bicategories} Gurski, N., Osorno, A., Infinite loop spaces and coherence for symmetric monoidal bicategories, \textit{Advances in Mathematics}, Vol. 246, pp. 1–32, (2013)
	\bibitem{Gurski Johnson Osorno K Theory for 2-categories} Gurski, N., Johnson, N., Osorno, A., K-theory for 2-categories, \textit{Advances in Mathematics}, Vol. 322, pp. 378–472, (2017)
	\bibitem{Gurski Loop Spaces} Gurski, N., Loop spaces, and coherence for monoidal and braided monoidal bicategories, \textit{Advances in Mathematics}, Vol. 226, pp 4225–4265, (2011)
	\bibitem{Hansen Shulman Constructing Symmetric Monoidal Bicategories Functorially} Hansen, L., Shulman, M., Constructing symmetric monoidal bicategories functorially, available at \url{https://arxiv.org/abs/1910.09240}
	\bibitem{Hasegawa Lemay Star Autonomous Monads} Hasegawa, M., Lemay, J., Linear Distributivity With Negation, Star-Autonomy, and Hopf Monads, \textit{Theory and Applications of Categories}, Vol. 33, No. 37, pp 1145-1157, (2018)
	\bibitem{Traced Monads and Hopf Monads} Hasegawa, M., Lemay, J., Traced Monads and Hopf Monads, \textit{Compositionality}, Vol. 5, Issue 10, (2023)
	\bibitem{Hyland Power Pseudo commutative monads and pseudo closed categories} Hyland, M., Power, A. J., Pseudo-commutative Monads and Pseudo-closed 2-categories, \textit{Journal of Pure and Applied Algebra}, Vol. 175, pp. 141–185, (2002)
	\bibitem{Kelly Doctrinal Adjunction} Kelly, G. M., Doctrinal Adjunction, \textit{Lecture Notes in Mathematics}, (1974), Vol. 420, pp. 257-280, doi:10.1007/BFb0063105
	\bibitem{Kelly Basic Concepts of Enriched Category Theory} Kelly, G. M., \textit{Basic concepts of enriched category theory}, Reprints in Theory and Applications of Categories, No. 10, (2005)
	\bibitem{Icons} Lack, S., Icons, \textit{Applied Categorical Structures}, Vol 18, pp 289–307, (2010)
	\bibitem{FTM2} Lack, S., Street, S., The Formal Theory of Monads II, \textit{Journal of Pure and Applied Algebra}, Vol. 175, pp 243-265, (2002)
	\bibitem{Lack Limits and Lax Morphisms} Lack, S. Limits for Lax Morphisms. \emph{Applied Categorical Structures}, Vol. 13, pp. 189–203 (2005). https://doi.org/10.1007/s10485-005-2958-5
	\bibitem{Marmolejo Pseudodistributive 1} Marmolejo, F., Distributive laws for pseudomonads I, \textit{Theory and Applications of Categories}, Vol. 5, No. 5, pp 81-147, (1999).
	\bibitem{McCrudden Opmonoidal Monads} McCrudden, P. Opmonoidal Monads, \textit{Theory and Applications of Categories}, Vol. 10, No.19, pp.469–485, (2002).
	\bibitem{Miranda strictifying operational coherences} Miranda, A., Strictifying operational coherences and weak functor classifiers in low dimensions, available at \url{https://arxiv.org/abs/2307.01498}.
	\bibitem{Miranda Enriched Kleisli objects for pseudomonads} Miranda, A., Enriched Kleisli objects for pseudomonads, available at \url{https://arxiv.org/abs/2311.15618}
	\bibitem{Miranda Extending Symmetric Monoidal Structures to Kleisli bicategories} Miranda, A., Extending Symmetric Monoidal Structures to Kleisli bicategories, \textit{In preparation}
	\bibitem{Miranda Tricategorical Limits and Colimits} Miranda, A., Tricategorical limits and colimits, \textit{In preparation}.
	\bibitem{Miranda PhD} Miranda, A., Topics in Low Dimensional Higher Category Theory, Ph. D. Thesis, Macquarie University, (2023)
	\bibitem{Moerdijk Monads on Tensor Categories} Moerdijk, I., Monads on Tensor Categories, \textit{Journal of Pure and Applied Algebra}, Vol. 168, Issues 2–3, pp. 189-208, (2002)
	\bibitem{Pacquet Saville Strong Pseudomonads and Premonoidal Bicategories} Paquet, H., Saville, P., Strong Pseudomonads and Premonoidal Bicategories, available at \url{https://arxiv.org/abs/2304.11014}
	\bibitem{Pastro Street closed categories star autonomy and monoidal comonads} Pastro, C., Street, R., Journal of Algebra,	Vol. 321, Issue 11, pp. 3494-3520, (2009)
	\bibitem{Schommer Piers Classification of 2D ETFTs} Schommer-Piers, C., The Classification of Two-Dimensional Extended Topological Field Theories (2014), available at \url{https://arxiv.org/abs/1112.1000#:~:text=Thereby%20we%20classify%20these%20types,a%20fixed%20commutative%20ground%20ring.} 
	\bibitem{Stell modelling term rewriting via sesquicategories} Stell. J., Modelling Term Rewriting Systems by Sesqui-Categories, \textit{Proc. Catégories, Algèbres, Esquisses et Néo-Esquisses} (1994)
	\bibitem{Formal Theory of Monads} Street R. The formal theory of monads, \textit{Journal of Pure and Applied Algebra}, Vol. 2, pp. 149-168, (1978)
	\bibitem{Zawadowski Monoidal Monads} Zawadowski, M., The Formal Theory of Monoidal Monads, \textit{Journal of Pure and Applied Algebra}, Vol. 216, pp. 1932-1942, (2012)
\end{thebibliography}
\end{document}